\documentclass[twocolumn]{autart}    

\usepackage{amsmath}
\usepackage{amssymb}
\usepackage{bm}
\usepackage{tikz}
\usetikzlibrary{matrix,arrows,decorations.pathmorphing}
\usepackage{url}

\usepackage{xcolor}

\newcommand{\R}{\mathbb{R}}
\newcommand{\N}{\mathbb{N}}

\newcommand{\E}{\mathbb{E}}

\newcommand{\Rn}{\mathbb{R}^{n}}
\newcommand{\Rpo}{\mathbb{R}^{p_1}}
\newcommand{\Rpt}{\mathbb{R}^{p_2}}
\newcommand{\Rnn}{\mathbb{R}^{n\times n}}

\newcommand{\Sn}{\mathcal{S}_n}
\newcommand{\Pn}{\mathcal{P}_{n}}
\newcommand{\Ppo}{\mathcal{P}_{p_1}}
\newcommand{\Ppt}{\mathcal{P}_{p_2}}
\newcommand{\PPn}{\mathcal{P}_{n}^{\succ}}
\newcommand{\PPpo}{\mathcal{P}_{p_1}^{\succ}}
\newcommand{\PPpt}{\mathcal{P}_{p_2}^{\succ}}
\newcommand{\Ppi}{\mathcal{P}_{p_1+p_2}}
\newcommand{\PPpi}{\mathcal{P}_{p_1+p_2}^{\succ}}

\newcommand{\xtrue}{x_{\mathrm{true}}}
\newcommand{\Cov}{\mathrm{Cov}}
\newcommand{\Pjoint}{P_{\mathrm{joint}}}

\newcommand{\sigmamax}{\sigma_{\mathrm{max}}}

\newcommand{\DomPhat}{\mathcal{D}_{\boldsymbol{\hat{P}}}}
\newcommand{\DomPhatMax}{\mathcal{D}_{\boldsymbol{\hat{P}}}^{\max}}
\newcommand{\Domcons}{\mathcal{D}_{\mathrm{cons}}}
\newcommand{\DomSDP}{\mathcal{D}_{\mathrm{SDP}}}
\newcommand{\DomF}{\mathcal{D}_{\bF}}
\newcommand{\Domalpha}{\mathcal{D}_{\alpha}}
\newcommand{\Domgamma}{\mathcal{D}_{\gamma}}

\newcommand{\bPhat}{\boldsymbol{\hat{P}}}
\newcommand{\bQ}{\boldsymbol{Q}}

\newcommand{\btau}{\boldsymbol{\tau}}

\newcommand{\bK}{\boldsymbol{K}}
\newcommand{\balpha}{\boldsymbol{\alpha}}
\newcommand{\bF}{\boldsymbol{F}}

\newcommand{\bSigma}{\boldsymbol{\Sigma}}
\newcommand{\bDelta}{\boldsymbol{\Delta}}

\newcommand{\CIa}{{\balpha}}
\newcommand{\CI}{{\balpha^\ast}}

\newcommand{\ro}{r_1^{\alpha^\ast}}
\newcommand{\rt}{r_2^{\alpha^\ast}}

\DeclareMathOperator{\Trace}{tr}

\DeclareMathOperator*{\argmin}{argmin}
\DeclareMathOperator{\rank}{rk}
\DeclareMathOperator{\adj}{adj}
\DeclareMathOperator{\vectorize}{vec}
\DeclareMathOperator{\dom}{dom}

\newcommand{\startM}{\left(\begin{smallmatrix}}
\newcommand{\stopM}{\end{smallmatrix}\right)}

\begin{document}

\begin{frontmatter}
\runtitle{Optimal conservative linear fusion}  

\title{Conservative fusion of unbiased partial state estimates:\\CI is optimal\thanksref{footnoteinfo}} 
\thanks[footnoteinfo]{This paper was not presented at any IFAC meeting. Corresponding author J.~Trumpf.}

\author[ANU]{Jochen Trumpf}\ead{Jochen.Trumpf@anu.edu.au},    
\author[UoM]{Behzad Zamani}\ead{behzad.zamani@unimelb.edu.au},               
\author[UoM]{Chris Manzie}\ead{manziec@unimelb.edu.au}  

\address[ANU]{Australian National University, Acton ACT 2601, Australia}  
\address[UoM]{University of Melbourne, VIC 3010, Australia}             


\begin{abstract}                          
We show that Covariance Intersection (CI) is optimal amongst all conservative unbiased linear fusion rules also in the general case of information fusion of two unbiased partial state estimates, significantly generalizing the known optimality result for fusion of full state estimates.
In fact, we prove the much stronger result that three different optimization problems are equivalent, namely the abstract optimal conservative unbiased linear information fusion problem with respect to a strictly isotone cost function, the scalar CI problem, and a simple semi-definite program.
We provide a general solvability condition for these problems as well as equations characterizing the optimal solutions for the matrix determinant and matrix trace cost functions.
\end{abstract}

\end{frontmatter}

\section{Introduction}

State estimation algorithms based on information fusion, sometimes also called data-fusion, provide an alternative to filtering algorithms such as the Kalman filter and its derivatives in situations where the correlation between the errors of a current state estimate and a second state estimate used to update this current estimate is not known or would be too difficult or too expensive to track.
This is particularly relevant in distributed estimation where tracking of cross-covariances between the states of different nodes in a network requires either centralised computation or reliable all-to-all inter-node communication~\cite{roumeliotis}.
Because of this, information fusion algorithms such as Covariance Intersection (CI)~\cite{Uhlmann1996} and its derivatives have been applied to distributed state estimation scenarios for many years, see the recent survey \cite{CI_survey} and the references therein.
A generalization of CI that allows fusion of two partial state estimates has been known for at least as far back as the work of Arambel et al. \cite{arambel2001} who applied this generalization to the distributed tracking of spacecraft.

It was shown by Reinhardt et al. in 2015 \cite{reinhardt2015} that in the case of fusion of two unbiased \emph{full} state estimates, CI is optimal amongst all conservative unbiased linear fusion rules where optimality means smallest size of the fused error covariance as measured by a strictly isotone cost function such as the matrix determinant or matrix trace.

In this paper, we generalize this optimality result to the case of fusion of two unbiased \emph{partial} state estimates.
In this case, the pre-image in state space of the (measurement) error covariance of a partial state estimate is typically singular, introducing significant mathematical challenges that mean we need to revisit the entire chain of proofs leading to the final optimality result.
The (somewhat lengthy) proofs unveil an intriguing interplay between statistics (Gau{\ss}-Markov theorem), matrix analysis (Schur complement, L{\"o}wner order), robust control (Petersen's lemma), geometry (ellipsoidal calculus) and optimization (semi-definite programming).
The main new contributions of the paper are as follows:
\vspace*{-2mm}
\begin{itemize}
\item The concept of a domain of conservativeness that is closed under diagonal enlargement (Definition~\ref{def:conservative}) and the associated disambiguation between the use of true error covariances and their estimates in the derivation of the inequalities that characterize conservativeness (see Proposition~\ref{prop:conservative} and the preceding discussion)
\item A stand-alone necessary and sufficient characterization of conservative fused error covariance estimates in three different forms (Theorem~\ref{th:conservative_3})
\item The generalization of the Bar-Schalom/Campo formula to the singular case (Theorem~\ref{th:Kopt})
\item The proper mathematical treatment of the parameterization of the Kahan-Uhlmann family, including in the singular case (see Sections~\ref{sec:KU} and \ref{sec:CI}: Theorems~\ref{th:kahan_parameter} and \ref{th:CI_implies_Kahan} have no equivalent in the existing literature)
\item A general sufficient condition for solvability of CI (Theorem~\ref{th:CI}) and a corrected characterizing equation for the optimal solutions with respect to the determinant cost (Theorem~\ref{th:CI_det})
\end{itemize}

The remainder of this paper is organized as follows.
In Section~\ref{sec:linear_fusion} we formally introduce the linear fusion problem.
Section~\ref{sec:conservative} introduces and formalizes the concept of conservative error covariance estimates and contains our first main result, a necessary and sufficient characterization of conservative fused error covariance estimates (Theorem~\ref{th:conservative_3}).
Section~\ref{sec:optimal} introduces and formalizes the optimal conservative unbiased linear fusion problem~\eqref{def:OF} (Definition~\ref{def:optimal_fusion}) and provides our generalization of the celebrated Bar-Shalom/Campo formula~\cite{Bar-Schalom-Campo} to the case of partial state estimates (Theorem~\ref{th:Kopt}).
Section~\ref{sec:ell} explores the relationship between linear fusion and ellipsoidal calculus and provides a geometric necessary condition for conservativeness of the fused error covariance estimate (Theorem~\ref{th:conservative_fusion}).
Section~\ref{sec:KU} introduces what we call the Kahan-Uhlmann family of linear fusion methods (Definition~\ref{def:kahan_family}) and contains results showing that problem~\eqref{def:OF} can be reduced to a search over this family (Theorem~\ref{th:Kahan_implies_opt} and Theorem~\ref{th:opt_implies_Kahan}).
Section~\ref{sec:CI} demonstrates that CI (Definition~\ref{def:covariance_intersection}) is optimal (Corollary~\ref{cor:reduction_to_CI}).
Section~\ref{sec:SDP} shows that a relatively simple semi-definite program~\eqref{def:SDP} (Definition~\ref{def:semidefinite_program}) is equivalent to both problem~\eqref{def:OF} (Theorem~\ref{th:reduction_to_SDP}) and \eqref{def:CI} (Theorem~\ref{th:SDP_equiv_CI}).
Section~\ref{sec:CI_solution} provides a general sufficient condition for when problem~\eqref{def:CI} has a solution (Theorem~\ref{th:CI}) as well as characterizing equations for the optimal solution for the matrix determinant (Theorem~\ref{th:CI_det}) and matrix trace (Theorem~\ref{th:CI_tr}) cost functions.
Section~\ref{sec:conclusion} concludes the paper.
The three appendices
contain detailed proofs of all lemmas (Appendix~\ref{sec:Lemmas}), auxiliary results (Appendix~\ref{app:proofs}) and main results (Appendix~\ref{sec:proofs_main}).

\section{Linear fusion of unbiased estimators}\label{sec:linear_fusion}

The general setup for linear information fusion of partial state estimates in $\Rn$ is as follows.
We are trying to estimate an \emph{unknown true value} $\xtrue\in\Rn$.
We have two separate estimators $\hat{x}_1$ for $H_1\xtrue$ and $\hat{x}_2$ for $H_2\xtrue$, where $H_1\in\R^{p_1\times n}$, $p_1\leq n$ and $H_2\in\R^{p_2\times n}$, $p_2\leq n$ are \emph{full row rank} matrices.
Note that the assumption of full row rank is made for mathematical convenience.
A change of basis followed by a projection onto the non-zero component otherwise reduces the problem to the full row rank case.
We model $\hat{x}_1$ (resp. $\hat{x}_2$) as random variables with values in $\R^{p_1}$ (resp. $\R^{p_2})$.
We now seek to combine these two estimators into a \emph{fused estimator} $\hat{x}$, a random variable with values in $\Rn$, in a way that guarantees a number of conditions that we will spell out progressively in what follows.
For this to make sense, the random variables $\hat{x}_1$, $\hat{x}_2$ and $\hat{x}$ are modelled on the same underlying probability space but we will have little occasion to spell this out explicitly.
We will also tacitly assume that all our random variables have finite first order moments (mean or expected value) and second order moments (covariance).

We consider \emph{linear fused estimators} $\hat{x}$ of the form
\begin{equation}\label{eq:linear_fusion}
 \hat{x}=K_1\hat{x}_1+K_2\hat{x}_2,
\end{equation}
where $K_1\in\R^{n\times p_1}$ and $K_2\in\R^{n\times p_2}$ are constant matrix coefficients.
It follows immediately that
\begin{equation}\label{eq:mean_linear_fusion}
 \E[\hat{x}]=K_1\E[\hat{x}_1]+K_2\E[\hat{x}_2].
\end{equation}
Of specific interest is the situation where both $\hat{x}_1$ and $\hat{x}_2$ are \emph{unbiased} estimators. It follows from equation~\eqref{eq:mean_linear_fusion} that then $\E[\hat{x}]=(K_1H_1+K_2H_2)\xtrue$, and hence $\hat{x}$ is again unbiased if
and only if $(K_1H_1+K_2H_2)\xtrue=\xtrue$.
It follows that $\hat{x}$ is unbiased for all possible values of the unknown true value $\xtrue\in\Rn$ if and only if
\begin{equation}\label{eq:linear_unbiased_fusion}
 K_1H_1+K_2H_2=I.
\end{equation}
Since we are mostly interested in unbiased estimators, we will make the additional assumption that $H_1$ and $H_2$ are such that an unbiased estimator of the form~\eqref{eq:linear_fusion} exists.
Examining condition~\eqref{eq:linear_unbiased_fusion}, this is equivalent to the existence of a left inverse
$\begin{pmatrix} K_1 & K_2 \end{pmatrix}$ of the block matrix
\begin{equation}\label{eq:H}
 H := \begin{pmatrix} H_1 \\ H_2 \end{pmatrix}
\end{equation}
and hence to $\rank(H)=n$.

For ease of reference, we collect the main assumptions into the following list.
\begin{equation}\label{eq:A1}\tag{A1}
\begin{split}
 &\rank(H_1)=p_1,\ \rank(H_2)=p_2,\ \rank(H)=n, \\
 &\hat{x}_1\text{ resp. }\hat{x}_2\text{ are unbiased estimators of } \\
 &H_1\xtrue\text{ resp. }H_2\xtrue.
\end{split}
\end{equation}

We now introduce the \emph{error estimators} $\hat{e}_1:=\hat{x}_1-H_1\xtrue$, $\hat{e}_2:=\hat{x}_2-H_2\xtrue$ and $\hat{e}:=\hat{x}-\xtrue$.
If $\hat{x}_1$ and $\hat{x}_2$ are unbiased, it follows immediately that
\begin{align*}
 P_i &:= \Cov(\hat{e}_i)=\E[(\hat{x}_i-H_i\xtrue)(\hat{x}_i-H_i\xtrue)^\top] \\
 &= E[(\hat{x}_i-\E[\hat{x}_i])(\hat{x}_i-\E[\hat{x}_i])^\top] \\
 &= \Cov(\hat{x}_i),
\end{align*}
where the matrices $P_i\in\R^{p_i\times p_i}$ are symmetric positive semi-definite for $i\in\{1,2\}$.
In general, the two unbiased estimators $\hat{x}_1$ and $\hat{x}_2$ are \emph{correlated} and hence their error estimators $\hat{e}_1$ and $\hat{e}_2$ are correlated in the same way.
Specifically,
\begin{equation}\label{eq:joint_variance}\begin{split}
 \Pjoint &:= \Cov(\hat{e}_1,\hat{e}_2)=\Cov(\hat{x}_1,\hat{x}_2) \\
 &= \begin{pmatrix} P_1 & P_{12} \\ P_{12}^\top & P_2 \end{pmatrix},
\end{split}\end{equation}
where $P_{12}\in\R^{p_1\times p_2}$ is such that $\Pjoint$ is symmetric positive semi-definite.
It follows by direct calculation that condition~\eqref{eq:linear_unbiased_fusion} then implies that
\begin{equation}\label{eq:fusion_variance}\begin{split}
 P &:= \Cov(\hat{e})=\Cov(\hat{x}) \\
 &= \begin{pmatrix} K_1 & K_2 \end{pmatrix}
    \begin{pmatrix} P_1 & P_{12} \\ P_{12}^\top & P_2 \end{pmatrix}
    \begin{pmatrix} K_1^\top \\ K_2^\top \end{pmatrix}.
\end{split}\end{equation}

\section{Conservative error covariance estimates}\label{sec:conservative}

In the following, denote
\begin{align*}
 \Sn &:=\{P\in\R^{n\times n}\,|\,P^\top=P\}, \\
 \Pn &:=\{P\in\Sn\,|\,P\succeq 0\}
 \quad\text{and}\\
 \PPn &:=\{P\in\Pn\,|\,P\succ 0\},
\end{align*}
where $\preceq$ denotes the L{\"o}wner order.

In applications of information fusion rules, it is of paramount importance to have useful estimates of the error covariance $P$ of the fused estimator $\hat{x}$.
We would hence like to augment the fused estimator~\eqref{eq:linear_fusion} with property~\eqref{eq:linear_unbiased_fusion} with a way to combine the values of the unbiased estimators $\hat{x}_1$ and $\hat{x}_2$ with estimates $\hat{P}_1$ of $P_1$ and $\hat{P}_2$ of $P_2$ into an estimate $\hat{P}$ of $P$ \emph{without knowledge of the unknown cross-covariance $P_{12}$}.
These error covariance estimates are deterministic estimates and \emph{not} modelled as stochastic estimators (random variables) themselves.
In fact, what we are seeking is a map
\begin{align*}
 \bPhat\colon&\DomPhat\subseteq\DomPhatMax\longrightarrow\Pn, \\
 &(K_1,K_2,\hat{x}_1,\hat{x}_2,\hat{P}_1,\hat{P}_2)\mapsto\bPhat(K_1,K_2,\hat{x}_1,\hat{x}_2,\hat{P}_1,\hat{P}_2)
\end{align*}
that is computable by some algorithm.
In keeping with property~\eqref{eq:linear_unbiased_fusion}, the maximal domain $\DomPhatMax$ of $\bPhat$ is defined as
\begin{equation}\label{eq:DomPhatMax}\begin{split}
 \DomPhatMax := \{ &(K_1,K_2,\hat{x}_1,\hat{x}_2,\hat{P}_1,\hat{P}_2)\in \\
 &\R^{n\times p_1}\times\R^{n\times p_2}\times\Rpo\times\Rpt\times\Ppo\times\Ppt \\
 &|\, K_1H_1+K_2H_2=I\}.
\end{split}\end{equation}
Note that assumption~\eqref{eq:A1} implies that $\DomPhatMax\not=\emptyset$.
We will see later that it is sometimes convenient to place additional restrictions on the domain of $\bPhat$, for example that $\hat{P}_1$ and $\hat{P}_2$ are invertible and hence positive definite, or that there is an equation relating $(K_1,K_2)$ and $(\hat{P}_1,\hat{P}_2)$.
Where this is necessary we will define $\DomPhat$ as the appropriate subset.

In order to increase readability of what follows, we do not write out the arguments of $\bPhat$ explicitly unless this is necessary to avoid confusion.
In particular, the short-cut statement $\bPhat^{\mathrm{a}}\succeq\bPhat^{\mathrm{b}}$ is to be read as
$\bPhat^{\mathrm{a}}(K_1,K_2,\hat{x}_1,\hat{x}_2,\hat{P}_1,\hat{P}_2)\succeq\bPhat^{\mathrm{b}}(K_1,K_2,\hat{x}_1,\hat{x}_2,\hat{P}_1,\hat{P}_2)$ for all $(K_1,K_2,\hat{x}_1,\hat{x}_2,\hat{P}_1,\hat{P}_2)\in\DomPhat$.
Similarly, the statement $\bPhat\preceq\hat{P}_1$ is to be read as a uniform bound $\bPhat(K_1,K_2,\hat{x}_1,\hat{x}_2,\hat{P}_1,\hat{P}_2)\preceq\hat{P}_1$ and $\bPhat\in\PPn$ means the uniform statement $\bPhat(K_1,K_2,\hat{x}_1,\hat{x}_2,\hat{P}_1,\hat{P}_2)\in\PPn$.
The only situation where this interpretation does not work is when such statements appear inside conditional statements with two or more cases.
As an example, the statement
\begin{align*}
 \bPhat^{\mathrm{a}}\succeq
 \begin{cases}
  \bPhat^{\mathrm{b}}, & \bPhat^{\mathrm{a}}\succeq I \\
  \bPhat^{\mathrm{c}}, & \bPhat^{\mathrm{a}}\prec I \\
  \bPhat^{\mathrm{d}}, & \text{otherwise}
 \end{cases}
\end{align*}
is to be read pointwise in $(K_1,K_2,\hat{x}_1,\hat{x}_2,\hat{P}_1,\hat{P}_2)\in\DomPhat$, i.e. the first  case holds for those $(K_1,K_2,\hat{x}_1,\hat{x}_2,\hat{P}_1,\hat{P}_2)\in\DomPhat$ for which
$\bPhat^{\mathrm{a}}(K_1,K_2,\hat{x}_1,\hat{x}_2,\hat{P}_1,\hat{P}_2)\succeq I$, etc.
This interpretation should be self-evident from context.

To simplify notation, we will introduce other such maps wherever this helps to improve readability.
The use of a {\bf bold-face} font to indicate that an object is a map should help to avoid confusion.
Note that a map defined on $\DomPhat$ does not need to \emph{explicitly} depend on all the variables
$(K_1,K_2,\hat{x}_1,\hat{x}_2,\hat{P}_1,\hat{P}_2)$.
In fact, we will see later that some of our maps $\bPhat$ only explicitly depend on $K_1$, $K_2$, $\hat{P}_1$ and $\hat{P}_2$ or even only on $\hat{P}_1$ and $\hat{P}_2$, but this does not change the fact that they need to be treated as maps with the full signature as above in the logic of the proofs.

The proofs often proceed pointwise and case-by-case, i.e. by first selecting $(K_1,K_2,\hat{x}_1,\hat{x}_2,\hat{P}_1,\hat{P}_2)\in\DomPhat$ by some criterion and then proving the desired result for the value $\hat{P}:=\bPhat(K_1,K_2,\hat{x}_1,\hat{x}_2,\hat{P}_1,\hat{P}_2)$.
Where the different cases cover all of $\DomPhat$, this constitutes a valid proof for the entire map $\bPhat$. The proofs would be near unreadable if we carried the full signature of $\bPhat$ in every computation step and we hope that the distinction between bold-face $\bPhat$ for the map and normal-face $\hat{P}$ for its value at a given point in $\DomPhat$ alleviates this problem.

We have deliberately not chosen the alternative approach where $(K_1,K_2,\hat{x}_1,\hat{x}_2,\hat{P}_1,\hat{P}_2)\in\DomPhat$ (or the equivalent for other maps) is chosen at the beginning of the paper and all proofs are given for that chosen value.
In our opinion, this obscures the fact that we are after a general computational rule and its domain of definition plays into the theory in very subtle ways.
This will be most evident for the definition of conservative error covariance estimate (see Definition~\ref{def:conservative}) and its first characterization in Proposition~\ref{prop:conservative} but actually permeates the entire theory.
Changing $\DomPhat$ potentially changes the theory.

We can now introduce the main property that we will require of the error covariance estimate $\bPhat$.
The literature calls $\hat{P}$ a \emph{conservative} estimate for $P$ if $\hat{P}\succeq P$.
Given that $\hat{P}$ will be computed based on $\hat{P}_1$ and $\hat{P}_2$ that are themselves assumed to be conservative estimates of $P_1$ resp. $P_2$, we believe it to be more appropriate to speak of the computational rule $\bPhat$ as \emph{preserving conservativeness}.
Of course, the output $\hat{P}$ of such a rule is then a conservative estimate in the usual sense.
The following definition captures this idea.

\begin{defn}\label{def:conservative}
 Let $\Domcons\subseteq\Ppi$ be \emph{closed under diagonal enlargement}, i.e. let
 \begin{multline}\label{eq:diag_dom}
  \begin{pmatrix} P_1 & P_{12} \\ P_{12}^\top & P_2 \end{pmatrix} \in \Domcons
  \ \& \ P_1\preceq\hat{P}_1 \ \&\ P_2\preceq\hat{P}_2 \\
  \Rightarrow
  \begin{pmatrix} \hat{P}_1 & P_{12} \\ P_{12}^\top & \hat{P}_2 \end{pmatrix} \in \Domcons.
 \end{multline}
 We say that the estimate $\bPhat$ \emph{preserves conservativeness on $\Domcons$} if
 \begin{multline}\label{eq:conservative}
   P_1\preceq\hat{P}_1 \ \&\ P_2\preceq\hat{P}_2 \\
   \Rightarrow
   \begin{pmatrix} K_1 & K_2 \end{pmatrix}
   \begin{pmatrix} P_1 & P_{12} \\ P_{12}^\top & P_2 \end{pmatrix}
   \begin{pmatrix} K_1^\top \\ K_2^\top \end{pmatrix}
   \preceq\bPhat
 \end{multline}
 for all $(K_1,K_2,\hat{x}_1,\hat{x}_2,\hat{P}_1,\hat{P}_2)\in\DomPhat$ and all
 \begin{align*}
  \Pjoint=\begin{pmatrix} P_1 & P_{12} \\ P_{12}^\top & P_2 \end{pmatrix}\in\Domcons.
 \end{align*}
\end{defn}

Under assumption~\eqref{eq:A1}, the meaning of this definition is evident from equations~\eqref{eq:joint_variance} and \eqref{eq:fusion_variance}: No matter what the actual value of the joint error covariance $\Pjoint$ of $\hat{x}_1$ and $\hat{x}_2$ is, if we have conservative estimates $\hat{P}_1\succeq P_1$ and $\hat{P}_2\succeq P_2$ of the error covariances of our two unbiased estimators $\hat{x}_1$ and $\hat{x}_2$, we want to combine them into a conservative estimate $\hat{P}\succeq P$ of the error covariance of the unbiased fused estimator $\hat{x}$, at least as long as the joint error covariance $\Pjoint$ was drawn from the set $\Domcons\subseteq\Ppi$.
Note that Definition~\ref{def:conservative} articulates a worst case principle in the sense that the same estimate $\hat{P}$ is provided for all possible values of $\Pjoint\in\Domcons$.
The set $\Domcons$ is ideally all of $\Ppi$, however, the information fusion literature is also interested in cases where we have \emph{a priori structural} information about $\Pjoint$, and restricting the set $\Domcons$ to encode this information may lead to less conservative solutions, in particular of the optimal conservative information fusion problem discussed in Section~\ref{sec:optimal}.
Note, however, that this construction does not allow to model \emph{stochastic prior} information on $\Pjoint$ (such as an assumption that $\Pjoint$ is drawn from a probability distribution).

It is an immediate consequence of the above definition that an estimate $\bPhat$ that preserves conservativeness on $\Domcons$ also preserves conservativeness on any subset of $\Domcons$ that is closed under diagonal enlargement.
It follows from Lemma~\ref{lem:diag_dom} in Appendix~\ref{sec:Lemmas} and transitivity of the L{\"o}wner order that $\Ppi$ and $\PPpi\subseteq\Ppi$  are both closed under diagonal enlargement and are hence both valid candidates for $\Domcons$.

The following result is sometimes tacitly assumed in the existing information fusion literature. It allows to ``forget'' the conceptual distinction between the \emph{true} error covariances $P_1$ and $P_2$ of the two estimators $\hat{x}_1$ and $\hat{x}_2$ and our conservative estimates $\hat{P}_1\succeq P_1$ and $\hat{P}_2\succeq P_2$ of these quantities when studying the conservativeness of the fused estimate $\hat{P}$ (the corresponding value of the map $\bPhat$).
Depending on the source, sometimes the theory is developed in terms of $P_1$ and $P_2$ (e.g. \cite{reinhardt2015,Uhlmann1996,WuCaiFu2018}), other times in terms of $\hat{P}_1$ and $\hat{P}_2$ (e.g. \cite{arambel2001,chen2002}).
The following result shows that those two approaches are equivalent.
We will make explicit use of this result to preserve the conceptual distinction between the true error covariances and our estimates of these quantities in our exposition.

\begin{prop}\label{prop:conservative}
 The estimate $\bPhat$ preserves conservativeness on $\Domcons$ if and only if
 \begin{equation}\label{eq:conservative_hat}
  \begin{pmatrix} K_1 & K_2 \end{pmatrix}
  \begin{pmatrix} \hat{P}_1 & \hat{P}_{12} \\
   \hat{P}_{12}^\top & \hat{P}_2 \end{pmatrix}
  \begin{pmatrix} K_1^\top \\ K_2^\top \end{pmatrix}
  \preceq\bPhat
 \end{equation}
 for all $(K_1,K_2,\hat{x}_1,\hat{x}_2,\hat{P}_1,\hat{P}_2)\in\DomPhat$ and $\hat{P}_{12}\in\R^{p_1\times p_2}$ with
 \begin{align*}
  \begin{pmatrix}
   \hat{P}_1 & \hat{P}_{12} \\
   \hat{P}_{12}^\top & \hat{P}_2
  \end{pmatrix} \in \Domcons.
 \end{align*}
\end{prop}

\begin{pf}
See Appendix~\ref{app:proofs}.
\end{pf}

The variable $\hat{P}_{12}$ in condition~\eqref{eq:conservative_hat} should not be confused for an estimate of the unknown cross-correlation $P_{12}$.
It merely serves to enumerate all possible off-diagonal terms that complete our estimates $\hat{P}_1$ and $\hat{P}_2$ to a matrix in $\Domcons$.
The following result shows that it can be equivalently replaced by a spectral norm-bounded matrix parameter $X$.
For the cases $\Domcons=\Ppi$ and $\Domcons=\PPpi$ this parameter is otherwise free.
For more general domains $\Domcons\subseteq\Ppi$, an interesting question is how the structural constraints on $\Pjoint$ encoded in the set $\Domcons$ translate into spectral conditions on the parameter $X$.
We suspect that structural constraints where this translation is simple should be of particular interest (cf. also \cite{WuCaiFu2018}).

\begin{thm}\label{th:conservative_1}
 The estimate $\bPhat$ preserves conservativeness on $\Ppi$ (resp. $\PPpi$) if and only if
 \begin{equation}\label{eq:conservative_1}
  \begin{pmatrix} K_1 & K_2 \end{pmatrix}
  \begin{pmatrix} \hat{P}_1 & \hat{P}_1^{\frac{1}{2}}X\hat{P}_2^{\frac{1}{2}} \\
   \hat{P}_2^{\frac{1}{2}}X^\top\hat{P}_1^{\frac{1}{2}} & \hat{P}_2 \end{pmatrix}
  \begin{pmatrix} K_1^\top \\ K_2^\top \end{pmatrix}
  \preceq\bPhat
 \end{equation}
 for all $X\in\R^{p_1\times p_2}$ with $\sigmamax(X)\leq 1$ (resp. $<1$).
\end{thm}

\begin{pf}
 This follows immediately from Proposition~\ref{prop:conservative} and Lemma~\ref{lem:HJ} in Appendix~\ref{sec:Lemmas}.
\qed\end{pf}

Define the maps
\begin{equation}\label{eq:Q1}\begin{split}
 \bQ_1\colon&\DomPhat\longrightarrow\R^{n\times p_1}, \\
 &(K_1,K_2,\hat{x}_1,\hat{x}_2,\hat{P}_1,\hat{P}_2) \mapsto K_1\hat{P}_1^{\frac{1}{2}}
\end{split}\end{equation}
and
\begin{equation}\label{eq:Q2}\begin{split}
 \bQ_2\colon&\DomPhat\longrightarrow\R^{n\times p_2}, \\
 &(K_1,K_2,\hat{x}_1,\hat{x}_2,\hat{P}_1,\hat{P}_2) \mapsto K_2\hat{P}_2^{\frac{1}{2}}.
\end{split}\end{equation}

With the help of these maps, Theorem~\ref{th:conservative_1} can be stated more succinctly as follows.

\begin{cor}\label{cor:conservative_2}
 The estimate $\bPhat$ preserves conservativeness on $\Ppi$ (resp. $\PPpi$) if and only if
 \begin{equation}\label{eq:conservative_2}
   \begin{pmatrix} \bQ_1 & \bQ_2 \end{pmatrix}
   \begin{pmatrix} I & X \\ X^\top & I \end{pmatrix}
   \begin{pmatrix} \bQ_1^\top \\ \bQ_2^\top \end{pmatrix}
   \preceq\bPhat
 \end{equation}
 for all $X\in\R^{p_1\times p_2}$ with $\sigmamax(X)\leq 1$ (resp. $<1$).
\end{cor}

\begin{pf}
See Appendix~\ref{app:proofs}.
\end{pf}

Somewhat surprisingly, the matrix parameter $X$ in Corollary~\ref{cor:conservative_2} can be equivalently replaced by a single scalar parameter through manipulation of the matrix inequality.
We show three different but equivalent versions of this result.
As a further consequence, we also show that there is no difference between the two choices $\Domcons :=\Ppi$ and $\Domcons :=\PPpi$.

\begin{thm}\label{th:conservative_3}
 The following conditions are equivalent.
 \begin{enumerate}
  \item
   $\bPhat$ preserves conservativeness on $\Domcons:=\Ppi$.
  \item
   $\bPhat$ preserves conservativeness on $\Domcons:=\PPpi$.
  \item
   There exists a positive function $\btau\colon\DomPhat\longrightarrow\R$, $\btau>0$ such that
   \begin{equation}\label{eq:conservative_4}
    \bPhat \succeq
    \begin{cases}
     \bQ_2\bQ_2^\top, & \bQ_1=0, \\
     \bQ_1\bQ_1^\top, & \bQ_2=0, \\
     \left\{\begin{tabular}{@{}l@{}}
      $(1+\btau)\bQ_1\bQ_1^\top$ \\ $\quad +(1+\frac{1}{\btau})\bQ_2\bQ_2^\top$
     \end{tabular}\right\}, & \text{otherwise}.
    \end{cases}
   \end{equation}
  \item
   There exists a function $\balpha\colon\DomPhat\longrightarrow\R$, $\balpha\in(0,1)$ such that
   \begin{equation}\label{eq:conservative_5.2}
    \bPhat \succeq
    \begin{cases}
     \bQ_2\bQ_2^\top, & \bQ_1=0, \\
     \bQ_1\bQ_1^\top, & \bQ_2=0, \\
     \frac{1}{\balpha}\bQ_1\bQ_1^\top+\frac{1}{1-\balpha}\bQ_2\bQ_2^\top, & \text{otherwise}.
    \end{cases}
   \end{equation}
  \item
   There exists a function $\balpha\colon\DomPhat\longrightarrow\R$, $\balpha\in[0,1]$ such that
   \begin{equation}\label{eq:conservative_5.3}
    \begin{pmatrix}
     \bPhat     & \bQ_1     & \bQ_2 \\
     \bQ_1^\top & \balpha I & 0 \\
     \bQ_2^\top & 0         & (1-\balpha)I
    \end{pmatrix}
    \succeq 0.
   \end{equation}
 \end{enumerate}
\end{thm}

\begin{pf}
 See Appendix~\ref{sec:proofs_main}.
\end{pf}

We will see in the following sections that (generalized) CI \cite{Uhlmann1996,arambel2001} computes an optimal conservativeness-preserving covariance estimate in the case where assumption~\eqref{eq:A1} holds and $\Pjoint\succ 0$.

\section{Optimal error covariance estimates}\label{sec:optimal}

It is clear from Theorem~\ref{th:conservative_3} that given an estimator~\eqref{eq:linear_fusion}, there is an uncountably infinite family of choices for the error covariance estimate $\hat{P}$ even if we require that it preserves conservativeness.
A natural additional requirement is then that $\hat{P}$ should not be \emph{too conservative}.
One way to formalize this requirement is described in this section.

\begin{defn}\label{def:isotone}
 A function $J\colon\dom(J)\subseteq\mathcal{S}_n\longrightarrow\R$ is called \emph{isotone} if for all $A,B\in\mathcal\dom(J)$,
 \begin{align*}
  A\preceq B \quad\Rightarrow\quad J(A)\leq J(B).
 \end{align*}
 An isotone function $J\colon\dom(J)\subseteq\mathcal{S}_n\longrightarrow\R$ is called \emph{strictly isotone} if for all $A,B\in\mathcal\dom(J)$,
 \begin{align*}
  (A\preceq B \ \& \ A\not=B) \quad\Rightarrow\quad J(A)<J(B).
 \end{align*}
\end{defn}

The matrix determinant ($\det$) on $\PPn$ and the matrix trace ($\Trace$) on $\mathcal{S}_n$ are strictly isotone functions \cite[Theorem~8.4.9]{bernstein}, and so are the restriction of a strictly isotone function to a smaller domain and the composition of a strictly increasing function $\R\longrightarrow\R$ with a strictly isotone function.

We will use strictly isotone \emph{cost functions} to compare the relative levels of conservativeness between different error covariance estimates.
Optimality will then be minimality with respect to the chosen cost function.

The advantage of working with strictly isotone cost functions as per Definition~\ref{def:isotone} is that if two error covariance estimates are comparable in the L{\"o}wner order, their relative rank with respect to any such cost function is the same.
We can hence potentially find an optimal solution that is optimal with respect to any strictly isotone cost function by finding a minimizer in the L{\"o}wner order.
Since the L{\"o}wner order is only a partial order, it is not clear that this can be done.
Indeed it turns out that in the case of interest to us, namely optimal conservative unbiased linear fusion of two (partial) state estimates with \emph{unknown} error cross-covariance, no such minimizer exists.
We can, however, try to find a parameterized family of \emph{lower bounds} for the error covariance estimate, where members of the family are not necessarily mutually comparable in the L{\"o}wner order.
As long as this family is sufficiently rich and contains any optimal solution with respect to a given strictly isotone cost function, an optimal solution can then be found by optimization over the family parameter(s).
As we will see in Section~\ref{sec:KU}, this is indeed possible and even with a family that has a single functional parameter, where the parameter is the function $\balpha\in[0,1]$ as it appears in Part~\textit{(5)} of Theorem~\ref{th:conservative_3} (this latter statement will be proved in Section~\ref{sec:SDP}).

It should be somewhat intuitive at this point why such a result might exist: Clearly, the best choice for $\bPhat$ with respect to the L{\"o}wner order is equality in equation~\eqref{eq:conservative_4} or equation~\eqref{eq:conservative_5.2}.
We will see later that equation~\eqref{eq:conservative_5.2} is particularly useful for understanding the connection to (generalized) CI.
The right hand-side of equation~\eqref{eq:conservative_5.2} depends on $K_1$ and $K_2$ through $\bQ_1$ resp. $\bQ_2$, though, so there is an additional step of finding the coefficients $K_1$ and $K_2$ in equation~\eqref{eq:linear_fusion} that have property~\eqref{eq:linear_unbiased_fusion} and lead to a right-hand side in equation~\eqref{eq:conservative_5.2} that is minimal with respect to the L{\"o}wner order.

Conceptually, what we are trying to do is to find a \emph{linear fusion rule}
\begin{equation}\label{eq:fusion_rule}\begin{split}
 \bF\colon &\DomF \subseteq \Rpo\times\Rpt\times\Ppo\times\Ppt \\
 &\qquad \longrightarrow \bF(\DomF)\subseteq\R^{n\times p_1}\times\R^{n\times p_2}\times\Pn, \\
 &(\hat{x}_1,\hat{x}_2,\hat{P}_1,\hat{P}_2) \mapsto (K_1,K_2,\hat{P})
\end{split}\end{equation}
that is computable by some algorithm and has the unbiasedness and minimality properties that we will spell out in detail below.
We emphasize again that in this framework, $\bF=(\bK_1,\bK_2,\bPhat)$ is a \emph{deterministic} map and not a stochastic estimator (cf. the discussion at the beginning of Section~\ref{sec:conservative}).
It turns out that the optimal conservative unbiased linear fusion rules do not explicitly depend on the values of $\hat{x}_1$ and $\hat{x}_2$ but that is not clear a priori for the $\bPhat$ part of the rule.
What is already clear on an intuitive level is that in such an optimal conservative unbiased linear fusion rule, the relationship between
$(\bK_1,\bK_2)$ and $\bPhat$ will be given by equality in equation~\eqref{eq:conservative_5.2}.

\begin{defn}\label{def:unbiased_fusion_rule}
 The linear fusion rule~\eqref{eq:fusion_rule} is called \emph{unbiased} if for all $(K_1,K_2,\hat{P})\in\bF(\DomF)$, $\hat{x}=K_1\hat{x}_1+K_2\hat{x}_2$ is an unbiased estimator of $\xtrue$ for all $\xtrue\in\Rn$, i.e. if $K_1H_1+K_2H_2=I$.
\end{defn}

In order to apply the results from Section~\ref{sec:conservative} to the estimate $\bPhat$ produced by an unbiased linear fusion rule $(\bK_1,\bK_2,\bPhat)$, we need to interpret $\bPhat$ as a map
\begin{align*}
 \bPhat\colon \DomPhat\subseteq\DomPhatMax\longrightarrow\Pn
\end{align*}
defined on the domain
\begin{align}
 \DomPhat := \{ &(K_1,K_2,\hat{x}_1,\hat{x}_2,\hat{P}_1,\hat{P}_2)\in\DomPhatMax \nonumber\\
             &\quad |\, (\hat{x}_1,\hat{x}_2,\hat{P}_1,\hat{P}_2)\in\DomF \label{eq:DomPhat} \quad \&\\
             &\quad\ \, (K_1,K_2)=(\bK_1,\bK_2)(\hat{x}_1,\hat{x}_2,\hat{P}_1,\hat{P}_2)\}.\nonumber
\end{align}
We will do this tacitly where it is necessary.

Recalling the definition of the maps $\bQ_1$ \eqref{eq:Q1} and $\bQ_2$ \eqref{eq:Q2}, the relationship between $(K_1,K_2)$ and the values $(Q_1,Q_2)$ of these maps is only bijective if we assume that $\hat{P}_1$ and $\hat{P}_2$ are non-singular, i.e. if $\hat{P}_1\succ 0$ and $\hat{P}_2\succ 0$.
Since these are assumed to be conservative estimates of the true covariances $P_1\succeq 0$ and $P_2\succeq 0$, respectively, we can always choose them to be positive definite.
To avoid unnecessary mathematical complications, we will hence define the domain of the linear fusion rule~\eqref{eq:fusion_rule} as
\begin{equation}\label{fusion_rule_domain}
 \DomF := \Rpo\times\Rpt\times\PPpo\times\PPpt.
\end{equation}

\begin{defn}\label{def:conservative_fusion_rule}
 Let $\Domcons\subseteq\Ppi$ be closed under diagonal enlargement.
 The linear fusion rule~\eqref{eq:fusion_rule} \emph{preserves conservativeness on $\Domcons$} or, for brevity, is called \emph{conservative} if
 \begin{multline}\label{eq:conservative_fusion_rule}
   P_1\preceq\hat{P}_1 \ \&\ P_2\preceq\hat{P}_2 \\
   \Rightarrow
   \begin{pmatrix} \bK_1 & \bK_2 \end{pmatrix}
   \begin{pmatrix} P_1 & P_{12} \\ P_{12}^\top & P_2 \end{pmatrix}
   \begin{pmatrix} \bK_1^\top \\ \bK_2^\top \end{pmatrix}
   \preceq\bPhat
 \end{multline}
 for all $(\hat{x}_1,\hat{x}_2,\hat{P}_1,\hat{P}_2)\in\DomF$ and all
 \begin{equation}\label{eq:Pjoint}
  \Pjoint=\begin{pmatrix} P_1 & P_{12} \\ P_{12}^\top & P_2 \end{pmatrix}\in\Domcons.
 \end{equation}
\end{defn}

The reader is invited to check that this definition of preserving conservativeness matches Definition~\ref{def:conservative} when we are using the interpretation of $\bPhat$ as a map defined on $\DomPhat$ \eqref{eq:DomPhat}.
We acknowledge that on a first read it may appear as if we could have avoided the repetitive definitions of unbiasedness and conservativeness for error covariance estimates $\bPhat$ and linear fusion rules $\bF$, however, this would have obscured the fact that the results of Section~\ref{sec:optimal} are results for a \emph{given} set of coefficients $\begin{pmatrix} K_1 & K_2 \end{pmatrix}$ while in this and the following sections, we are trying to \emph{find} these (optimal) coefficients as part of a linear fusion rule.

To match our choice for $\DomF$, we will set $\Domcons:=\PPpi$ meaning that $\Pjoint\succ 0$ for all $\Pjoint\in\Domcons$ and hence in particular $P_1\succ 0$ and $P_2\succ 0$ (this follows from Lemma~\ref{lem:strict_Schur_complement} in Appendix~\ref{sec:Lemmas}).
A further reason why $\Domcons=\PPpi$ is a good choice is as follows.
Consider a conservative unbiased linear fusion rule \eqref{eq:fusion_rule}.
Since $\bK_1H_1+\bK_2H_2=I$ it follows that $\bK:=\begin{pmatrix} \bK_1 & \bK_2 \end{pmatrix}$ has full row rank and hence Lemma~\ref{lem:sandwich} in Appendix~\ref{sec:Lemmas} implies that $\bK\Pjoint\bK^\top\succ 0$.
By inequality~\eqref{eq:conservative_fusion_rule} this implies that $\bPhat\succ 0$ and such a linear fusion rule can hence be applied iteratively, where $\hat{P}$ (the value of $\bPhat$) from the previous iteration becomes $\hat{P}_1$ or $\hat{P}_2$ in the next iteration.
This allows the use of fusion rules as an alternative to dynamic state estimation by filtering (iterative fusion of sensor measurements obtained at different times with a state estimate propagated based on a dynamic system model between measurements) in cases where subsequent measurements are correlated with unknown correlation. This for example occurs in robotics applications involving the modular simultaneous estimation of robot states and unknown environmental features \cite{cdc2026}.

We are now in a position to formally state the optimal conservative unbiased linear information fusion problem for two (partial) state estimates.

\begin{defn}\label{def:optimal_fusion}
 Let $J\colon\PPn\longrightarrow\R$ be strictly isotone.
 Let Assumption~\eqref{eq:A1} hold and denote by $\mathcal{F}_H$ the set of linear fusion rules
 \begin{multline*}
  (\bK_1,\bK_2,\bPhat)\colon \Rpo\times\Rpt\times\PPpo\times\PPpt \\
  \longrightarrow
  \R^{n\times p_1}\times\R^{n\times p_2}\times\PPn
 \end{multline*}
 that are unbiased and preserve conservativeness on $\PPpi$.
 The \emph{optimal conservative unbiased linear information fusion problem} is the optimization problem
 \begin{equation}\label{def:OF}\tag{OF}
  (\bK_1^\ast,\bK_2^\ast,\bPhat^\ast) = \argmin_{(\bK_1,\bK_2,\bPhat)\in\mathcal{F}_H} J(\bPhat).
 \end{equation}
 A solution $(\bK_1^\ast,\bK_2^\ast,\bPhat^\ast)$ of problem~\eqref{def:OF} is called \emph{$J$-optimal} or \emph{optimal} if the choice of $J$ is clear from context.
\end{defn}

Note that as per our usual convention, problem~\eqref{def:OF} is to be read as a pointwise problem in $(\hat{x}_1,\hat{x}_2,\hat{P}_1,\hat{P}_2)\in\Rpo\times\Rpt\times\PPpo\times\PPpt$ and that
a linear fusion rule $(\bK_1^\ast,\bK_2^\ast,\bPhat^\ast)\in\mathcal{F}_H$ is a solution of problem~\eqref{def:OF}, and hence $J$-optimal, if and only if $J(\bPhat^\ast)\leq J(\bPhat)$ for all $(\bK_1,\bK_2,\bPhat)\in\mathcal{F}_H$.

As a first step towards a solution of problem~\eqref{def:OF}, the following theorem provides a formula for optimal unbiased linear fusion of two partial state estimates with \emph{known} joint error covariance $\Pjoint$ \eqref{eq:Pjoint}, including the error cross-covariance $P_{12}$.
For this problem, there exists a unique minimizer with respect to the L{\"o}wner order although the L{\"o}wner order is only a partial order.
This result (including for the cases that we do not treat here, where $H$ does not have full column rank and/or $\Pjoint$ is singular) is known as \emph{Gau{\ss}-Markov theorem} in the applied mathematics literature.
To prove the following theorem, we use a formulation of the Gau{\ss}-Markov theorem due to Albert \cite{albert1973} that is tailor-made for our purposes (see also \cite[Section~8.3]{Ben-Israel_Greville}).

\begin{thm}\label{th:Kopt}
 Let Assumption~\eqref{eq:A1} hold and assume that $\Pjoint\succ 0$.
 Consider the set of error covariances
 \begin{align*}
  \mathcal{P} := \{K\Pjoint K^\top \,|\, KH = I \}
 \end{align*}
 associated with unbiased estimators of the form~\eqref{eq:linear_fusion}, cf. equation~\eqref{eq:fusion_variance} and equation~\eqref{eq:linear_unbiased_fusion}.
 The set $\mathcal{P}$ contains a unique minimal element $P^\ast=K^\ast\Pjoint (K^\ast)^\top$ with respect to the L{\"o}wner order.
 The matrices $K^\ast$ and $P^\ast$ can be computed as follows.
 \begin{equation}\label{eq:Kopt}\begin{split}
  K^\ast &= (H^\top\Pjoint^{-1}H)^{-1}H^\top\Pjoint^{-1} \quad \text{and} \\
  P^\ast &= (H^\top\Pjoint^{-1}H)^{-1}.
 \end{split}\end{equation}
 Furthermore,
 \begin{equation}\label{eq:within_intersection}\begin{split}
  \left(P^\ast\right)^{-1}&\succeq H_1^\top P_1^{-1}H_1 \quad\text{and} \\
  \left(P^\ast\right)^{-1}&\succeq H_2^\top P_2^{-1}H_2.
 \end{split}\end{equation}
\end{thm}

\begin{pf}
 See Appendix~\ref{app:proofs}.
\end{pf}

The next result shows that in the special case $p_1=p_2=n$ and $H_1=H_2=I$, the formula \eqref{eq:Kopt} for the optimal fused estimator is equivalent to the  Bar-Schalom/Campo formula \cite{Bar-Schalom-Campo}.

\begin{cor}\label{cor:BSC}
 In the setting of Theorem~\ref{th:Kopt}, let $p_1=p_2=n$ and $H_1=H_2=I$. Then
 \begin{align*}
  K_1^\ast &= I - (P_1-P_{12})\Delta^{-1}, \\
  K_2^\ast &= (P_1-P_{12})\Delta^{-1} \text{ and } \\
  P^\ast &= P_1 - (P_1-P_{12})\Delta^{-1}(P_1-P_{12}^\top),
 \end{align*}
 where $\Delta:=P_1+P_2-P_{12}-P_{12}^\top$.
\end{cor}

\begin{pf}
 See Appendix~\ref{app:proofs}.
\end{pf}

\section{Ellipsoidal calculus}\label{sec:ell}

Using Theorem~\ref{th:Kopt}, we can now derive a geometric necessary condition for an unbiased linear fusion rule to preserve conservativeness on $\PPpi$.

The geometric interpretation of equation~\eqref{eq:within_intersection} in terms of error ellipsoids is as follows.
For a given \emph{shape matrix} $\Sigma\in\Pn$, the set
\begin{align*}
 \mathcal{E}(\Sigma) := \{ x\in\Rn \,|\, x^\top\Sigma x\leq 1 \}
\end{align*}
is a (possibly degenerate) ellipsoid centered at $0\in\Rn$.
Using Lemma~\ref{lem:kahan} in Appendix~\ref{sec:Lemmas}, equation~\eqref{eq:within_intersection} is equivalent to
\begin{align*}
 \mathcal{E}((P^\ast)^{-1}) \subseteq
 \mathcal{E}(H_1^\top P_1^{-1}H_1) \cap \mathcal{E}(H_2^\top P_2^{-1}H_2).
\end{align*}
The ellipsoid on the left of this inclusion is the error ellipsoid associated with the optimal fused estimator, the two (possibly degenerate) ellipsoids on the right of the inclusion are the two prior error ellipsoids in $\Rn$ that are mapped to the error ellipsoids in $\R^{p_1}$ resp. $\R^{p_2}$ by $H_1$ resp. $H_2$.

The following result shows that when the cross-co\-variance term $P_{12}$ in $\Pjoint$ (cf. equation~\eqref{eq:joint_variance}) ranges over all possible values that lead to a positive definite joint error covariance $\Pjoint$, then all the error ellipsoids associated with the corresponding optimal fused estimators not only lie in the intersection of those two prior error ellipsoids but cover the interior of the intersection.
Our proof is inspired by the proof of statement 2 on Page~1880 of \cite{chen2002}.

\begin{thm}\label{th:fill_intersection}
 Let Assumption~\eqref{eq:A1} hold and let $P_i\in\R^{p_i\times p_i}$, $P_i\succ 0$, $i=1,2$.
 For given $P_{12}\in\R^{p_1\times p_2}$ let
 \begin{align*}
  \Pjoint(P_{12}) := \begin{pmatrix} P_1 & P_{12} \\ P_{12}^\top & P_2 \end{pmatrix}
 \end{align*}
 and consider the set
 \begin{align*}
  \mathcal{P}_{12}^{\succ} := \{ P_{12}\in\R^{p_1\times p_2} \,|\, \Pjoint(P_{12})\succ 0 \}.
 \end{align*}
 For given $P_{12}\in\mathcal{P}_{12}^{\succ}$, denote by
 \begin{align*}
  P^\ast(P_{12}):=(H^\top(\Pjoint(P_{12}))^{-1}H)^{-1}
 \end{align*}
 the optimal error covariance matrix from Theorem~\ref{th:Kopt}.
 Then
 \begin{multline}\label{eq:fill_intersection}
  \mathrm{int}(\mathcal{E}(H_1^\top P_1^{-1}H_1) \cap \mathcal{E}(H_2^\top P_2^{-1}H_2)) \\
  \subseteq
  \bigcup_{P_{12}\in\mathcal{P}_{12}^{\succ}} \mathcal{E}((P^\ast(P_{12}))^{-1}),
 \end{multline}
 where $\mathrm{int}(S)$ denotes the interior of the set $S\subseteq\Rn$.
\end{thm}

\begin{pf}
 See Appendix~\ref{app:proofs}.
\end{pf}

The following example shows that the family of optimal error ellipsoids in Theorem~\ref{th:fill_intersection} does not necessarily cover the entire intersection of the prior error ellipsoids.
It will turn out, though, that covering (at least) the interior of the intersection as per Theorem~\ref{th:fill_intersection} is sufficient for our purposes.

\begin{exmp}\label{ex:not_full_cover}
 Let $n=2$ and $p_1=p_2=1$.
 Let $H_1=\begin{pmatrix} 1 & 0 \end{pmatrix}$ and $H_2=\begin{pmatrix} 0 & 1 \end{pmatrix}$ then $H=I\in\R^{2\times 2}$.
 Let $P_1=P_2=1$ then the conditions of Theorem~\ref{th:fill_intersection} hold, $\mathcal{P}_{12}^{\succ}=(-1,1)$ and
 $\mathcal{E}(H_1^\top P_1^{-1}H_1) \cap \mathcal{E}(H_2^\top P_2^{-1}H_2)=[-1,1]\times[-1,1]\subset\R^2$.
 Let $P_{12}\in\mathcal{P}_{12}^{\succ}=(-1,1)$ then
 \begin{align*}
  (P^\ast(P_{12}))^{-1}&=H^\top(\Pjoint(P_{12}))^{-1}H=(\Pjoint(P_{12}))^{-1} \\
  &= \frac{1}{1-P_{12}^2} \begin{pmatrix} 1\hspace*{-3mm} & -P_{12} \\ -P_{12}\hspace*{-3mm} & 1 \end{pmatrix}
 \end{align*}
 and hence
 \begin{align*}
  \begin{pmatrix} 1 & 1 \end{pmatrix} (P^\ast(P_{12}))^{-1} \begin{pmatrix} 1 \\ 1 \end{pmatrix}
  = \frac{2}{1+P_{12}} > 1.
 \end{align*}
 It follows that $x=\begin{pmatrix} 1 & 1 \end{pmatrix}^\top$ is an element of the intersection
 $\mathcal{E}(H_1^\top P_1^{-1}H_1) \cap \mathcal{E}(H_2^\top P_2^{-1}H_2)$ but not of the union
 $\bigcup_{P_{12}\in\mathcal{P}_{12}^{\succ}} \mathcal{E}((P^\ast(P_{12}))^{-1})$.

 A more careful analysis would show that in this case, the family of optimal error ellipsoids almost fully covers the intersection of the prior error ellipsoids except for the four corner points.
 This observation appears to be of limited use, however.
\end{exmp}

The following result is a consequence of the previous two theorems and provides a geometric necessary condition for an unbiased linear fusion rule to preserve conservativeness on $\PPpi$.

\begin{thm}\label{th:conservative_fusion}
 Let Assumption~\eqref{eq:A1} hold and assume that the linear fusion rule
  \begin{multline*}
  (\bK_1,\bK_2,\bPhat)\colon \Rpo\times\Rpt\times\PPpo\times\PPpt \\
  \longrightarrow
  \R^{n\times p_1}\times\R^{n\times p_2}\times\PPn
 \end{multline*}
 is unbiased and preserves conservativeness on $\PPpi$.
 Then
 \begin{align*}
  \mathcal{E}(\bPhat^{-1})\supseteq\mathcal{E}(H_1^\top\hat{P}_1^{-1}H_1) \cap \mathcal{E}(H_2^\top\hat{P}_2^{-1}H_2).
 \end{align*}
\end{thm}

\begin{pf}
 See Appendix~\ref{sec:proofs_main}.
\end{pf}

The final ingredient for the derivation of lower bounds for the error covariance of any conservative unbiased linear fusion rule is the following theorem due to Kahan \cite{kahan1968,kahan}.\footnote{The paper \cite{kahan} is an updated version of Kahan's famous paper on ellipsoid intersections \cite{kahan1968}. The updated paper has more complete proofs, including for the case of singular shape matrices (degenerate ellipsoids). The papers \cite{kahan1968} and \cite{kahan} do not provide numberings for results and not all the important results are broken out as statements of theorems. For a proof of this statement see \cite[Construction of $M$ on page~438]{kahan1968} and \cite[Interposition Procedure on page~5]{kahan}.}

\begin{thm}\label{th:kahan}
 Let $\Sigma_1,\Sigma_2,\Sigma\in\Pn$.
 Assume that $\mathcal{E}(\Sigma_1)\not\subseteq\mathcal{E}(\Sigma_2)$ and $\mathcal{E}(\Sigma_2)\not\subseteq\mathcal{E}(\Sigma_1)$.
 Then $\mathcal{E}(\Sigma)\supseteq\mathcal{E}(\Sigma_1)\cap\mathcal{E}(\Sigma_2)$ implies that there exists $\alpha\in[0,1]$ such that
 \begin{align*}
  \mathcal{E}(\Sigma)\supseteq\mathcal{E}(\alpha\Sigma_1+(1-\alpha)\Sigma_2)\supseteq\mathcal{E}(\Sigma_1)\cap\mathcal{E}(\Sigma_2).
 \end{align*}
\end{thm}

\section{The Kahan-Uhlmann family}\label{sec:KU}

The next result provides lower bounds for the error covariance of any conservative unbiased linear fusion rule in terms of the L{\"o}wner order.
In contrast to the lower bounds provided by Theorem~\ref{th:conservative_3}, these lower bounds are independent of the coefficients $(K_1,K_2)$.
We will show that there exists a family of fusion rules corresponding to these lower bounds that contains all solutions of problem~\eqref{def:OF}.

\begin{thm}\label{th:lower_bound}
 Let Assumption~\eqref{eq:A1} hold and assume that the linear fusion rule
  \begin{multline*}
  (\bK_1,\bK_2,\bPhat)\colon \Rpo\times\Rpt\times\PPpo\times\PPpt \\
  \longrightarrow
  \R^{n\times p_1}\times\R^{n\times p_2}\times\PPn
 \end{multline*}
 is unbiased and preserves conservativeness on $\PPpi$.
 Then
 \begin{align*}
  \bPhat \succeq \begin{cases}
   (H_2^\top\hat{P}_2^{-1}H_2)^{-1}, &\hspace*{-4em} H_2^\top\hat{P}_2^{-1}H_2 \succeq H_1^\top\hat{P}_1^{-1}H_1, \\
   (H_1^\top\hat{P}_1^{-1}H_1)^{-1}, &\hspace*{-4em} H_1^\top\hat{P}_1^{-1}H_1 \succeq H_2^\top\hat{P}_2^{-1}H_2, \\
   \left\{\begin{tabular}{@{}l@{}}
    $(\alpha H_1^\top\hat{P}_1^{-1}H_1 +$ \\
    $(1-\alpha)H_2^\top\hat{P}_2^{-1}H_2)^{-1}$ \\
    for some $\alpha\in[0,1]$
   \end{tabular}\right\},
   & \text{otherwise}.
  \end{cases}
 \end{align*}
\end{thm}

\begin{pf}
 See Appendix~\ref{sec:proofs_main}.
\end{pf}

Examining the lower bound provided by Theorem~\ref{th:lower_bound} suggests to introduce a one-parameter family of linear fusion rules for two (partial) state estimates as follows.

\begin{defn}\label{def:kahan_family}
 Denote by $\mathcal{K}_H^{\balpha}$ the set of all functions
 \begin{align*}
  \balpha\colon &\Rpo\times\Rpt\times\PPpo\times\PPpt \longrightarrow [0,1], \\
  &(\hat{x}_1,\hat{x}_2,\hat{P}_1,\hat{P}_2) \mapsto \balpha(\hat{x}_1,\hat{x}_2,\hat{P}_1,\hat{P}_2)
 \end{align*}
 with
 \begin{align*}
  \balpha = \begin{cases}
   0, &\hspace*{-4em} H_2^\top\hat{P}_2^{-1}H_2 \succ H_1^\top\hat{P}_1^{-1}H_1, \\
   \alpha\in[0,1], &\hspace*{-4em} H_2^\top\hat{P}_2^{-1}H_2 = H_1^\top\hat{P}_1^{-1}H_1, \\
   1, &\hspace*{-4em} H_2^\top\hat{P}_2^{-1}H_2 \prec H_1^\top\hat{P}_1^{-1}H_1, \\
   \left\{\begin{tabular}{@{}l@{}}
    $\alpha\in[0,1]$ such that \\
    $\alpha H_1^\top\hat{P}_1^{-1}H_1 +$ \\
    $(1-\alpha)H_2^\top\hat{P}_2^{-1}H_2$ \\
    is non-singular
   \end{tabular}\right\},
   & \text{otherwise}.
  \end{cases}
 \end{align*}
 For $\balpha\in\mathcal{K}_H^{\balpha}$ define the linear fusion rule
 \begin{align*}
  (\bK_1^\CIa,\bK_2^\CIa,\bPhat^\CIa)\colon &\Rpo\times\Rpt\times\PPpo\times\PPpt \\
  &\qquad \longrightarrow
  \R^{n\times p_1}\times\R^{n\times p_2}\times\PPn, \\
  &(\hat{x}_1,\hat{x}_2,\hat{P}_1,\hat{P}_2) \mapsto (K_1,K_2,\hat{P})
 \end{align*}
 by
 \begin{equation}\label{eq:CIa}\begin{split}
  \bPhat^\CIa &:= (\balpha H_1^\top\hat{P}_1^{-1}H_1 + (1-\balpha)H_2^\top\hat{P}_2^{-1}H_2)^{-1}, \\
  \bK_1^\CIa &:= \balpha\bPhat^\CIa H_1^\top\hat{P}_1^{-1}, \\
  \bK_2^\CIa &:= (1-\balpha)\bPhat^\CIa H_2^\top\hat{P}_2^{-1}.
 \end{split}\end{equation}
 The set
 \begin{align*}
  \mathcal{K}_H:=\{(\bK_1^\CIa,\bK_2^\CIa,\bPhat^\CIa)\,|\,\balpha\in\mathcal{K}_H^{\balpha}\}
 \end{align*}
 of all such linear fusion rules is called the \emph{Kahan-Uhlmann family}.
\end{defn}

The following result shows that the Kahan-Uhlmann family is contained in the feasible set for problem~\eqref{def:OF}.
Its proof depends on Theorem~\ref{th:conservative_3}.

\begin{thm}\label{th:kahan_family}
 Assume that Assumption~\eqref{eq:A1} holds then $\mathcal{K}_H^{\balpha}\not=\emptyset$ and $(\bK_1^\CIa,\bK_2^\CIa,\bPhat^\CIa)$ is well-defined, unbiased, and preserves conservativeness on $\PPpi$ for all $\balpha\in\mathcal{K}_H^{\balpha}$, i.e. $\emptyset\not=\mathcal{K}_H\subseteq\mathcal{F}_H$.
\end{thm}

\begin{pf}
 See Appendix~\ref{sec:proofs_main}.
\end{pf}

The following result shows that the parametrization of the Kahan-Uhlmann family $\mathcal{K}_H$ through $\mathcal{K}_H^{\balpha}$ is one-to-one.

\begin{thm}\label{th:kahan_parameter}
 Let $\balpha,\balpha'\in\mathcal{K}_H^{\balpha}$ then $(\bK_1^\CIa,\bK_2^\CIa,\bPhat^\CIa)=(\bK_1^{\CIa'},\bK_2^{\CIa'},\bPhat^{\CIa'})$ if and only if $\balpha=\balpha'$.
\end{thm}

\begin{pf}
 See Appendix~\ref{app:proofs}.
\end{pf}

The following theorem shows that any linear fusion rule that is optimal within the Kahan-Uhlmann family is in fact optimal amongst \emph{all} conservative unbiased linear fusion rules and is hence a solution of problem~\eqref{def:OF}.

\begin{thm}\label{th:Kahan_implies_opt}
 Let $J\colon\PPn\longrightarrow\R$ be strictly isotone and let Assumption~\eqref{eq:A1} hold.
 If the linear fusion rule $(\bK_1^\ast,\bK_2^\ast,\bPhat^\ast)\in\mathcal{K}_H$ is $J$-optimal over $\mathcal{K}_H$ then it is $J$-optimal.
\end{thm}

\begin{pf}
See Appendix~\ref{sec:proofs_main}.
\end{pf}

The following theorem provides a converse result to the previous theorem: any solution of problem~\eqref{def:OF} is necessarily a member of the Kahan-Uhlmann family.

\begin{thm}\label{th:opt_implies_Kahan}
 Let $J\colon\PPn\longrightarrow\R$ be strictly isotone and let Assumption~\eqref{eq:A1} hold.
 If the linear fusion rule $(\bK_1^\ast,\bK_2^\ast,\bPhat^\ast)\in\mathcal{F}_H$ is $J$-optimal then it is a member of the Kahan-Uhlmann family $\mathcal{K}_H$ and is $J$-optimal over $\mathcal{K}_H$.
\end{thm}

\begin{pf}
 See Appendix~\ref{sec:proofs_main}.
\end{pf}

\section{Covariance Intersection}\label{sec:CI}

Theorem~\ref{th:Kahan_implies_opt} and Theorem~\ref{th:opt_implies_Kahan} imply that solutions of problem~\eqref{def:OF} correspond to optimal choices $\balpha^\ast$ for the family parameter $\balpha\in\mathcal{K}_H^{\balpha}$ in the Kahan-Uhlmann family $\mathcal{K}_H\subseteq\mathcal{F}_H$.
In fact, we can even show that such an optimal choice $\balpha^\ast$ is optimal amongst \emph{all} functions $\balpha\in[0,1]$, not just those in $\mathcal{K}_H^{\balpha}$.
In order to formalize this statement, we will extend the calculus of strictly isotone cost functions to also cover the case where
$\alpha H_1^\top\hat{P}_1^{-1}H_1 + (1-\alpha)H_2^\top\hat{P}_2^{-1}H_2$ is singular and hence has no (finite) matrix inverse.

\begin{defn}\label{def:isotone_ext}
 Let $J\colon\PPn\longrightarrow\R$ be (strictly) isotone.
 The \emph{extension} of $J$ is the function
 \begin{align*}
  \bar{J}\colon &\PPn\cup\{+\infty\}\longrightarrow\R\cup\{+\infty\}, \\
  & P\mapsto\begin{cases}
   +\infty, & P=+\infty, \\
   J(P), & \text{otherwise}.
  \end{cases}
 \end{align*}
 For any \emph{singular} $P\in\Pn\subseteq\mathcal{S}_n$ define $P^{-1}:=+\infty\in\PPn\cup\{+\infty\}$.
\end{defn}

Using Definition~\ref{def:isotone_ext}, we can now formally state the optimization problem known as (generalized) CI in the information fusion literature \cite{Uhlmann1996,arambel2001}.

\begin{defn}\label{def:covariance_intersection}
 Let $J\colon\PPn\longrightarrow\R$ be strictly isotone and let $\bar{J}$ denote its extension.
 The \emph{Covariance Intersection problem} is the optimization problem
 \begin{equation}\label{def:CI}\tag{CI}\begin{split}
  \balpha^\ast = \argmin_{\balpha\in[0,1]} \bar{J}((&\balpha H_1^\top\hat{P}_1^{-1}H_1 \\
  &\quad + (1-\balpha)H_2^\top\hat{P}_2^{-1}H_2)^{-1}),
 \end{split}\end{equation}
 where the optimization runs over all functions
 \begin{align*}
  \balpha\colon &\Rpo\times\Rpt\times\PPpo\times\PPpt \longrightarrow [0,1], \\
  &(\hat{x}_1,\hat{x}_2,\hat{P}_1,\hat{P}_2) \mapsto \balpha(\hat{x}_1,\hat{x}_2,\hat{P}_1,\hat{P}_2).
 \end{align*}
 Given a solution $\balpha^\ast\in[0,1]$ of problem~\eqref{def:CI}, the \emph{Covariance Intersection algorithm} computes the linear fusion rule
 $(\bK_1^\CI,\bK_2^\CI,\bPhat^\CI)$ by
 \begin{align*}
  \bPhat^\CI &:= (\balpha^\ast H_1^\top\hat{P}_1^{-1}H_1 + (1-\balpha^\ast)H_2^\top\hat{P}_2^{-1}H_2)^{-1}, \\
  \bK_1^\CI &:= \balpha^\ast\bPhat^\CI H_1^\top\hat{P}_1^{-1}, \\
  \bK_2^\CI &:= (1-\balpha^\ast)\bPhat^\CI H_2^\top\hat{P}_2^{-1}.
 \end{align*}
\end{defn}

Similar to problem~\eqref{def:OF}, problem~\eqref{def:CI} is to be read as a pointwise problem in $(\hat{x}_1,\hat{x}_2,\hat{P}_1,\hat{P}_2)\in\Rpo\times\Rpt\times\PPpo\times\PPpt$.
Note that the formula for the linear fusion rule computed by the Covariance Intersection algorithm is identical to the parametrization formula~\eqref{eq:CIa} of the Kahan-Uhlmann family $\mathcal{K}_H$, but at this point we have not yet shown that a solution $\balpha^\ast\in[0,1]$ of problem~\eqref{def:CI} is a valid Kahan-Uhlmann parameter.
The following result implies that the linear fusion rules associated with solutions of problem~\eqref{def:CI} (if they exist) are always members of the Kahan-Uhlmann family and are hence well-defined.

\begin{thm}\label{th:CI_implies_Kahan}
 Let $J\colon\PPn\longrightarrow\R$ be strictly isotone and let Assumption~\eqref{eq:A1} hold.
 Let the function $\balpha'\colon\Rpo\times\Rpt\times\PPpo\times\PPpt \longrightarrow [0,1]$ be such that
 \begin{align*}
  \bar{J}((\balpha' H_1^\top\hat{P}_1^{-1}H_1 + (1-\balpha')H_2^\top\hat{P}_2^{-1}H_2)^{-1})
  \leq \bar{J}(\bPhat^\CIa)
 \end{align*}
 for some $\balpha\in\mathcal{K}_H^{\balpha}$ then $\balpha'\in\mathcal{K}_H^{\balpha}$.
 In particular, if $\balpha^\ast\in[0,1]$ is a solution of problem~\eqref{def:CI} then $\balpha^\ast\in\mathcal{K}_H^{\balpha}$.
\end{thm}

\begin{pf}
 See Appendix~\ref{sec:proofs_main}.
\end{pf}

The following corollary sums up the resulting relationship between the optimization problems \eqref{def:OF} and \eqref{def:CI}.

\begin{cor}\label{cor:reduction_to_CI}
 Let $J\colon\PPn\longrightarrow\R$ be strictly isotone and let Assumption~\eqref{eq:A1} hold.
 Then problem~\eqref{def:OF} and problem~\eqref{def:CI} are equivalent in the sense that $(\bK_1^\ast,\bK_2^\ast,\bPhat^\ast)\in\mathcal{F}_H$ is a solution of problem~\eqref{def:OF} if and only if there exists a solution $\balpha^\ast\in[0,1]$ of problem~\eqref{def:CI} such that $(\bK_1^\ast,\bK_2^\ast,\bPhat^\ast)=(\bK_1^\CI,\bK_2^\CI,\bPhat^\CI)$.
\end{cor}

\begin{pf}
 See Appendix~\ref{sec:proofs_main}.
\end{pf}

\section{A semi-definite program}\label{sec:SDP}

We will now show that problem~\eqref{def:OF}, and hence by Corollary~\ref{cor:reduction_to_CI} also problem~\eqref{def:CI}, is equivalent to the following family of semi-definite programs.

\begin{defn}\label{def:semidefinite_program}
 Let $J\colon\PPn\longrightarrow\R$ be strictly isotone, define
 $\DomSDP := \R \times \R^{n\times p_1} \times \R^{n\times p_2} \times \PPn$
 and the family of semi-definite programs
 \begin{equation}\label{def:SDP}\tag{SDP}\begin{split}
  &(\balpha^\ast,\bK_1^\ast,\bK_2^\ast,\bPhat^\ast) := \argmin_{(\balpha,\bK_1,\bK_2,\bPhat)\in\DomSDP} J(\bPhat) \\
  &\quad \begin{pmatrix}
    \bPhat                            & \bK_1\hat{P}_1^{\frac{1}{2}} & \bK_2\hat{P}_2^{\frac{1}{2}} \\
    \hat{P}_1^{\frac{1}{2}}\bK_1^\top & \balpha I                    & 0 \\
    \hat{P}_2^{\frac{1}{2}}\bK_2^\top & 0                            & (1-\balpha)I
   \end{pmatrix}
   \succeq 0, \\
   &\qquad \bK_1H_1+\bK_2H_2=I.
 \end{split}\end{equation}
\end{defn}

As with problem~\eqref{def:OF} and problem~\eqref{def:CI}, problem~\eqref{def:SDP} is to be read as a pointwise problem in $(\hat{x}_1,\hat{x}_2,\hat{P}_1,\hat{P}_2)\in\Rpo\times\Rpt\times\PPpo\times\PPpt$ with each being a separate semi-definite program.
Note that the above semi-definite program is different to both the optimization problems studied in \cite[Eq. (7)-(9)]{WangShenZhu2019} and
\cite[Eq. (27)-(28)]{WangShenZhu2019}.

\begin{thm}\label{th:reduction_to_SDP}
 Let $J\colon\PPn\longrightarrow\R$ be strictly isotone and let Assumption~\eqref{eq:A1} hold.
 Then problem~\eqref{def:OF} and problem~\eqref{def:SDP} are equivalent in the sense that $(\bK_1^\ast,\bK_2^\ast,\bPhat^\ast)\in\mathcal{F}_H$ is a solution of problem~\eqref{def:OF} if and only if there exists $\balpha^\ast\in\R$ such that $(\balpha^\ast,\bK_1^\ast,\bK_2^\ast,\bPhat^\ast)$ is a solution of problem~\eqref{def:SDP}.
\end{thm}

\begin{pf}
 See Appendix\ref{sec:proofs_main}.
\end{pf}

Theorem~\ref{th:reduction_to_SDP} and Corollary~\ref{cor:reduction_to_CI} imply that any solution $(\balpha^\ast,\bK_1^\ast,\bK_2^\ast,\bPhat^\ast)\in\DomSDP$ of problem \eqref{def:SDP} corresponds to a solution $\balpha'\in[0,1]$ of problem \eqref{def:CI} such that
$(\bK_1^\ast,\bK_2^\ast,\bPhat^\ast)=(\bK_1^{\balpha'},\bK_2^{\balpha'},\bPhat^{\balpha'})$.
It is, however, not clear if this necessarily implies that $\balpha^\ast=\balpha'$.
The following two results answer this question in the affirmative.

\begin{thm}\label{th:alpha}
 Let Assumption~\eqref{eq:A1} hold, let $\balpha\in\mathcal{K}_H^{\balpha}$, let $\balpha'\colon\Rpo\times\Rpt\times\PPpo\times\PPpt\longrightarrow\R$ and let
 \begin{equation}\label{eq:alpha_alpha'}
  \begin{pmatrix}
    \bPhat^{\balpha}                              & \bK_1^{\balpha}\hat{P}_1^{\frac{1}{2}} & \bK_2^{\balpha}\hat{P}_2^{\frac{1}{2}} \\
    \hat{P}_1^{\frac{1}{2}}(\bK_1^{\balpha})^\top & \balpha' I                         & 0 \\
    \hat{P}_2^{\frac{1}{2}}(\bK_2^{\balpha})^\top & 0                                      & (1-\balpha')I
  \end{pmatrix}
  \succeq 0.
 \end{equation}
 Then $\balpha=\balpha'$.
\end{thm}

\begin{pf}
 See Appendix~\ref{app:proofs}.
\end{pf}

\begin{thm}\label{th:SDP_equiv_CI}
 Let $J\colon\PPn\longrightarrow\R$ be strictly isotone and let Assumption~\eqref{eq:A1} hold.
 Then problem~\eqref{def:SDP} and problem~\eqref{def:CI} are equivalent in the sense that
 $(\balpha^\ast\!,\bK_1^\ast,\bK_2^\ast,\bPhat^\ast)\in\DomSDP$ is a solution of problem \eqref{def:SDP} if and only if $\balpha^\ast$ is a solution of problem~\eqref{def:CI} and $(\bK_1^\ast,\bK_2^\ast,\bPhat^\ast)=(\bK_1^\CI,\bK_2^\CI,\bPhat^\CI)$.
\end{thm}

\begin{pf}
 See Appendix~\ref{sec:proofs_main}.
\end{pf}

Note that none of the above results implies that optimal solutions for the three equivalent optimization problems \eqref{def:OF}, \eqref{def:CI} and \eqref{def:SDP} exist.
We will tackle this question in the following section by studying problem~\eqref{def:CI} in more detail.

\section{Optimal solutions for Covariance Intersection}\label{sec:CI_solution}

The following result provides a fairly general condition on the cost function $J\colon\mathcal{S}_n\longrightarrow\R$ under which problem~\eqref{def:CI} is guaranteed to have a solution.
It is likely that this condition is only sufficient and not necessary but it is sufficiently general to cover the two most important choices $J:=\Trace$ and $J:=\det$.

\begin{defn}\label{def:spectrally_unbounded}
 A function $J\colon\mathcal{S}_n\longrightarrow\R$ is called \emph{spectrally unbounded on $\Pn$} if
 \begin{align*}
  \lim_{n\to\infty}J(P_n)=+\infty
 \end{align*}
 for every sequence $\{P_n\}_{n\in\N}\subset\PPn$ with $\lambda_{\min}(P_n)\geq\epsilon$ for some $\epsilon>0$ and all $n\in\N$ and $\lim\limits_{n\to\infty}\lambda_{\max}(P_n)=+\infty$.
 Here, $\lambda_{\min}(P_n)$ resp. $\lambda_{\max}(P_n)$ denotes the smallest resp. largest eigenvalue of $P_n$.
\end{defn}

Both the determinant and trace cost functions are strictly isotone, continuous and spectrally unbounded on $\Pn$.

\begin{thm}\label{th:CI}
 Let $J\colon\mathcal{S}_n\longrightarrow\R$ be strictly isotone, continuous and spectrally unbounded on $\Pn$ and let Assumption~\eqref{eq:A1} hold.
 Then problem~\eqref{def:CI} has a solution.
\end{thm}

\begin{pf}
 See Appendix~\ref{sec:proofs_main}.
\end{pf}

For the cost function $J=\det$, we will show below that problem~\eqref{def:CI} has a unique solution except for the special case where $H_2^\top\hat{P}_2^{-1}H_2=H_1^\top\hat{P}_1^{-1}H_1$ because in this case, the function
$\alpha \mapsto \alpha H_1^\top\hat{P}_1^{-1}H_1+(1-\alpha)H_2^\top\hat{P}_2^{-1}H_2$ is constant and hence any value of $\alpha\in[0,1]$ gives the same value for $\bar{J}((\hat{P}^\alpha)^{-1})=J((\hat{P}^\alpha)^{-1})$, meaning that problem~\eqref{def:CI} never has a unique solution unless this special case is excluded from consideration.
Note that by Part~\textit{(3)} of Lemma~\ref{lem:corner_cases} in Appendix~\ref{sec:Lemmas}, $\alpha H_1^\top\hat{P}_1^{-1}H_1+(1-\alpha)H_2^\top\hat{P}_2^{-1}H_2$ is invertible for all $\alpha\in(0,1)$ and hence it is invertible for all $\alpha\in[0,1]$ in the special case where $H_2^\top\hat{P}_2^{-1}H_2=H_1^\top\hat{P}_1^{-1}H_1$.

An equation characterizing the unique solution for the case $J=\det$ was derived in \cite[Theorem~3 and Corollary~1]{sabater2002}, however, as the following theorem shows this equation only characterizes the part of the solution inside the open interval $(0,1)$.
Counter to intuition and contrary to the statement of \cite[Theorem~3]{sabater2002}, the optimal value for $\alpha\in[0,1]$ can be $\alpha^\ast=0$ or $\alpha^\ast=1$ even when there is no inclusion relationship between the two prior error ellipsoids.
We provide a concrete example for this following the theorem.

\begin{thm}\label{th:CI_det}
 Let $J:=\det$ and let Assumption~\eqref{eq:A1} hold.
 For $\alpha\in[0,1]$ define
 \begin{align*}
  \bSigma_\alpha &:= \alpha H_1^\top\hat{P}_1^{-1}H_1+(1-\alpha)H_2^\top\hat{P}_2^{-1}H_2 \quad \text{and} \\
  \bDelta(\alpha) &:= \Trace(\adj(\bSigma_\alpha)(\bSigma_1-\bSigma_0))
 \end{align*}
 and note that $\bSigma_0=H_2^\top\hat{P}_2^{-1}H_2$, $\bSigma_1=H_1^\top\hat{P}_1^{-1}H_1$ and $\Delta(\alpha)$ is a polynomial of degree at most $n-1$.
 The solution set of problem~\eqref{def:CI} consists of the functions
 \begin{align*}
  \balpha^\ast = \begin{cases}
   0, &\hspace*{-3em} \bDelta(0)\leq 0\ \&\ \bSigma_0\not=\bSigma_1, \\
   \alpha^\ast\in[0,1], &\hspace*{-3em} \bSigma_0=\bSigma_1, \\
   1, &\hspace*{-3em} \bDelta(1)\geq 0\ \&\ \bSigma_0\not=\bSigma_1, \\
   \left\{\begin{tabular}{@{}l@{}}
    the unique solution of\mbox{ } \\
    $\bDelta(\alpha)=0$, $\alpha\in(0,1)$
   \end{tabular}\right\},
   & \text{otherwise}.
  \end{cases}
 \end{align*}
\end{thm}

\begin{pf}
 See Appendix~\ref{sec:proofs_main}.
\end{pf}

\begin{exmp}\label{ex:CI_det}
 Let $n=p_1=p_2=2$, $H_1=H_2=I$ and
 \begin{align*}
  \hat{P}_1=\begin{pmatrix} 1 & 0 \\ 0 & 1 \end{pmatrix}\succ 0
  \quad \text{and} \quad
  \hat{P}_2=\begin{pmatrix} 1.25 & 0 \\ 0 & 0.1 \end{pmatrix}\succ 0.
 \end{align*}
 Then $\Sigma_0\not=\Sigma_1$ and $\Delta(\alpha)=-3.6\alpha-5.2<0$ for all $\alpha\in[0,1]$, and hence $\alpha^\ast=0$ although
 \begin{align*}
 \Sigma_1-\Sigma_0=\begin{pmatrix} 0.2 & 0 \\ 0 & -9 \end{pmatrix}
 \end{align*}
 and therefore $\Sigma_0\not\succeq\Sigma_1$ and $\Sigma_0\not\preceq\Sigma_1$.
 The derivative of $\alpha\mapsto\det(\hat{P}^\alpha)=\det((\Sigma_\alpha)^{-1})$ is $\frac{10(9\alpha+13)}{(9\alpha-10)^2(\alpha+4)^2}>0$ for all $\alpha\in[0,1]$ in this case.
\end{exmp}

To the best of the authors' knowledge, a characterization similar to Theorem~\ref{th:CI_det} is not known for the case $J:=\Trace$ and would be much harder to obtain since the matrix trace does not play nicely with matrix inversion and matrix multiplication.
We can, however, derive a different equation that must be fulfilled by solutions to problem~\eqref{def:CI} in this case.
A similar equation can be found in \cite{chen2002} and sheds some light on how the two prior pieces of information are weighted relative to each other in $\Trace$-optimal solutions to problem~\eqref{def:CI}.
See \cite[Theorem~1]{WuCaiFu2018} for an interpretation of the denominator that occurs in this equation.
Note, however, that solving this implicit equation for $\balpha^\ast$ is no simpler than solving problem~\eqref{def:CI} directly.

\begin{thm}\label{th:CI_tr}
 Let Assumption~\eqref{eq:A1} hold and let $J:=\Trace$.
 Then
 \begin{align*}
  \balpha^\ast =
   \frac{\sqrt{\Trace(\bK_1^{\balpha^\ast}\hat{P}_1(\bK_1^{\balpha^\ast})^\top)}}
        {\sqrt{\Trace(\bK_1^{\balpha^\ast}\hat{P}_1(\bK_1^{\balpha^\ast})^\top)}
         + \sqrt{\Trace(\bK_2^{\balpha^\ast}\hat{P}_2(\bK_2^{\balpha^\ast})^\top)}}
 \end{align*}
 for all solutions $\balpha^\ast\in[0,1]$ of problem~\eqref{def:CI}.
\end{thm}

\begin{pf}
 See Appendix~\ref{sec:proofs_main}.
\end{pf}

\section{Conclusion}\label{sec:conclusion}

We have generalized the known optimality result for CI of full state estimates to the case of partial state estimates, a case of particular importance in distributed and network applications.
Enroute we provided several new results that are of independent interest, including three different necessary and sufficient characterizations of conservative fused error covariance estimates, a generalization of the Bar-Schalom/Campo formula to partial state estimates, and a correction to the previously known characterizing equation of optimal solutions for the matrix determinant cost function.
The more formal treatment of the domain of preservation of conservativeness and its relation to spectral constraints on the error cross-covariance opens a path to new ways of incorporating prior knowledge of the structure of the joint error covariance that is the subject of ongoing and future research.

\section*{Acknowledgements}

The authors would like to thank the anonymous reviewers for pointing out two embarrassing errors in a previous version of the paper and for valuable suggestions on how to improve the presentation of the material.

%

\appendix
\section{Lemmas}\label{sec:Lemmas}

In this appendix we provide a series of lemmas that are well-known but where we require very specific versions for our purposes that can be difficult to track down in the literature.

\begin{lem}\label{lem:diag_dom}
 Let $P_1,\hat{P}_1,P_2,\hat{P}_2\in\Sn$.
 \begin{enumerate}
  \item
   If $P_1\prec\hat{P}_1$ and $P_2\prec\hat{P}_2$ then
   \begin{align*}
    \begin{pmatrix} P_1 & P_{12} \\ P_{12}^\top & P_2 \end{pmatrix} \prec
    \begin{pmatrix} \hat{P}_1 & P_{12} \\ P_{12}^\top & \hat{P}_2 \end{pmatrix}.
   \end{align*}
  \item
   If $P_1\preceq\hat{P}_1$ and $P_2\preceq\hat{P}_2$ then
   \begin{align*}
    \begin{pmatrix} P_1 & P_{12} \\ P_{12}^\top & P_2 \end{pmatrix} \preceq
    \begin{pmatrix} \hat{P}_1 & P_{12} \\ P_{12}^\top & \hat{P}_2 \end{pmatrix}.
   \end{align*}
 \end{enumerate}
\end{lem}

\begin{pf}
 Compute
 \begin{align*}
  \Delta :=
  \begin{pmatrix} P_1 & P_{12} \\ P_{12}^\top & P_2 \end{pmatrix} -
  \begin{pmatrix} \hat{P}_1 & P_{12} \\ P_{12}^\top & \hat{P}_2 \end{pmatrix} =
  \begin{pmatrix} P_1-\hat{P}_1 & 0 \\ 0 &\! P_2-\hat{P}_2 \end{pmatrix}.
 \end{align*}
 Then $\Delta\prec 0$ if $P_1-\hat{P}_1\prec 0$ and $P_2-\hat{P}_2\prec 0$. This proves part \textit{(1)}.
 Similarly, $\Delta\preceq 0$ if $P_1-\hat{P}_1\preceq 0$ and $P_2-\hat{P}_2\preceq 0$. This proves part \textit{(2)}.
\qed\end{pf}

\begin{lem}\label{lem:sandwich}
 Let $P_1,P_2\in\Pn$ and $S\in\R^{m\times n}$.
 \begin{enumerate}
  \item
  If $P_2\succ P_1$ and $S$ has rank $m$ then $SP_2S^\top\succ SP_1S^\top$.
  \item
  If $P_2\succeq P_1$ then $SP_2S^\top\succeq SP_1S^\top$.
 \end{enumerate}
\end{lem}

\begin{pf}
 Let $P_2\succ P_1$ then it follows that $x^\top(S(P_2-P_1)S^\top)x=(S^\top x)^\top(P_2-P_1)(S^\top x) > 0$ for all $x\in\R^m$ with
 $S^\top x\not=0$.
 If $S^\top$ has full column rank, these are all $x\in\R^m$ with $x\not=0$.
 It then follows that $SP_2S^\top - SP_1S^\top=S(P_2-P_1)S^\top\succ 0$ and hence that $SP_2S^\top\succ SP_1S^\top$, proving part \textit{(1)}.

 Similarly, if $P_2\succeq P_1$ it follows that $(P_2-P_1)\succeq 0$ and hence that $x^\top(S(P_2-P_1)S^\top)x=(S^\top x)^\top(P_2-P_1)(S^\top x)\geq 0$ for all $x\in\R^m$.
 This implies that $SP_2S^\top - SP_1S^\top=S(P_2-P_1)S^\top\succeq 0$ and it follows that $SP_2S^\top\succeq SP_1S^\top$, proving part \textit{(2)}.
\qed\end{pf}

\begin{lem}\label{lem:strict_Schur_complement}
 Let $Q\in\R^{q\times q}$ with $Q^\top=Q$, $R\in\R^{r\times r}$ with $R^\top=R$ and $S\in\R^{q\times r}$ then the following are equivalent.
 \begin{enumerate}
  \item
   \begin{align*}
    T := \begin{pmatrix} Q & S \\ S^\top & R \end{pmatrix} \succ 0,
   \end{align*}
  \item $Q\succ 0$ and $R-S^\top Q^{-1} S\succ 0$,\text{ and}
  \item $R\succ 0$ and $Q-SR^{-1}S^\top\succ 0$.
 \end{enumerate}
\end{lem}

\begin{pf}
 The equivalence of \textit{(1)} and \textit{(2)} follows from \cite[Theorem~7.7.6]{HJ} by setting $A:=Q$, $B:=S$ and $C:=R$.

 Let
 \begin{align*}
  J:=\begin{pmatrix} 0 & I \\ I & 0 \end{pmatrix}
 \end{align*}
 then $J$ has full row rank and $J^\top=J=J^{-1}=(J^{-1})^\top$.
 Let $T\succ 0$ then Lemma~\ref{lem:sandwich} implies that
 \begin{align*}
  JTJ = \begin{pmatrix} R & S^\top \\ S & Q \end{pmatrix} \succ 0
 \end{align*}
 and, conversely, if $JTJ\succ 0$ then Lemma~\ref{lem:sandwich} implies that $T=J^{-1}(JTJ)J^{-1}\succ 0$.
 It follows that $T\succ 0$ if and only if $JTJ\succ 0$.

 The equivalence of \textit{(1)} and \textit{(3)} follows from \cite[Theorem~7.7.6]{HJ} applied to $JTJ$ by setting $A:=R$, $B:=S^\top$ and $C:=Q$.
\qed\end{pf}

\begin{lem}\label{lem:Schur_complement}
 Let $Q\in\R^{q\times q}$ with $Q^\top=Q$, $R\in\R^{r\times r}$ with $R^\top=R$ and $S\in\R^{q\times r}$ then the following are equivalent.
 \begin{enumerate}
  \item
   \begin{align*}
    T := \begin{pmatrix} Q & S \\ S^\top & R \end{pmatrix} \succeq 0,
   \end{align*}
  \item $Q\succeq 0$, $R-S^\top Q^+ S\succeq 0$ and $(QQ^+-I)S=0$,\text{ and}
  \item $R\succeq 0$, $Q-SR^+S^\top\succeq 0$ and $S(I-RR^+)=0$.
 \end{enumerate}
 Here, $Q^+$ and $R^+$ denote the Moore-Penrose pseudo-inverses of $Q$ and $R$, respectively.
\end{lem}

\begin{pf}
 The equivalence of \textit{(1)} and \textit{(2)} follows from \cite[Theorem~1(i)]{Albert1969} by setting $m:=r$, $n:=q$, $S_{11}:=Q$, $S_{12}:=S$ and $S_{22}:=R$.

 Let
 \begin{align*}
  J:=\begin{pmatrix} 0 & I \\ I & 0 \end{pmatrix}
 \end{align*}
 then $J^\top=J=J^{-1}=(J^{-1})^\top$.
 Let $T\succeq 0$ then Lemma~\ref{lem:sandwich} implies that
 \begin{align*}
  JTJ = \begin{pmatrix} R & S^\top \\ S & Q \end{pmatrix} \succeq 0
 \end{align*}
 and, conversely, if $JTJ\succeq 0$ then Lemma~\ref{lem:sandwich} implies that $T=J^{-1}(JTJ)J^{-1}\succeq 0$.
 It follows that $T\succeq 0$ if and only if $JTJ\succeq 0$.

 The equivalence of \textit{(1)} and \textit{(3)} follows from \cite[Theorem~1(i)]{Albert1969} by setting $m:=q$, $n:=r$, $S_{11}:=R$, $S_{12}:=S^\top$ and $S_{22}:=Q$ and using the fact that $(RR^+)^\top=RR^+$.
\qed\end{pf}

\begin{lem}\label{lem:spectral_radius}
  Let $P_1\in\R^{p_1\times p_1}$, $P_1\succ 0$, $P_2\in\R^{p_2\times p_2}$, $P_2\succ 0$ and $P_{12}\in\R^{p_1\times p_2}$.
  Let $\rho(.)$ denote the spectral radius and let $\sigmamax(.)$ denote the largest singular value.
  Then $\rho(P_{12}^\top P_1^{-1}P_{12}P_2^{-1})=(\sigmamax(P_1^{-\frac{1}{2}}P_{12}P_2^{-\frac{1}{2}}))^2$.
\end{lem}

\begin{pf}
 Compute
 \begin{align*}
  \rho(P_{12}^\top P_1^{-1}P_{12}P_2^{-1})
  &= \rho(P_2^{-\frac{1}{2}}P_{12}^\top P_1^{-1}P_{12}P_2^{-1}P_2^{\frac{1}{2}}) \\
  &= \rho(P_2^{-\frac{1}{2}}P_{12}^\top P_1^{-\frac{1}{2}}P_1^{-\frac{1}{2}}P_{12}P_2^{-\frac{1}{2}}) \\
  &= \rho((P_1^{-\frac{1}{2}}P_{12}P_2^{-\frac{1}{2}})^\top(P_1^{-\frac{1}{2}}P_{12}P_2^{-\frac{1}{2}})) \\
  &= (\sigmamax(P_1^{-\frac{1}{2}}P_{12}P_2^{-\frac{1}{2}}))^2.
 \end{align*}
\qed\end{pf}

\begin{lem}\label{lem:HJ}
 Let $P_1\in\R^{p_1\times p_1}$, $P_2\in\R^{p_2\times p_2}$ and $P_{12}\in\R^{p_1\times p_2}$ and let $\sigmamax(.)$ denote the largest singular value.
 \begin{enumerate}
 \item
  Then
  \begin{equation}\label{eq:HJ.1}
   \begin{pmatrix} P_1 & P_{12} \\ P_{12}^\top & P_2 \end{pmatrix} \succ 0
  \end{equation}
  if and only if $P_1\succ 0$, $P_2\succ 0$ and there exists $X\in\R^{p_1\times p_2}$ with $\sigmamax(X)<1$ such that $P_{12}=P_1^{\frac{1}{2}}XP_2^{\frac{1}{2}}$.
 \item
  Similarly,
  \begin{equation}\label{eq:HJ.2}
   \begin{pmatrix} P_1 & P_{12} \\ P_{12}^\top & P_2 \end{pmatrix} \succeq 0
  \end{equation}
  if and only if $P_1\succeq 0$, $P_2\succeq 0$ and there exists $X\in\R^{p_1\times p_2}$ with $\sigmamax(X)\leq 1$ such that $P_{12}=P_1^{\frac{1}{2}}XP_2^{\frac{1}{2}}$.
 \end{enumerate}
\end{lem}

\begin{pf}
 Assume that inequality~\eqref{eq:HJ.1} holds then Lemma~\ref{lem:strict_Schur_complement} implies that $P_1\succ 0$, $P_2\succ 0$ and $P_2-P_{12}^\top P_1^{-1}P_{12}\succ 0$.
 It follows that $P_2\succ P_{12}^\top P_1^{-1}P_{12}$ where $P_{12}^\top P_1^{-1}P_{12}\succeq 0$ by Lemma~\ref{lem:sandwich}.
 But then \cite[Theorem~7.7.3]{HJ} implies that $\rho(P_{12}^\top P_1^{-1}P_{12}P_2^{-1})<1$ and hence it follows from Lemma~\ref{lem:spectral_radius} that $\sigmamax(X)<1$, where $X:=P_1^{-\frac{1}{2}}P_{12}P_2^{-\frac{1}{2}}$.
 It follows that $P_{12}=P_1^{\frac{1}{2}}XP_2^{\frac{1}{2}}$.

 Conversely, let $P_1\succ 0$ and $P_2\succ 0$ and assume that there exists $X\in\R^{p_1\times p_2}$ with $\sigmamax(X)<1$ such that $P_{12}=P_1^{\frac{1}{2}}XP_2^{\frac{1}{2}}$.
 Then $X=P_1^{-\frac{1}{2}}P_{12}P_2^{-\frac{1}{2}}$ and it follows from Lemma~\ref{lem:spectral_radius} that
 $\rho(P_{12}^\top P_1^{-1}P_{12}P_2^{-1})<1$.
 It follows from Lemma~\ref{lem:sandwich} that $P_{12}^\top P_1^{-1}P_{12}\succeq 0$.
 But then \cite[Theorem~7.7.3]{HJ} implies that $P_2\succ P_{12}^\top P_1^{-1}P_{12}$ and hence that
 $P_2-P_{12}^\top P_1^{-1}P_{12}\succ 0$.
 It follows from Lemma~\ref{lem:strict_Schur_complement} that inequality~\eqref{eq:HJ.1} holds.
 This proves part \textit{(1)}.

 Assume that inequality~\eqref{eq:HJ.2} holds then Lemma~\ref{lem:Schur_complement} implies that $P_1\succeq 0$ and $P_2\succeq 0$.

 It follows from Lemma~\ref{lem:diag_dom} that
 \begin{align*}
   \begin{pmatrix} P_1+\frac{1}{n} I & P_{12} \\ P_{12}^\top & P_2+\frac{1}{n} I \end{pmatrix} \succ 0
 \end{align*}
 for all $n\in\N$.
 By part \textit{(1)} this implies that there exists a sequence $(X_n\in\R^{p_1\times p_2})_{n\in\N}$ with $\sigmamax(X_n)<1$ and such that $P_{12}=P_1^{\frac{1}{2}}X_nP_2^{\frac{1}{2}}$ for all $n\in\N$.
 By the Bolzano-Weierstra{\ss} theorem the sequence $(X_n)$ contains a convergent subsequence, and defining $X\in\R^{p_1\times p_2}$ as the limit of this subsequence it follows that $\sigmamax(X)\leq 1$ and $P_{12}=P_1^{\frac{1}{2}}XP_2^{\frac{1}{2}}$.

 Conversely, let $P_1\succeq 0$ and $P_2\succeq 0$ and assume that there exists $X\in\R^{p_1\times p_2}$ with $\sigmamax(X)\leq 1$ such that $P_{12}=P_1^{\frac{1}{2}}XP_2^{\frac{1}{2}}$.
 Let $(P_{1,n}\succ 0)_{n\in\N}$, $(P_{2,n}\succ 0)_{n\in\N}$ and $(X_n\in\R^{p_1\times p_2})_{n\in\N}$ be sequences with $\sigmamax(X_n)<1$ for all $n\in\N$, $\lim_{n\to\infty}P_{1,n}=P_1$, $\lim_{n\to\infty}P_{2,n}=P_2$ and $\lim_{n\to\infty}X_n=X$.
 Then it follows from part \textit{(1)} that
 \begin{align*}
   \begin{pmatrix} P_{1,n} & P_{1,n}^{\frac{1}{2}}X_nP_{2,n}^{\frac{1}{2}} \\
   (P_{1,n}^{\frac{1}{2}}X_nP_{2,n}^{\frac{1}{2}})^\top & P_{2,n} \end{pmatrix} \succ 0
 \end{align*}
 for all $n\in\N$. By taking the limit $n\to\infty$ it follows that inequality~\eqref{eq:HJ.2} holds.
 This proves part \textit{(2)}.
\qed\end{pf}

\begin{lem}\label{lem:Petersen}
 Let $G\in\Rnn$ with $G^\top=G$, $M\in\R^{n\times p_1}$ with $M\not=0$ and $N\in\R^{p_2\times n}$ with $N\not=0$.
 Then the following are equivalent.
 \begin{enumerate}
  \item $G+MYN+N^\top Y^\top M^\top \preceq 0$ for all $Y\in\R^{p_1\times p_2}$ with $\sigmamax(Y)<1$.
  \item $G+MYN+N^\top Y^\top M^\top \preceq 0$ for all $Y\in\R^{p_1\times p_2}$ with $\sigmamax(Y)\leq 1$.
  \item There exists $\epsilon>0$ such that $G+\epsilon MM^\top+\frac{1}{\epsilon}N^\top N \preceq 0$.
 \end{enumerate}
\end{lem}

\begin{pf}
 Assume that $G+MYN+N^\top Y^\top M^\top \preceq 0$ for all $Y\in\R^{p_1\times p_2}$ with $\sigmamax(Y)<1$.
 Let $\bar{Y}\in\R^{p_1\times p_2}$ with $\sigmamax(\bar{Y})=1$.
 There exists a sequence $(Y_n\in\R^{p_1\times p_2})_{n\in\N}$ such that $\sigmamax(Y_n)<1$ for all $n\in\N$ and $\lim_{n\to\infty}Y_n=\bar{Y}$.
 This implies that $G+MY_nN+N^\top Y_n^\top M^\top \preceq 0$ for all $n\in N$ and hence also that
 $G+M\bar{Y}N+N^\top \bar{Y}^\top M^\top \preceq 0$.
 It follows that $G+MYN+N^\top Y^\top M^\top \preceq 0$ for all $Y\in\R^{p_1\times p_2}$ with $\sigmamax(Y)\leq 1$.

 Conversely, assume that $G+MYN+N^\top Y^\top M^\top \preceq 0$ for all $Y\in\R^{p_1\times p_2}$ with $\sigmamax(Y)\leq 1$ then it follows in particular that $G+MYN+N^\top Y^\top M^\top \preceq 0$ for all $Y\in\R^{p_1\times p_2}$ with $\sigmamax(Y)<1$.
 It follows that \textit{(1)} and \textit{(2)} are equivalent.

 The equivalence of \textit{(2)} and \textit{(3)} is a restatement of Petersen's lemma as given on page~1933 of \cite{KS2008} in our notation.
\qed\end{pf}

\begin{lem}\label{lem:block_inv}
 Let $Q\in\R^{q\times q}$ with $Q^\top=Q$, $R\in\R^{r\times r}$ with $R^\top=R$ and $S\in\R^{q\times r}$ and let
 \begin{align*}
  T := \begin{pmatrix} Q & S \\ S^\top & R \end{pmatrix}.
 \end{align*}
 Define $\Delta_Q := R-S^\top Q^{-1}S$ and $\Delta_R := Q-SR^{-1}S^\top$.
 \begin{enumerate}
  \item If $T\succ0$ then $Q\succ 0$, $R\succ 0$, $\Delta_Q\succ 0$ and $\Delta_R\succ 0$.
  \item If $Q$ and $\Delta_Q$ are non-singular then $T$ is non-singular and
   \begin{align*}
    T^{-1} &= \begin{pmatrix}
              Q^{-1} + Q^{-1}S\Delta_Q^{-1}S^\top Q^{-1}  &  -Q^{-1}S\Delta_Q^{-1}  \\
              -\Delta_Q^{-1}S^\top Q^{-1}                 &  \Delta_Q^{-1}
             \end{pmatrix} \\
           &=: \begin{pmatrix} Q^{-1} & 0 \\ 0 & 0 \end{pmatrix} + \Gamma_Q.
   \end{align*}
   If $T\succ 0$ then $\Gamma_Q\succeq 0$.
  \item If $R$ and $\Delta_R$ are non-singular then $T$ is non-singular and
   \begin{align*}
    T^{-1} &= \begin{pmatrix}
              \Delta_R^{-1}         &  -\Delta_R^{-1}SR^{-1}  \\
              -R^{-1}S^\top\Delta_R^{-1}  &  R^{-1}S^\top\Delta_R^{-1}SR^{-1} + R^{-1}
             \end{pmatrix} \\
           &=: \Gamma_R + \begin{pmatrix} 0 & 0 \\ 0 & R^{-1} \end{pmatrix}.
   \end{align*}
   If $T\succ 0$ then $\Gamma_R\succeq 0$.
  \item If the matrices $Q$, $R$, $\Delta_Q$ and $\Delta_R$ are non-singular then $\Delta_Q^{-1} = R^{-1}S^\top\Delta_R^{-1}SR^{-1}+R^{-1}$ and $\Delta_R^{-1} = Q^{-1}+Q^{-1}S\Delta_Q^{-1}S^\top Q^{-1}$.
 \end{enumerate}
\end{lem}

\begin{pf}
 Part \textit{(1)} follows directly from Lemma~\ref{lem:strict_Schur_complement}. The formulas for $T^{-1}$ in Parts \textit{(2)} and \textit{(3)} are a reformulation of \cite[Proposition~2.8.7]{bernstein} in our notation.
 That $T\succ0$ implies $\Gamma_Q\succeq 0$ and $\Gamma_R\succeq 0$ follows from part \textit{(1)} and Lemma~\ref{lem:Schur_complement}.
 Part \textit{(4)} follows by comparing the relevant blocks in the formulas for $T^{-1}$ from Parts \textit{(2)} and \textit{(3)}.
\qed\end{pf}

\begin{lem}\label{lem:corner_cases}
 Let Assumption~\eqref{eq:A1} hold and let $\hat{P}_1\succ 0$ and $\hat{P}_2\succ 0$.
 \begin{enumerate}
  \item If $H_2^\top\hat{P}_2^{-1}H_2\succeq H_1^\top\hat{P}_1^{-1}H_1$ then $p_2=n$ and $H_2$ is non-singular square.
  \item If $H_1^\top\hat{P}_1^{-1}H_1\succeq H_2^\top\hat{P}_2^{-1}H_2$ then $p_1=n$ and $H_1$ is non-singular square.
  \item $\alpha H_1^\top\hat{P}_1^{-1}H_1 + (1-\alpha)H_2^\top\hat{P}_2^{-1}H_2\succ 0$ for all $\alpha\in(0,1)$.
 \end{enumerate}
\end{lem}

\begin{pf}
 Let $H_2^\top\hat{P}_2^{-1}H_2\succeq H_1^\top\hat{P}_1^{-1}H_1$.
 Since $H_2$ has full row rank, there exists a non-singular $S\in\Rnn$ such that $H_2 S=\begin{pmatrix} 0 & I \end{pmatrix}$ where $I\in\R^{p_2\times p_2}$ is the identity matrix.
 Let $\begin{pmatrix} M & N \end{pmatrix} := H_1 S$, where $M\in\R^{p_1\times p_2}$.
 It follows from Lemma~\ref{lem:sandwich} that
 \begin{align*}
  \begin{pmatrix} 0 & 0 \\ 0 & \hat{P}_2^{-1} \end{pmatrix} &= S^\top H_2^\top\hat{P}_2^{-1}H_2 S \\
  &\succeq S^\top H_1^\top\hat{P}_1^{-1}H_1 S \\
  &= \begin{pmatrix} M^\top\hat{P}_1^{-1}M &\ M^\top\hat{P}_1^{-1}N \\ N^\top\hat{P}_1^{-1}M &\ N^\top\hat{P}_1^{-1}N \end{pmatrix}
 \end{align*}
 and hence that
 \begin{align*}
  \begin{pmatrix}
   -M^\top\hat{P}_1^{-1}M & -M^\top\hat{P}_1^{-1}N \\
   -N^\top\hat{P}_1^{-1}M &\ \hat{P}_2^{-1}\!-N^\top\hat{P}_1^{-1}N
  \end{pmatrix} \succeq 0.
 \end{align*}
 It follows from Lemma~\ref{lem:Schur_complement} that $-M^\top\hat{P}_1^{-1}M\succeq 0$.
 Since $\hat{P}_1^{-1}\succeq 0$, Lemma~\ref{lem:sandwich} implies that also $M^\top\hat{P}_1^{-1}M\succeq 0$ and by anti-symmetry of the L{\"o}wner order it follows that $M^\top\hat{P}_1^{-1}M=0$ and hence that $M=0$.
 But then it follows that
 \begin{align*}
  \rank\begin{pmatrix} 0 & 0 \\ 0 & I \end{pmatrix} =
  \rank\begin{pmatrix} 0 & N \\ 0 & I \end{pmatrix} =
  \rank\begin{pmatrix} H_1 S \\ H_2 S \end{pmatrix} =
  \rank\begin{pmatrix} H_1 \\ H_2 \end{pmatrix} = n
 \end{align*}
 and this implies $p_2=n$.
 It follows that $H_2$ is square, and since it has full row rank, it follows that $H_2$ is non-singular.
 This completes the proof of part \textit{(1)}.

 The proof of part \textit{(2)} follows along the same lines swapping the roles of $H_1$ (resp. $\hat{P}_1$) and $H_2$ (resp. $\hat{P}_2$).

 It follows from Lemma~\ref{lem:strict_Schur_complement} that
 \begin{align*}
  \begin{pmatrix} \alpha\hat{P}_1 & 0 \\ 0 & (1-\alpha)\hat{P_2} \end{pmatrix} \succ 0
 \end{align*}
 for all $\alpha\in(0,1)$.
 By Lemma~\ref{lem:sandwich} this implies that
 \begin{multline*}
  \alpha H_1^\top\hat{P}_1^{-1}H_1 + (1-\alpha)H_2^\top\hat{P}_2^{-1}H_2 \\
  =
  \begin{pmatrix} H_1^\top & H_2^\top \end{pmatrix}
  \begin{pmatrix} \alpha\hat{P}_1 & 0 \\ 0 & (1-\alpha)\hat{P_2} \end{pmatrix}
  \begin{pmatrix} H_1 \\ H_2 \end{pmatrix}
  \succ 0
 \end{multline*}
 for all $\alpha\in(0,1)$.
 This completes the proof of part \textit{(3)}.
\qed\end{pf}

\begin{lem}\label{lem:kahan}
 Let $\Sigma_1,\Sigma_2\in\Pn$. Then $\Sigma_1\succeq\Sigma_2$ if and only if $\mathcal{E}(\Sigma_1)\subseteq\mathcal{E}(\Sigma_2)$.
\end{lem}

\begin{pf}
 This is restatement of \cite[Corollary on Page~4]{kahan} in our language.
\qed\end{pf}

\begin{lem}\label{lem:swap12}
 Let Assumption~\eqref{eq:A1} hold and assume that $\Pjoint\succ 0$.
 Let $P^\ast=(H^\top\Pjoint^{-1}H)^{-1}$ be the optimal error covariance from Theorem~\ref{th:Kopt}.
 Let
 \begin{align*}
  &J := \begin{pmatrix} 0 & I \\ I & 0 \end{pmatrix}, \quad
  \bar{H} := JH = \begin{pmatrix} H_2 \\ H_1 \end{pmatrix} \quad \text{and} \\
  &\bar{P}_{\mathrm{joint}} := J\Pjoint J^\top = \begin{pmatrix} P_2 & P_{12}^\top \\ P_{12} & P_1 \end{pmatrix},
 \end{align*}
 where the first (resp. second) block row in the matrix $J$ has dimension $p_1$ (resp. $p_2$).
 Then $\rank(\bar{H})=n$, $\bar{P}_{\mathrm{joint}}\succ 0$,
 $\bar{P}^\ast := (\bar{H}^\top\bar{P}_{\mathrm{joint}}^{-1}\bar{H})^{-1} = P^\ast$ and
 $\mathcal{E}((\bar{P}^\ast)^{-1})=\mathcal{E}((P^\ast)^{-1})$.
\end{lem}

\begin{pf}
 The matrix $J$ is non-singular with $J^{-1}=J^\top$ hence $\rank(\bar{H})=\rank(JH)=\rank(H)=n$ and $\bar{P}_{\mathrm{joint}}\succ 0$ by Lemma~\ref{lem:sandwich}.
 Since $(\bar{P}_{\mathrm{joint}})^{-1}=(J\Pjoint J^\top)^{-1}=J\Pjoint^{-1}J^\top$, it follows that
 \begin{align*}
  (\bar{P}^\ast)^{-1} &= \bar{H}^\top\bar{P}_{\mathrm{joint}}^{-1}\bar{H}
  = H^\top J^\top J\Pjoint^{-1}J^\top JH \\
  &= H^\top\Pjoint^{-1}H
  = (P^\ast)^{-1}
 \end{align*}
 and hence that $\bar{P}^\ast = P^\ast$ and $\mathcal{E}((\bar{P}^\ast)^{-1})=\mathcal{E}((P^\ast)^{-1})$.\!\!\!\!\qed
\end{pf}

\section{Proofs of auxiliary results}\label{app:proofs}

\textit{Proof of Proposition~\ref{prop:conservative}.}
 Assume that $\bPhat$ fulfils condition~\eqref{eq:conservative} for all $(K_1,K_2,\hat{x}_1,\hat{x}_2,\hat{P}_1,\hat{P}_2)\in\DomPhat$ and all $\Pjoint\in\Domcons$.
 Let $(K_1,K_2,\hat{x}_1,\hat{x}_2,\hat{P}_1,\hat{P}_2)\in\DomPhat$ and let $\hat{P}_{12}\in\R^{p_1\times p_2}$ be such that
 \begin{align*}
  \begin{pmatrix}
   \hat{P}_1 & \hat{P}_{12} \\
   \hat{P}_{12}^\top & \hat{P}_2
  \end{pmatrix} \in \Domcons.
 \end{align*}

 Since $\hat{P}_1\preceq\hat{P}_1$ and $\hat{P}_2\preceq\hat{P}_2$, this matrix fulfils the antecedent in the implication~\eqref{eq:conservative}.
 It follows that the consequent holds as well and that is condition~\eqref{eq:conservative_hat}.

 Conversely, assume that $\bPhat$ fulfils condition~\eqref{eq:conservative_hat} for all
 $(K_1,K_2,\hat{x}_1,\hat{x}_2,\hat{P}_1,\hat{P}_2)\in\DomPhat$ and $\hat{P}_{12}\in\R^{p_1\times p_2}$ with
 \begin{align*}
  \begin{pmatrix}
   \hat{P}_1 & \hat{P}_{12} \\
   \hat{P}_{12}^\top & \hat{P}_2
  \end{pmatrix} \in \Domcons.
 \end{align*}
 Let $(K_1,K_2,\hat{x}_1,\hat{x}_2,\hat{P}_1,\hat{P}_2)\in\DomPhat$ and let
 \begin{align*}
  \Pjoint=\begin{pmatrix} P_1 & P_{12} \\ P_{12}^\top & P_2 \end{pmatrix} \in \Domcons
 \end{align*}
 and such that $P_1\preceq\hat{P}_1$ and $P_2\preceq\hat{P}_2$.
 Since $\Domcons$ is closed under diagonal enlargement, it follows that
 \begin{align*}
  \begin{pmatrix} \hat{P}_1 & P_{12} \\ P_{12}^\top & \hat{P}_2 \end{pmatrix} \in \Domcons
 \end{align*}
 and by our assumption of condition~\eqref{eq:conservative_hat} this implies that
 \begin{equation}\label{eq:cons_hat}
  \begin{pmatrix} K_1 & K_2 \end{pmatrix}
  \begin{pmatrix} \hat{P}_1 & P_{12} \\
   P_{12}^\top & \hat{P}_2 \end{pmatrix}
  \begin{pmatrix} K_1^\top \\ K_2^\top \end{pmatrix}
  \preceq\bPhat.
 \end{equation}
 Furthermore,
 \begin{align*}
  \begin{pmatrix} P_1 & P_{12} \\ P_{12}^\top & P_2 \end{pmatrix} -
  \begin{pmatrix} \hat{P}_1 & P_{12} \\ P_{12}^\top & \hat{P}_2 \end{pmatrix} =
  \begin{pmatrix} P_1-\hat{P}_1 & 0 \\ 0 & P_2-\hat{P}_2 \end{pmatrix} \preceq 0
 \end{align*}
 and hence
 \begin{align*}
  \begin{pmatrix} P_1 & P_{12} \\ P_{12}^\top & P_2 \end{pmatrix} \preceq
  \begin{pmatrix} \hat{P}_1 & P_{12} \\ P_{12}^\top & \hat{P}_2 \end{pmatrix}.
 \end{align*}
 By Lemma~\ref{lem:sandwich} it follows that
 \begin{multline*}
  \begin{pmatrix} K_1 & K_2 \end{pmatrix}
  \begin{pmatrix} P_1 & P_{12} \\ P_{12}^\top & P_2 \end{pmatrix}
  \begin{pmatrix} K_1^\top \\ K_2^\top \end{pmatrix} \\
  \preceq
  \begin{pmatrix} K_1 & K_2 \end{pmatrix}
  \begin{pmatrix} \hat{P}_1 & P_{12} \\
   P_{12}^\top & \hat{P}_2 \end{pmatrix}
  \begin{pmatrix} K_1^\top \\ K_2^\top \end{pmatrix}
  \preceq\bPhat,
 \end{multline*}
 where the last inequality is~\eqref{eq:cons_hat}.
 Transitivity of the Loewner order now implies that condition~\eqref{eq:conservative} holds.
\qed

\textit{Proof of Corollary~\ref{cor:conservative_2}.}
 This follows directly from
 \begin{multline*}
  \begin{pmatrix} K_1 & K_2 \end{pmatrix}
  \begin{pmatrix} \hat{P}_1 & \hat{P}_1^{\frac{1}{2}}X\hat{P}_2^{\frac{1}{2}} \\
   \hat{P}_2^{\frac{1}{2}}X^\top\hat{P}_1^{\frac{1}{2}} & \hat{P}_2 \end{pmatrix}
  \begin{pmatrix} K_1^\top \\ K_2^\top \end{pmatrix} \\
  \begin{aligned}
   &=
   \begin{pmatrix} K_1 & K_2 \end{pmatrix}
   \begin{pmatrix} \hat{P}_1^{\frac{1}{2}}  & 0 \\ 0 & \hat{P}_2^{\frac{1}{2}} \end{pmatrix}
   \begin{pmatrix} I & X \\ X^\top & I \end{pmatrix}
   \begin{pmatrix} \hat{P}_1^{\frac{1}{2}}  & 0 \\ 0 & \hat{P}_2^{\frac{1}{2}} \end{pmatrix}
   \begin{pmatrix} K_1^\top \\ K_2^\top \end{pmatrix}\!\!\!\!\! \\
   &=
   \begin{pmatrix} K_1\hat{P}_1^{\frac{1}{2}} & K_2\hat{P}_2^{\frac{1}{2}} \end{pmatrix}
   \begin{pmatrix} I & X \\ X^\top & I \end{pmatrix}
   \begin{pmatrix} \hat{P}_1^{\frac{1}{2}}K_1^\top \\ \hat{P}_2^{\frac{1}{2}}K_2^\top \end{pmatrix}
  \end{aligned}
 \end{multline*}
 for all $(K_1,K_2,\hat{x}_1,\hat{x}_2,\hat{P}_1,\hat{P}_2)\in\DomPhatMax$.
\qed

\textit{Proof of Theorem~\ref{th:Kopt}.}
 Assume that $P^{\ast,1},P^{\ast,2}\in\mathcal{P}$ are such that $P\succeq P^{\ast,1}$ and $P\succeq P^{\ast,2}$ for all $P\in\mathcal{P}$ then in particular $P^{\ast,2}\succeq P^{\ast,1}$ and $P^{\ast,1}\succeq P^{\ast,2}$.
 It follows from the anti-symmetry of the L{\"o}wner order that then $P^{\ast,1}=P^{\ast,2}$ and hence $\mathcal{P}$ contains at most one minimal element with respect to the L{\"o}wner order, i.e. if $\mathcal{P}$ contains such a minimal element it is necessarily unique.

 Since $H$ has full column rank by Assumption~\eqref{eq:A1}, it follows that the range of $H^\top$ equals (all of) $\R^n$.
 Applying \cite[equation~(6)]{albert1973} with $z:=\begin{pmatrix}\hat{x}_1^\top & \hat{x}_2^\top\end{pmatrix}^\top$ and $x:=\xtrue$ (and hence $V^2=\Pjoint$) implies that $\mathrm{e}_i^\top\hat{x}$, where $\hat{x}:=K^\ast z$ and
 $K^\ast:=(H^\top\Pjoint^{-1}H)^\dagger H^\top\Pjoint^{-1}$, is an unbiased linear estimator for $\mathrm{e}_i^\top\xtrue$ for all $i\in\{1,\dots,n\}$.
 Here, $(\cdot)^\dagger$ denotes the Moore-Penrose pseudoinverse and $\mathrm{e}_i$, $i=1,\dots,n$ denotes the $i$-th standard basis vector of $\R^n$.\footnote{See the reference to Scheffe on Page~183 of \cite{albert1973} to check that $\mathrm{e}_i^\top\xtrue$ is estimable for all $i\in\{1,\dots,n\}$.}
 It follows that $\hat{x}=K^\ast z$ is an unbiased linear estimator for $\xtrue$ with the property that if $\bar{x}=\bar{K}z$ is any other unbiased linear estimator (of $\xtrue$) then $\Cov(\bar{x})\succeq\Cov(\hat{x})$ \cite[equation~(3)]{albert1973}.\footnote{Note the unfortunate typo in \cite[equation~(3)]{albert1973}, where the right-hand side is missing the ``hat'' on top of $x$. See Page~185 of \cite{albert1973} for what is actually being proved.}
 Since unbiasedness of a linear estimator $Kz$ for $\xtrue$ is equivalent to $KH=I$, see equation \eqref{eq:linear_unbiased_fusion}, it follows that $P^\ast:=K^\ast\Pjoint(K^\ast)^\top=\Cov(\hat{x})$ is a minimal element of $\mathcal{P}$ with respect to the L{\"o}wner order and hence the unique minimal element of $\mathcal{P}$.

 Since $H^\top$ has full row rank, it follows from Lemma~\ref{lem:sandwich} that $H^\top\Pjoint^{-1}H\succ 0$ and hence $(H^\top\Pjoint^{-1}H)^\dagger=(H^\top\Pjoint^{-1}H)^{-1}$.
 It follows that
 \begin{align*}
  K^\ast = (H^\top\Pjoint^{-1}H)^{-1}H^\top\Pjoint^{-1}
 \end{align*}
 and
 \begin{align*}
  P^\ast &= K^\ast\Pjoint(K^\ast)^\top \\
  &= (H^\top\Pjoint^{-1}H)^{-1}H^\top\Pjoint^{-1}H(H^\top\Pjoint^{-1}H)^{-1} \\
  &= (H^\top\Pjoint^{-1}H)^{-1}.
 \end{align*}

 But then
 \begin{align*}
  (P^\ast)^{-1} &= H^\top\Pjoint^{-1}H \\
  &= \begin{pmatrix} H_1^\top & H_2^\top \end{pmatrix}
     \left(\begin{pmatrix} P_1^{-1} & 0 \\ 0 & 0 \end{pmatrix} + \Gamma_{P_1}\right)
     \begin{pmatrix} H_1 \\ H_2 \end{pmatrix} \\
  &= H_1^\top P_1^{-1}H_1 + H^\top\Gamma_{P_1}H
  \succeq H_1^\top P_1^{-1}H_1,
 \end{align*}
 where we have used Part \textit{(2)} of Lemma~\ref{lem:block_inv} and the fact that $\Gamma_{P_1}\succeq 0$ implies $H^\top\Gamma_{P_1}H\succeq 0$ by Lemma~\ref{lem:sandwich}.
 Similarly,
 \begin{align*}
  (P^\ast)^{-1} &= H^\top\Pjoint^{-1}H \\
  &= \begin{pmatrix} H_1^\top & H_2^\top \end{pmatrix}
     \left(\Gamma_{P_2} + \begin{pmatrix} 0 & 0 \\ 0 & P_2^{-1} \end{pmatrix}\right)
     \begin{pmatrix} H_1 \\ H_2 \end{pmatrix} \\
  &=  H^\top\Gamma_{P_2}H + H_2^\top P_2^{-1}H_2
  \succeq H_2^\top P_2^{-1}H_2,
 \end{align*}
 where we have used Part \textit{(3)} of Lemma~\ref{lem:block_inv} and the fact that $\Gamma_{P_2}\succeq 0$ implies $H^\top\Gamma_{P_2}H\succeq 0$ by Lemma~\ref{lem:sandwich}.
\qed

\textit{Proof of Corollary~\ref{cor:BSC}.}
 $\Pjoint\succ 0$ implies $P_1\succ 0$ and $\Delta_{P_1} := P_2-P_{12}^\top P_1^{-1}P_{12}\succ 0$ by Lemma~\ref{lem:strict_Schur_complement}.
 It hence follows from Lemma~\ref{lem:block_inv} \textit{(2)} applied to $\Pjoint$ that
 \begin{align*}
  \Pjoint^{-1}
  &= \begin{pmatrix} P_1 & P_{12} \\ P_{12}^\top & P_2 \end{pmatrix}^{-1} \\
  &= \begin{pmatrix}
     P_1^{-1} + P_1^{-1}P_{12}\Delta_{P_1}^{-1}P_{12}^\top P_1^{-1} & - P_1^{-1}P_{12}\Delta_{P_1}^{-1} \\
     - \Delta_{P_1}^{-1}P_{12}^\top P_1^{-1} & \Delta_{P_1}^{-1}
    \end{pmatrix}
 \end{align*}
 and hence
 \begin{align*}
  H^\top\Pjoint^{-1}H
  &= \begin{pmatrix} I &\ I \end{pmatrix}
     \begin{pmatrix} P_1 & P_{12} \\ P_{12}^\top & P_2 \end{pmatrix}^{-1}
     \begin{pmatrix} I \\ I \end{pmatrix} \\
  &= P_1^{-1}
     + P_1^{-1}P_{12}\Delta_{P_1}^{-1}P_{12}^\top P_1^{-1} \\
     &\qquad\quad - P_1^{-1}P_{12}\Delta_{P_1}^{-1}
     - \Delta_{P_1}^{-1}P_{12}^\top P_1^{-1}
     + \Delta_{P_1}^{-1} \\
  &= P_1^{-1} - P_1^{-1}P_{12}\Delta_{P_1}^{-1}(I-P_{12}^\top P_1^{-1}) \\
     &\qquad\quad + \Delta_{P_1}^{-1}(I-P_{12}^\top P_1^{-1}) \\
  &= P_1^{-1} + (I-P_1^{-1}P_{12})\Delta_{P_1}^{-1}(I-P_{12}^\top P_1^{-1}).
 \end{align*}
 Define $Q:=P_1^{-1}\succ 0$, $R:=-\Delta_{P_1}\prec 0$ and $S:=(I-P_1^{-1}P_{12})$ then
 $S^\top Q^{-1}S\succeq 0$ and $SR^{-1}S^\top\preceq 0$ by Lemma~\ref{lem:sandwich}.
 This implies $\Delta_Q:=R-S^\top Q^{-1}S\prec 0$ and $\Delta_R:=Q-SR^{-1}S^\top=H^\top\Pjoint^{-1}H\succ 0$.
 Compute
 \begin{align*}
  \Delta_Q
  &= -(\Delta_{P_1}+(I-P_{12}^\top P_1^{-1})P_1(I-P_1^{-1}P_{12})) \\
  &= -(P_2-P_{12}^\top P_1^{-1}P_{12} + P_1 \\
     &\qquad\qquad -P_{12}-P_{12}^\top+P_{12}^\top P_1^{-1}P_{12}) \\
  &= -(P_1+P_2-P_{12}-P_{12}^\top)\prec 0
 \end{align*}
 and define $\Delta:=-\Delta_Q=P_1+P_2-P_{12}-P_{12}^\top\succ 0$.
 By Lemma~\ref{lem:block_inv} \textit{(4)} it follows that
 \begin{align*}
  \Delta_R^{-1}
  &= Q^{-1}+Q^{-1}S\Delta_Q^{-1}S^\top Q^{-1} \\
  &= P_1-P_1(I-P_1^{-1}P_{12})\Delta^{-1}(I-P_{12}^\top P_1^{-1})P_1 \\
  &= P_1 - (P_1-P_{12})\Delta^{-1}(P_1-P_{12}^\top)
 \end{align*}
 and hence it follows from equation~\eqref{eq:Kopt} that
 \begin{align*}
  P^\ast
  &= (H^\top\Pjoint^{-1}H)^{-1} = \Delta_R^{-1} \\
  &= P_1 - (P_1-P_{12})\Delta^{-1}(P_1-P_{12}^\top).
 \end{align*}
 As a further consequence of equation~\eqref{eq:Kopt} it follows that
 \begin{align*}
  \begin{pmatrix} K_1^\ast & K_2^\ast \end{pmatrix}
  &= (H^\top\Pjoint^{-1}H)^{-1}\begin{pmatrix} I & I \end{pmatrix}\begin{pmatrix} P_1 & P_{12} \\ P_{12}^\top & P_2 \end{pmatrix}^{-1}
 \end{align*}
 and hence that
 \begin{align*}
  &K_2^\ast\Delta_{P_1} = \\
  &(H^\top\Pjoint^{-1}H)^{-1}(-P_1^{-1}P_{12}\Delta_{P_1}^{-1} + \Delta_{P_1}^{-1})\Delta_{P_1} = \\
  &(P_1-(P_1-P_{12})\Delta^{-1}(P_1-P_{12}^\top))P_1^{-1}(P_1-P_{12}) = \\
  &((P_1-P_{12}) \\
  &\qquad -(P_1-P_{12})\Delta^{-1}(P_1-P_{12}^\top)P_1^{-1}(P_1-P_{12})) = \\
  &(P_1-P_{12})(I-\Delta^{-1}(P_1-P_{12}^\top)P_1^{-1}(P_1-P_{12})) = \\
  &(P_1-P_{12})(I-\Delta^{-1}(P_1-P_{12}-P_{12}^\top+P_{12}^\top P_1^{-1}P_{12})) = \\
  &(P_1-P_{12})(I-\Delta^{-1}(\Delta-P_2+P_{12}^\top P_1^{-1}P_{12})) = \\
  &(P_1-P_{12})(I-\Delta^{-1}(\Delta-\Delta_{P_1})) = \\
  &(P_1-P_{12})\Delta^{-1}\Delta_{P_1}.
 \end{align*}
 It follows that $K_2^\ast=(P_1-P_{12})\Delta^{-1}$.
 The formula for $K_1^\ast$ now follows from $K_1^\ast+K_2^\ast=K_1^\ast H_1+K_2^\ast H_2=K^\ast H=I$.
\qed

 \begin{table*}
 \begin{minipage}{\textwidth}
 \hrule
 \begin{equation}\label{eq:intersect_1}\begin{split}
  (\Pjoint(P_{12}))^{-1} &=
  \begin{pmatrix}
   \Delta_{P_2}^{-1} & -\Delta_{P_2}^{-1}P_{12}P_2^{-1} \\[2mm]
   -P_2^{-1}P_{12}^\top\Delta_{P_2}^{-1} &\quad P_2^{-1}P_{12}^\top\Delta_{P_2}^{-1}P_{12}P_2^{-1} + P_2^{-1}
  \end{pmatrix} \\
  &= \frac{1}{1-\lambda^2}
  \begin{pmatrix}
   P_1^{-1} & -\lambda P_1^{-\frac{1}{2}} U_1 P_2^{-\frac{1}{2}} \\[1mm]
   -\lambda P_2^{-\frac{1}{2}} U_1^\top P_1^{-\frac{1}{2}} &\quad \lambda^2 P_2^{-\frac{1}{2}}U_1^\top U_1 P_2^{-\frac{1}{2}}
  \end{pmatrix}
  +
  \begin{pmatrix}
   0 & 0 \\ 0 & P_2^{-1}
  \end{pmatrix}
 \end{split}\end{equation}
 \begin{equation}\label{eq:intersect_2}\begin{split}
  x^\top(P^\ast(P_{12}))^{-1}x &= x^\top(H^\top(\Pjoint(P_{12}))^{-1}H)x \\
  &= \frac{1}{1-\lambda^2} \cdot
     (x^\top H_1^\top P_1^{-1}H_1 x
      - \lambda x^\top H_1^\top P_1^{-\frac{1}{2}} U_1 P_2^{-\frac{1}{2}} H_2 x
      - \lambda x^\top H_2^\top P_2^{-\frac{1}{2}} U_1^\top P_1^{-\frac{1}{2}} H_1 x \\
      &\hspace*{7em} + \lambda^2 x^\top H_2^\top P_2^{-\frac{1}{2}}U_1^\top U_1 P_2^{-\frac{1}{2}} H_2 x)
     + x^\top H_2^\top P_2^{-1} H_2 x \\
  &= x^\top H_2^\top P_2^{-1} H_2 x,
 \end{split}\end{equation}
 \begin{equation}\label{eq:intersect_3}\begin{split}
  (\Pjoint(P_{12}^\epsilon))^{-1} &=
  \begin{pmatrix}
   \Delta_{P_2}^{-1} & -\Delta_{P_2}^{-1}P_{12}^\epsilon P_2^{-1} \\[2mm]
   -P_2^{-1}(P_{12}^\epsilon)^\top\Delta_{P_2}^{-1} &\quad P_2^{-1}(P_{12}^\epsilon)^\top\Delta_{P_2}^{-1}P_{12}^\epsilon P_2^{-1} + P_2^{-1}
  \end{pmatrix} \\
  &= \frac{1}{1-(1-\epsilon)^2}
  \begin{pmatrix}
   P_1^{-1} & -(1-\epsilon)P_1^{-\frac{1}{2}} U_1 P_2^{-\frac{1}{2}} \\[1mm]
   -(1-\epsilon)P_2^{-\frac{1}{2}} U_1^\top P_1^{-\frac{1}{2}} &\ (1-\epsilon)^2 P_2^{-\frac{1}{2}}U_1^\top U_1 P_2^{-\frac{1}{2}}
  \end{pmatrix}
  +
  \begin{pmatrix}
   0 & 0 \\ 0 & P_2^{-1}
  \end{pmatrix}
 \end{split}\end{equation}
 \begin{equation}\label{eq:intersect_4}\begin{split}
  x^\top(P^\ast(P_{12}^\epsilon))^{-1}x &= x^\top(H^\top(\Pjoint(P_{12}^\epsilon))^{-1}H)x \\
  &= \frac{1}{1-(1-\epsilon)^2} \cdot
     (x^\top H_1^\top P_1^{-1}H_1 x
      - (1-\epsilon)x^\top H_1^\top P_1^{-\frac{1}{2}} U_1 P_2^{-\frac{1}{2}} H_2 x
      - (1-\epsilon)x^\top H_2^\top P_2^{-\frac{1}{2}} U_1^\top P_1^{-\frac{1}{2}} H_1 x \\
      &\hspace*{9em} + (1-\epsilon)^2 x^\top H_2^\top P_2^{-\frac{1}{2}} U_1^\top U_1 P_2^{-\frac{1}{2}} H_2 x)
     + x^\top H_2^\top P_2^{-1} H_2 x \\
  &= \frac{\epsilon^2}{1-(1-\epsilon)^2} \cdot x^\top H_1^\top P_1^{-1}H_1 x
     + x^\top H_2^\top P_2^{-1} H_2 x.
 \end{split}\end{equation}
 \hrule
 \end{minipage}
 \end{table*}

\textit{Proof of Theorem~\ref{th:fill_intersection}.}
 By Lemma~\ref{lem:swap12}, we can assume without loss of generality that $p_2\geq p_1$ since otherwise we can swap and re-label $\hat{x}_1$ and $\hat{x}_2$ before proceeding.
 A point $x\in\Rn$ lies in the open set
 \begin{align*}
  \mathrm{int}(\mathcal{E}(H_1^\top P_1^{-1}H_1) \cap \mathcal{E}(H_2^\top P_2^{-1}H_2))
 \end{align*}
 if and only if $0\leq x^\top(H_1^\top P_1^{-1}H_1)x < 1$ and $0\leq x^\top(H_2^\top P_2^{-1}H_2)x < 1$.

 First, consider the case where $0\leq x^\top(H_1^\top P_1^{-1}H_1)x < x^\top(H_2^\top P_2^{-1}H_2)x < 1$.
 Then it follows that
 \begin{align*}
  0 \leq \lambda^2 := \frac{x^\top(H_1^\top P_1^{-1}H_1)x}{x^\top(H_2^\top P_2^{-1}H_2)x} < 1
 \end{align*}
 and
 \begin{align*}
  \lVert\lambda P_2^{-\frac{1}{2}}H_2 x\rVert_2^2
  &= \lambda^2\cdot x^\top(H_2^\top P_2^{-1}H_2)x \\
  &= x^\top(H_1^\top P_1^{-1}H_1)x
  = \lVert P_1^{-\frac{1}{2}}H_1 x\rVert_2^2.
 \end{align*}
 This implies that
 \begin{align*}
  \lVert\lambda P_2^{-\frac{1}{2}}H_2 x\rVert_2
  = \left\lVert\begin{pmatrix} P_1^{-\frac{1}{2}}H_1 x \\ 0 \end{pmatrix}\right\rVert_2,
 \end{align*}
 where we have used the assumption that $p_2\geq p_1$ to zero-pad the vector on the right-hand side to dimension $p_2$.
 But then there exists an orthogonal matrix
 \begin{align*}
  U \in \R^{p_2\times p_2}, \quad U =: \begin{pmatrix} U_1 \\ U_2 \end{pmatrix}
 \end{align*}
 with $U_1\in\R^{p_1\times p_2}$ such that
 \begin{equation}\label{eq:U_action}\begin{split}
  U\lambda P_2^{-\frac{1}{2}}H_2 x
  &= \begin{pmatrix} U_1 \\ U_2 \end{pmatrix} \lambda P_2^{-\frac{1}{2}}H_2 x \\
  &= \begin{pmatrix} P_1^{-\frac{1}{2}}H_1 x \\ 0 \end{pmatrix}.
 \end{split}\end{equation}
 It follows from $U^{-1}=U^\top$ that $U_1U_1^\top=I$ and from equation~\eqref{eq:U_action} that
 $U_1\lambda P_2^{-\frac{1}{2}}H_2 x=P_1^{-\frac{1}{2}}H_1 x$.

 Let $P_{12}:=\lambda P_1^{\frac{1}{2}} U_1 P_2^{\frac{1}{2}}$ then it follows that
 \begin{align*}
  \Delta_{P_2} &:= P_1-P_{12}P_2^{-1}P_{12}^\top \\
  &= P_1-\lambda^2 P_1^{\frac{1}{2}}U_1 P_2^{\frac{1}{2}}P_2^{-1}P_2^{\frac{1}{2}} U_1^\top P_1^{\frac{1}{2}} \\
  &= (1-\lambda^2)P_1
  \succ 0
 \end{align*}
 and hence that $\Pjoint(P_{12})\succ 0$ by Lemma~\ref{lem:strict_Schur_complement} \textit{(3)}.
 Moreover, it follows from Lemma~\ref{lem:block_inv} \textit{(3)} that equation~\eqref{eq:intersect_1} holds for $(\Pjoint(P_{12}))^{-1}$ and this implies equation~\eqref{eq:intersect_2},  where the last equality follows from $U_1\lambda P_2^{-\frac{1}{2}}H_2 x = P_1^{-\frac{1}{2}}H_1 x$, used twice to cancel the terms in the parentheses.

 It follows that $x^\top(P^\ast(P_{12}))^{-1} x = x^\top H_2^\top P_2^{-1} H_2 x < 1$ and hence that $x\in\mathcal{E}((P^\ast(P_{12}))^{-1})$.

 For the case $0\leq x^\top(H_2^\top\! P_2^{-1}H_2)x < x^\top(H_1^\top\! P_1^{-1}H_1)x < 1$, the proof follows along the same lines but with the roles of $P_1$ (resp. $H_1$) and $P_2$ (resp. $H_2$) reversed, and using Lemma~\ref{lem:strict_Schur_complement} \textit{(2)} resp. Lemma~\ref{lem:block_inv} \textit{(2)} instead of Lemma~\ref{lem:strict_Schur_complement} \textit{(3)} resp. Lemma~\ref{lem:block_inv} \textit{(3)}.

 In the case $0\leq x^\top(H_1^\top P_1^{-1}H_1)x=x^\top(H_2^\top P_2^{-1}H_2)x < 1$ it follows that
 $\lVert P_2^{-\frac{1}{2}}H_2 x\rVert_2^2=x^\top(H_2^\top P_2^{-1}H_2)x=x^\top(H_1^\top P_1^{-1}H_1)x=\lVert P_1^{-\frac{1}{2}}H_1 x\rVert_2^2$ and, using the same reasoning as in the first case, this implies that there exists a matrix $U\in\R^{p_1\times p_2}$ such that $U_1U_1^\top=I$ and $U_1P_2^{-\frac{1}{2}}H_2 x = P_1^{-\frac{1}{2}}H_1 x$.
 For $0<\epsilon\leq 1$ let $P_{12}^\epsilon:=(1-\epsilon)P_1^{\frac{1}{2}} U_1 P_2^{\frac{1}{2}}$ then it follows that
 $\Delta_{P_2}:=P_1-P_{12}^\epsilon P_2^{-1}(P_{12}^\epsilon)^\top=P_1-(1-\epsilon)^2 P_1^{\frac{1}{2}}U_1 P_2^{\frac{1}{2}}P_2^{-1}P_2^{\frac{1}{2}} U_1^\top P_1^{\frac{1}{2}}=(1-(1-\epsilon)^2)P_1\succ 0$ and hence that
 $\Pjoint(P_{12}^\epsilon)\succ 0$ by Lemma~\ref{lem:strict_Schur_complement} \textit{(3)}.
 Moreover, it follows from Lemma~\ref{lem:block_inv} \textit{(3)} that equation~\eqref{eq:intersect_3} holds for $(\Pjoint(P_{12}^\epsilon))^{-1}$ and this implies equation~\eqref{eq:intersect_4}, where the last equality follows from repeated use of
 $U_1P_2^{-\frac{1}{2}}H_2 x = P_1^{-\frac{1}{2}}H_1 x$ to combine the terms in the parentheses.

 Since $x^\top H_2^\top P_2^{-1} H_2 x < 1$ and
 \begin{align*}
  \lim_{\epsilon\to 0} \frac{\epsilon^2}{1-(1-\epsilon)^2} = \lim_{\epsilon\to 0} \frac{\epsilon^2}{\epsilon(2-\epsilon)} = 0,
 \end{align*}
 there exists $\epsilon_0\in(0,1]$ such that $x^\top(P^\ast(P_{12}^{\epsilon_0}))^{-1}x < 1$ and hence $x\in\mathcal{E}((P^\ast(P_{12}^{\epsilon_0}))^{-1})$.
 This completes the proof.
\qed

\textit{Proof of Theorem~\ref{th:kahan_parameter}.}
 Assume that $\balpha=\balpha'$ then $(\bK_1^\CIa,\bK_2^\CIa,\bPhat^\CIa)=(\bK_1^{\CIa'},\bK_2^{\CIa'},\bPhat^{\CIa'})$ by definition.

 Conversely, let $(\bK_1^\CIa\!,\bK_2^\CIa,\bPhat^\CIa)=(\bK_1^{\CIa'}\!\!,\bK_2^{\CIa'}\!\!,\bPhat^{\CIa'})$.
 For the case $H_2^\top\hat{P}_2^{-1}H_2=H_1^\top\hat{P}_1^{-1}H_1$ it follows from Lemma~\ref{lem:corner_cases} that $p_1=n$ and $H_1$ is non-singular and hence that
 \begin{align*}
  \alpha I
  &= K_1^\alpha\hat{P}_1(H_1^\top)^{-1}(\hat{P}^\alpha)^{-1} \\
  &= K_1^{\alpha'}\hat{P}_1(H_1^\top)^{-1}(\hat{P}^{\alpha'})^{-1}
  = \alpha' I.
 \end{align*}
 It follows that $\alpha=\alpha'$ in this case.
 Now assume that $H_2^\top\hat{P}_2^{-1}H_2\not=H_1^\top\hat{P}_1^{-1}H_1$ then there exist $i,j\in\{1,\dots,n\}$ such that
 \begin{align*}
  z_2:=(H_2^\top\hat{P}_2^{-1}H_2)_{ij}\not=(H_1^\top\hat{P}_1^{-1}H_1)_{ij}=:z_1
 \end{align*}
 and $\alpha z_1+(1-\alpha) z_2=((\hat{P}^\alpha)^{-1})_{ij}=((\hat{P}^{\alpha'})^{-1})_{ij}=\alpha' z_1+(1-\alpha')z_2$.
 It follows that $\alpha(z_1-z_2)=\alpha'(z_1-z_2)$ and hence that $\alpha=\alpha'$ also in this case.
\qed

\textit{Proof of Theorem~\ref{th:alpha}.}
 It follows from repeated application of Lemma~\ref{lem:Schur_complement} \textit{(3)} that $\balpha' I\succeq 0$ and $(1-\balpha')I\succeq 0$ and hence that $\balpha'\in[0,1]$.
 The proof proceeds case by case.

 Case 1: $\alpha'=0$.
 In this case it follows from Lemma~\ref{lem:Schur_complement} \textit{(3)} that
 \begin{align*}
  0 &=
  \begin{pmatrix} K_1^{\alpha}\hat{P}_1^{\frac{1}{2}} & K_2^{\alpha}\hat{P}_2^{\frac{1}{2}} \end{pmatrix}
  \left(\begin{pmatrix} I & 0 \\ 0 & I \end{pmatrix} -
   \begin{pmatrix} 0 & 0 \\ 0 & I \end{pmatrix}\begin{pmatrix} 0 & 0 \\ 0 & I \end{pmatrix}^\dagger\right) \\
  &= K_1^{\alpha}\hat{P}_1^{\frac{1}{2}}
  = \alpha\hat{P}^\alpha H_1^\top\hat{P}_1^{-\frac{1}{2}}.
 \end{align*}
 Since $\hat{P}^\alpha\succ 0$ this implies that $\alpha H_1^\top=0$ and hence that $\alpha=0$ or $H_1=0$.
 In the latter case it follows from Assumption~\eqref{eq:A1} that $p_2=n$ and $H_2$ is non-singular square.
 But then it follows from Lemma~\ref{lem:sandwich} that $H_2^\top\hat{P}_2^{-1}H_2\succ 0=H_1^\top\hat{P}_1^{-1}H_1$ and hence that $\alpha=0$ also in this case.
 In summary, it follows that $\alpha=0=\alpha'$.

 Case 2: $\alpha'=1$.
 In this case it follows from Lemma~\ref{lem:Schur_complement} \textit{(3)} that
 \begin{align*}
  0 &=
  \begin{pmatrix} K_1^{\alpha}\hat{P}_1^{\frac{1}{2}} & K_2^{\alpha}\hat{P}_2^{\frac{1}{2}} \end{pmatrix}
  \left(\begin{pmatrix} I & 0 \\ 0 & I \end{pmatrix} -
   \begin{pmatrix} I & 0 \\ 0 & 0 \end{pmatrix}\begin{pmatrix} I & 0 \\ 0 & 0 \end{pmatrix}^\dagger\right) \\
  &= K_2^{\alpha}\hat{P}_2^{\frac{1}{2}}
  = (1-\alpha)\hat{P}^\alpha H_2^\top\hat{P}_2^{-\frac{1}{2}}.
 \end{align*}
 Since $\hat{P}^\alpha\succ 0$ this implies that $(1-\alpha)H_2^\top=0$ and hence that $\alpha=1$ or $H_2=0$.
 In the latter case it follows from Assumption~\eqref{eq:A1} that $p_1=n$ and $H_1$ is non-singular square.
 But then it follows from Lemma~\ref{lem:sandwich} that $H_2^\top\hat{P}_2^{-1}H_2=0\prec H_1^\top\hat{P}_1^{-1}H_1$ and hence that $\alpha=1$ also in this case.
 In summary, it follows that $\alpha=1=\alpha'$.

 Case 3: $\alpha'\in(0,1)$.
 In this case it follows from Lemma~\ref{lem:Schur_complement} \textit{(3)} that
 \begin{align*}
  &\hat{P}^\alpha -
  \begin{pmatrix} K_1^{\alpha}\hat{P}_1^{\frac{1}{2}} & K_2^{\alpha}\hat{P}_2^{\frac{1}{2}} \end{pmatrix}
  \begin{pmatrix} \alpha' I & 0 \\ 0 & (1-\alpha')I \end{pmatrix}^\dagger
  \!\!\begin{pmatrix} \hat{P}_1^{\frac{1}{2}}(K_1^{\alpha})^\top \\ \hat{P}_2^{\frac{1}{2}}(K_2^{\alpha})^\top \end{pmatrix} \\
  &\quad \succeq 0
 \end{align*}
 and hence that
 \begin{align*}
  \hat{P}^\alpha &\succeq
  \frac{1}{\alpha'}K_1^\alpha\hat{P}_1(K_1^\alpha)^\top + \frac{1}{1-\alpha'}K_2^\alpha\hat{P}_2(K_2^\alpha)^\top \\
  &= \frac{\alpha^2}{\alpha'}\hat{P}^\alpha H_1^\top\hat{P}_1^{-1}H_1\hat{P}^\alpha
     + \frac{(1-\alpha)^2}{1-\alpha'}\hat{P}^\alpha H_2^\top\hat{P}_2^{-1}H_2\hat{P}^\alpha\!\!.
 \end{align*}
 By Lemma~\ref{lem:sandwich} this implies that
 \begin{multline*}
  (\hat{P}^\alpha)^{-1} =
  \alpha H_1^\top\hat{P}_1^{-1}H_1 + (1-\alpha) H_2^\top\hat{P}_2^{-1}H_2 \\
  \succeq
  \frac{\alpha^2}{\alpha'}H_1^\top\hat{P}_1^{-1}H_1 + \frac{(1-\alpha)^2}{1-\alpha'}H_2^\top\hat{P}_2^{-1}H_2.
 \end{multline*}
 Subtracting the right-hand side of this inequality from the left-hand side and applying algebraic manipulations then yields
 \begin{equation}\label{eq:Delta_alpha}
  \left(\frac{\alpha}{\alpha'}H_1^\top\hat{P}_1^{-1}H_1 - \frac{1-\alpha}{1-\alpha'} H_2^\top\hat{P}_2^{-1}H_2\right) \Delta_\alpha
  \succeq 0,
 \end{equation}
 where $\Delta_\alpha:=\alpha'-\alpha$.
 Assume, to arrive at a contradiction, that $\Delta_\alpha\not=0$.

 Case 3a: $\Delta_\alpha<0$.
 In this case it follows that
 \begin{align*}
  H_1^\top\hat{P}_1^{-1}H_1 &\prec \frac{\alpha}{\alpha'}H_1^\top\hat{P}_1^{-1}H_1
  \preceq \frac{1-\alpha}{1-\alpha'} H_2^\top\hat{P}_2^{-1}H_2 \\
  &\prec H_2^\top\hat{P}_2^{-1}H_2,
 \end{align*}
 where the first and last inequality follow from $\alpha'<\alpha$.
 It follows that $\alpha=0$ and hence that $\alpha'<0$, a contradiction.

 Case 3b: $\Delta_\alpha>0$.
 In this case it follows that
 \begin{align*}
  H_1^\top\hat{P}_1^{-1}H_1 &\succ \frac{\alpha}{\alpha'}H_1^\top\hat{P}_1^{-1}H_1
  \succeq \frac{1-\alpha}{1-\alpha'} H_2^\top\hat{P}_2^{-1}H_2 \\
  &\succ H_2^\top\hat{P}_2^{-1}H_2,
 \end{align*}
 where the first and last inequality follow from $\alpha'>\alpha$.
 It follows that $\alpha=1$ and hence that $\alpha'>1$, a contradiction.

 In summary, it follows that $\Delta_\alpha=0$ and hence that $\alpha=\alpha'$.
\qed

\section{Proofs of main results}\label{sec:proofs_main}

\textit{Proof of Theorem~\ref{th:conservative_3}.}
 First consider the case $Q_1=0$.
 In this case, condition~\eqref{eq:conservative_2} is equivalent to $Q_2Q_2^\top\preceq\hat{P}$ which is condition~\eqref{eq:conservative_4} in this case.
 Similarly, in the case $Q_2=0$, condition~\eqref{eq:conservative_2} is equivalent to $Q_1Q_1^\top\preceq\hat{P}$ which is condition~\eqref{eq:conservative_4} in that case.
 In both these cases, the equivalence holds irrespective of the choice of $\Domcons$.
 In the case $Q_1\not=0$ and $Q_2\not=0$, condition~\eqref{eq:conservative_2} is equivalent to
 $Q_1Q_1^\top+Q_1XQ_2^\top+Q_2XQ_1^\top+Q_2Q_2^\top\preceq\hat{P}$.
 Setting $G:=-\hat{P}+Q_1Q_1^\top+Q_2Q_2^\top$, $M:=Q_1$, $N:=Q_2^\top$ and $Y:=X$ it follows that $G=G^\top$, $M\not=0$ and $N\not=0$ and that condition~\eqref{eq:conservative_2} is equivalent to $G+MYN+N^\top Y^\top M^\top \preceq 0$.
 But then Lemma~\ref{lem:Petersen} implies that condition~\eqref{eq:conservative_2} is equivalent to the existence of $\epsilon>0$ such that $-\hat{P}+Q_1Q_1^\top+Q_2Q_2^\top+\epsilon Q_1Q_1^\top+\frac{1}{\epsilon}Q_2Q_2^\top\preceq 0$, again irrespective of the choice of $\Domcons$.
 Setting $\tau:=\epsilon$ it follows that condition~\eqref{eq:conservative_2} is equivalent to condition~\eqref{eq:conservative_4} also in this case.
 In summary, it follows that \textit{(1)}, \textit{(2)} and \textit{(3)} are equivalent.

 Setting $\alpha:=\frac{1}{1+\tau}$ for $\tau\not=-1$ it follows that $\tau>0$ if and only $\alpha\in(0,1)$. Moreover,
 $(1+\tau)=\frac{1}{\alpha}$, $\tau=\frac{1}{\alpha}-1$ and
 $(1+\frac{1}{\tau})=1+\frac{1}{\frac{1}{\alpha}-1}=\frac{1-\alpha}{1-\alpha}+\frac{\alpha}{1-\alpha}=\frac{1}{1-\alpha}$.
 The equivalence of \textit{(3)} and \textit{(4)} now follows immediately.

 Setting
 \begin{align*}
  Q:=\hat{P},
  R:=\begin{pmatrix} \alpha I & 0 \\ 0 & (1-\alpha)I \end{pmatrix} \text{ and }
  S:=\begin{pmatrix} Q_1 & Q_2 \end{pmatrix}
 \end{align*}
 it follows that $R\succeq 0$,
 \begin{align*}
  R^+ =
  \begin{cases}
   \startM 0 & 0 \\ 0 & I \stopM, & \alpha=0, \\
   \startM I & 0 \\ 0 & 0 \stopM, & \alpha=1, \\
   \startM \frac{1}{\alpha}I & 0 \\ 0 & \frac{1}{1-\alpha}I \stopM, & \alpha\in(0,1),
  \end{cases}
 \end{align*}
 and
 \begin{align*}
  S(I-RR^+) =
  \begin{cases}
   Q_1 , & \alpha=0, \\
   Q_2 , & \alpha=1, \\
   0 , & \alpha\in(0,1),
  \end{cases}
 \end{align*}
 and
 \begin{align*}
  Q-SR^+S^\top =
  \begin{cases}
   \hat{P}-Q_2Q_2^\top, & \!\!\!\!\alpha=0, \\
   \hat{P}-Q_1Q_1^\top, & \!\!\!\!\alpha=1, \\
   \hat{P}-\frac{1}{\alpha}Q_1Q_1^\top-\frac{1}{1-\alpha}Q_2Q_2^\top, & \!\!\!\!\alpha\in(0,1).
  \end{cases}
 \end{align*}

 Lemma~\ref{lem:Schur_complement} \textit{(3)} now implies that condition~\eqref{eq:conservative_5.3} is equivalent to condition~\eqref{eq:conservative_5.2}, noting that $\alpha=0$ (resp. $\alpha=1$) in condition~\eqref{eq:conservative_5.3} is equivalent to $Q_1=0$ (resp. $Q_2=0$) in condition~\eqref{eq:conservative_5.2}.
 This implies the equivalence of \textit{(4)} and \textit{(5)}.
\qed

\textit{Proof of Theorem~\ref{th:conservative_fusion}.}
 Let $(\bK_1,\bK_2,\bPhat)$ be unbiased and preserve conservativeness on $\PPpi$.
 By Proposition~\ref{prop:conservative} this implies that
 \begin{align*}
  \bPhat \succeq \begin{pmatrix} \bK_1 & \bK_2 \end{pmatrix}
                 \begin{pmatrix} \hat{P}_1 & \hat{P}_{12} \\ \hat{P}_{12}^\top & \hat{P}_1 \end{pmatrix}
                 \begin{pmatrix} \bK_1^\top \\ \bK_2^\top \end{pmatrix}
 \end{align*}
 for all $\hat{P}_{12}\in\R^{p_1\times p_2}$ such that
 \begin{align*}
  \hat{P}_{\rm joint} := \begin{pmatrix} \hat{P}_1 & \hat{P}_{12} \\ \hat{P}_{12}^\top & \hat{P}_1 \end{pmatrix} \in \PPpi.
 \end{align*}
 Note that by our notational convention, the variables $\hat{P}_1$ and $\hat{P}_2$ on the right hand-side of this inequality are (part of) the arguments of the function $\bPhat$ on the left.

 By Theorem~\ref{th:Kopt}\footnote{Note that while Theorem~\ref{th:Kopt} and Theorem~\ref{th:fill_intersection} are written in terms of generic variables $P_1\succ 0$ and $P_2\succ 0$, they get applied here with the specific choice $P_1:=\hat{P}_1$ and $P_2:=\hat{P}_2$.} this implies that $\bPhat\succeq (H^\top\hat{P}_{\rm joint}^{-1}H)^{-1}$ for all
 $\hat{P}_{\rm joint}\in\PPpi$.
 It follows from \cite[Proposition~8.5.5]{bernstein} that $\bPhat^{-1}\preceq H^\top\hat{P}_{\rm joint}^{-1}H$ for all
 $\hat{P}_{\rm joint}\in\PPpi$.
 Using the notation of Theorem~\ref{th:fill_intersection}, it hence follows from Lemma~\ref{lem:kahan} that
 $\mathcal{E}(\bPhat^{-1})\supseteq\mathcal{E}((P^\ast(\hat{P}_{12}))^{-1})$ for all $\hat{P}_{12}\in\mathcal{P}_{12}^{\succ}$.
 But then it follows from Theorem~\ref{th:fill_intersection} that
 $\mathcal{E}(\bPhat^{-1})\supseteq\mathrm{int}(\mathcal{E}(H_1^\top\hat{P}_1^{-1}H_1) \cap \mathcal{E}(H_2^\top\hat{P}_2^{-1}H_2))$ and since both $\mathcal{E}(\bPhat^{-1})$ and $\mathcal{E}(H_1^\top\hat{P}_1^{-1}H_1) \cap \mathcal{E}(H_2^\top\hat{P}_2^{-1}H_2))$ are closed sets, this implies that
 \begin{align*}
  \mathcal{E}(\bPhat^{-1})\supseteq\mathcal{E}(H_1^\top\hat{P}_1^{-1}H_1) \cap \mathcal{E}(H_2^\top\hat{P}_2^{-1}H_2).
 \end{align*}
 This completes the proof.
\qed

\textit{Proof of Theorem~\ref{th:lower_bound}.}
 In this proof we will make frequent tacit use of Lemma~\ref{lem:kahan} to switch between inequalities in the L{\"o}wner order and inclusion relations between the corresponding error ellipsoids.
 By Theorem~\ref{th:conservative_fusion},
 \begin{equation}\label{eq:contains_intersection}
  \mathcal{E}(\bPhat^{-1})\supseteq\mathcal{E}(H_1^\top\hat{P}_1^{-1}H_1) \cap \mathcal{E}(H_2^\top\hat{P}_2^{-1}H_2).
 \end{equation}

 Consider the case where $H_2^\top\hat{P}_2^{-1}H_2 \succeq H_1^\top\hat{P}_1^{-1}H_1$ and hence
 $\mathcal{E}(H_2^\top\hat{P}_2^{-1}H_2)\subseteq\mathcal{E}(H_1^\top\hat{P}_1^{-1}H_1)$.
 By equation~\eqref{eq:contains_intersection} this implies that $\mathcal{E}(\hat{P}^{-1})\supseteq\mathcal{E}(H_2^\top\hat{P}_2^{-1}H_2)$ and hence that $0\prec\hat{P}^{-1}\preceq H_2^\top\hat{P}_2^{-1}H_2$.
 It follows from transitivity of the L{\"o}wner order that $H_2^\top\hat{P}_2^{-1}H_2\succ 0$ and hence
 it follows from \cite[Proposition~8.5.5]{bernstein} that $\hat{P}\succeq(H_2^\top\hat{P}_2^{-1}H_2)^{-1}$.

 Consider the case where $H_1^\top\hat{P}_1^{-1}H_1 \succeq H_2^\top\hat{P}_2^{-1}H_2$ and hence
 $\mathcal{E}(H_1^\top\hat{P}_1^{-1}H_1)\subseteq\mathcal{E}(H_2^\top\hat{P}_2^{-1}H_2)$.
 By equation~\eqref{eq:contains_intersection} this implies that $\mathcal{E}(\hat{P}^{-1})\supseteq\mathcal{E}(H_1^\top\hat{P}_1^{-1}H_1)$ and hence that $0\prec\hat{P}^{-1}\preceq H_1^\top\hat{P}_1^{-1}H_1$.
 It follows from transitivity of the L{\"o}wner order that $H_1^\top\hat{P}_1^{-1}H_1\succ 0$ and hence
 it follows from \cite[Proposition~8.5.5]{bernstein} that $\hat{P}\succeq(H_1^\top\hat{P}_1^{-1}H_1)^{-1}$.

 Consider the case where $H_2^\top\hat{P}_2^{-1}H_2 \not\succeq H_1^\top\hat{P}_1^{-1}H_1$ and
 $H_1^\top\hat{P}_1^{-1}H_1 \not\succeq H_2^\top\hat{P}_2^{-1}H_2$.
 In this case it follows from Theorem~\ref{th:kahan} applied with $\Sigma:=\hat{P}^{-1}$, $\Sigma_1:=H_1^\top\hat{P}_1^{-1}H_1$ and $\Sigma_2:=H_2^\top\hat{P}_2^{-1}H_2$ that there exists $\alpha\in[0,1]$ such that
 \begin{align*}
  \mathcal{E}(\hat{P}^{-1})\supseteq\mathcal{E}(\alpha H_1^\top\hat{P}_1^{-1}H_1 + (1-\alpha)H_2^\top\hat{P}_2^{-1}H_2)
 \end{align*}
 and hence $0\prec\hat{P}^{-1}\preceq\alpha H_1^\top\hat{P}_1^{-1}H_1 + (1-\alpha)H_2^\top\hat{P}_2^{-1}H_2$.
 It follows from transitivity of the L{\"o}wner order that $\alpha H_1^\top\hat{P}_1^{-1}H_1 + (1-\alpha)H_2^\top\hat{P}_2^{-1}H_2\succ 0$
 and hence it follows from \cite[Proposition~8.5.5]{bernstein} that
 $\hat{P}\succeq(\alpha H_1^\top\hat{P}_1^{-1}H_1 + (1-\alpha)H_2^\top\hat{P}_2^{-1}H_2)^{-1}$.
 This completes the proof.
\qed

\textit{Proof of Theorem~\ref{th:kahan_family}.}
 It follows from Lemma~\ref{lem:corner_cases} that $\mathcal{K}_H^{\balpha}\not=\emptyset$.
 Let $\balpha\in\mathcal{K}_H^{\balpha}$ then it follows from Lemma~\ref{lem:corner_cases} that the linear fusion rule $(\bK_1^\CIa,\bK_2^\CIa,\bPhat^\CIa)$ is well-defined.
 It is unbiased since
 \begin{align*}
  &\bK_1^\CIa H_1+\bK_2^\CIa H_2 \\
  &\quad = \balpha\bPhat^\CIa H_1^\top\hat{P}_1^{-1}H_1 + (1-\balpha)\bPhat^\CIa H_2^\top\hat{P}_2^{-1}H_2 \\
  &\quad = \bPhat^\CIa(\balpha H_1^\top\hat{P}_1^{-1}H_1 + (1-\balpha)H_2^\top\hat{P}_2^{-1}H_2) = I.
 \end{align*}
 Let $\bQ_1^\CIa:=\bK_1^\CIa\hat{P}_1^{\frac{1}{2}}$ and $\bQ_2^\CIa:=\bK_2^\CIa\hat{P}_2^{\frac{1}{2}}$ then it follows that
 \begin{align*}
  \bQ_1^\CIa(\bQ_1^\CIa)^\top &= \balpha^2\bPhat^\CIa H_1^\top\hat{P}_1^{-1}H_1\bPhat^\CIa \quad \text{and} \\
  \bQ_2^\CIa(\bQ_2^\CIa)^\top &= (1-\balpha)^2\bPhat^\CIa H_2^\top\hat{P}_2^{-1}H_2\bPhat^\CIa.
 \end{align*}

 Consider the case $Q_1^\alpha=0$ which implies $\alpha=0$ or $H_1=0$.
 It follows from Lemma~\ref{lem:corner_cases} \textit{(1)} and Lemma~\ref{lem:sandwich} that $H_1=0$ implies
 $H_2^\top\hat{P}_2^{-1}H_2 \succ 0 = H_1^\top\hat{P}_1^{-1}H_1$ and hence also that $\alpha=0$.
 It follows that $\hat{P}^\alpha = (H_2^\top\hat{P}_2^{-1}H_2)^{-1}$ and
 $Q_2^\alpha(Q_2^\alpha)^\top=\hat{P}^\alpha H_2^\top\hat{P}_2^{-1}H_2\hat{P}^\alpha=\hat{P}^\alpha$.
 This implies that $\hat{P}^\alpha\succeq Q_2^\alpha(Q_2^\alpha)^\top$ in this case.

 Consider the case $Q_2^\alpha=0$ which implies $\alpha=1$ or $H_2=0$.
 It follows from Lemma~\ref{lem:corner_cases} \textit{(2)} and Lemma~\ref{lem:sandwich} that $H_2=0$ implies
 $H_1^\top\hat{P}_1^{-1}H_1 \succ 0 = H_2^\top\hat{P}_2^{-1}H_2$ and hence also that $\alpha=1$.
 This implies that $\hat{P}^\alpha = (H_1^\top\hat{P}_1^{-1}H_1)^{-1}$ and
 $Q_1^\alpha(Q_1^\alpha)^\top=\hat{P}^\alpha H_1^\top\hat{P}_1^{-1}H_1\hat{P}^\alpha=\hat{P}^\alpha$.
 This implies that $\hat{P}^\alpha\succeq Q_1^\alpha(Q_1^\alpha)^\top$ in this case.

 Consider the case where $Q_1^\alpha\not=0$ and $Q_2^\alpha\not=0$.
 This is the case if and only if $\alpha\in(0,1)$ and $H_1\not=0$ and $H_2\not=0$, in which case
 \begin{align*}
  &\frac{1}{\alpha}Q_1^\alpha(Q_1^\alpha)^\top+\frac{1}{1-\alpha}Q_2^\alpha(Q_2^\alpha)^\top \\
  &\quad = \frac{1}{\alpha}\alpha^2\hat{P}^\alpha H_1^\top\hat{P}_1^{-1}H_1\hat{P}^\alpha \\
  &\quad\ \qquad\qquad  + \frac{1}{1-\alpha}(1-\alpha)^2\hat{P}^\alpha H_2^\top\hat{P}_2^{-1}H_2\hat{P}^\alpha \\
  &\quad = \hat{P}^\alpha(\alpha H_1^\top\hat{P}_1^{-1}H_1 + (1-\alpha) H_2^\top\hat{P}_2^{-1}H_2)\hat{P}^\alpha \\
  &\quad = \hat{P}^\alpha.
 \end{align*}
 This implies that $\hat{P}^\alpha\succeq\frac{1}{\alpha}Q_1^\alpha(Q_1^\alpha)^\top+\frac{1}{1-\alpha}Q_2^\alpha(Q_2^\alpha)^\top$ in this case.

 It now follows from Theorem~\ref{th:conservative_3} that $\bPhat^\CIa$ preserves conservativeness on $\PPpi$.
 This implies that the linear fusion rule $(\bK_1^\CIa,\bK_2^\CIa,\bPhat^\CIa)$ preserves conservativeness on $\PPpi$
 and hence that $(\bK_1^\CIa,\bK_2^\CIa,\bPhat^\CIa)\in\mathcal{F}_H$.
\qed

\textit{Proof of Theorem~\ref{th:Kahan_implies_opt}.}
 Let $(\bK_1^\ast,\bK_2^\ast,\bPhat^\ast)\in\mathcal{K}_H\subseteq\mathcal{F}_H$ be $J$-optimal over $\mathcal{K}_H$ and let $(\bK_1,\bK_2,\bPhat)\in\mathcal{F}_H$.
 It follows from Theorem~\ref{th:lower_bound} that there exists $(\bK_1^\CIa,\bK_2^\CIa,\bPhat^\CIa)\in\mathcal{K}_H$ such that $\bPhat\succeq\bPhat^\CIa$.
 But then
 \begin{align*}
  J(\bPhat)\geq J(\bPhat^\CIa)\geq J(\bPhat^\ast),
 \end{align*}
 where the first inequality follows from isotony of $J$ and the second inequality follows from $J$-optimality of $\bPhat^\ast$ over $\mathcal{K}_H$.
 It follows that $(\bK_1^\ast,\bK_2^\ast,\bPhat^\ast)$ is $J$-optimal.
\qed

\textit{Proof of Theorem~\ref{th:opt_implies_Kahan}.}
 Let $(\bK_1^\ast,\bK_2^\ast,\bPhat^\ast)\in\mathcal{F}_H$ be $J$-optimal and assume, to arrive at a contradiction, that $(\bK_1^\ast,\bK_2^\ast,\bPhat^\ast)\not\in\mathcal{K}_H$.
 Then there exists $(\hat{x}_1,\hat{x}_2,\hat{P}_1,\hat{P}_2)\in\Rpo\times\Rpt\times\PPpo\times\PPpt$ such that one of the following three statements is true.
 \begin{enumerate}
  \item $H_2^\top\hat{P}_2^{-1}H_2 \succeq H_1^\top\hat{P}_1^{-1}H_1\ \&\ \hat{P}^\ast\not=(H_2^\top\hat{P}_2^{-1}H_2)^{-1}$.
  \item $H_1^\top\hat{P}_1^{-1}H_1 \succeq H_2^\top\hat{P}_2^{-1}H_2\ \&\ \hat{P}^\ast\not=(H_1^\top\hat{P}_1^{-1}H_1)^{-1}$.
  \item Neither $H_2^\top\hat{P}_2^{-1}H_2 \succeq H_1^\top\hat{P}_1^{-1}H_1$ nor
   $H_1^\top\hat{P}_1^{-1}H_1 \succeq H_2^\top\hat{P}_2^{-1}H_2$, and
   \begin{align*}
    \hat{P}^\ast\not=(\alpha H_1^\top\hat{P}_1^{-1}H_1 + (1-\alpha)H_2^\top\hat{P}_2^{-1}H_2)^{-1}
   \end{align*}
   for all $\alpha\in[0,1]$ such that
   \begin{align*}
    \alpha H_1^\top\hat{P}_1^{-1}H_1 + (1-\alpha)H_2^\top\hat{P}_2^{-1}H_2
   \end{align*}
   is non-singular.
 \end{enumerate}

 Consider the case where statement \textit{(1)} is true.
 It follows from Theorem~\ref{th:lower_bound} that $\hat{P}^\ast\succeq(H_2^\top\hat{P}_2^{-1}H_2)^{-1}$.
 Since $\hat{P}^\ast\not=(H_2^\top\hat{P}_2^{-1}H_2)^{-1}$ and $J$ is strictly isotone, it follows that \begin{align*}
  J(\hat{P}^\ast)>J((H_2^\top\hat{P}_2^{-1}H_2)^{-1}).
 \end{align*}
 However, $\hat{P}^\alpha=\bPhat^\CIa(\hat{x}_1,\hat{x}_2,\hat{P}_1,\hat{P}_2)=(H_2^\top\hat{P}_2^{-1}H_2)^{-1}$ for all $(\bK_1^\CIa,\bK_2^\CIa,\bPhat^\CIa)\in\mathcal{K}_H\subseteq\mathcal{F}_H$, and hence the $J$-optimality of $(\bK_1^\ast,\bK_2^\ast,\bPhat^\ast)$ implies that $J(\hat{P}^\ast)\leq J(\hat{P}^\alpha)=J((H_2^\top\hat{P}_2^{-1}H_2)^{-1})$, a contradiction.

 Consider the case where statement \textit{(2)} is true.
 It follows from Theorem~\ref{th:lower_bound} that $\hat{P}^\ast\succeq(H_1^\top\hat{P}_1^{-1}H_1)^{-1}$.
 Since $\hat{P}^\ast\not=(H_1^\top\hat{P}_1^{-1}H_1)^{-1}$ and $J$ is strictly isotone, it follows that
 \begin{align*}
  J(\hat{P}^\ast)>J((H_1^\top\hat{P}_1^{-1}H_1)^{-1}).
 \end{align*}
 However, $\hat{P}^\alpha=\bPhat^\CIa(\hat{x}_1,\hat{x}_2,\hat{P}_1,\hat{P}_2)=(H_1^\top\hat{P}_1^{-1}H_1)^{-1}$ for all $(\bK_1^\CIa,\bK_2^\CIa,\bPhat^\CIa)\in\mathcal{K}_H\subseteq\mathcal{F}_H$, and hence the $J$-optimality of $(\bK_1^\ast,\bK_2^\ast,\bPhat^\ast)$ implies that $J(\hat{P}^\ast)\leq J(\hat{P}^\alpha)=J((H_1^\top\hat{P}_1^{-1}H_1)^{-1})$, a contradiction.

 Consider the case where statement \textit{(3)} is true.
 It follows from Theorem~\ref{th:lower_bound} that there exists $\alpha\in[0,1]$ such that
 $\alpha H_1^\top\hat{P}_1^{-1}H_1 + (1-\alpha)H_2^\top\hat{P}_2^{-1}H_2$ is non-singular and
 $\hat{P}^\ast\succeq(\alpha H_1^\top\hat{P}_1^{-1}H_1 + (1-\alpha)H_2^\top\hat{P}_2^{-1}H_2)^{-1}$.
 Since $\hat{P}^\ast\not=(\alpha H_1^\top\hat{P}_1^{-1}H_1 + (1-\alpha)H_2^\top\hat{P}_2^{-1}H_2)^{-1}$ and $J$ is strictly isotone, it follows that
 \begin{align*}
  J(\hat{P}^\ast)>J((\alpha H_1^\top\hat{P}_1^{-1}H_1 + (1-\alpha)H_2^\top\hat{P}_2^{-1}H_2)^{-1}).
 \end{align*}
 However, there exists a function $\balpha\in\mathcal{K}_H^{\balpha}$ such that $\balpha(\hat{x}_1,\hat{x}_2,\hat{P}_1,\hat{P}_2)=\alpha$.
 It follows that $(\bK_1^\CIa,\bK_2^\CIa,\bPhat^\CIa)\in\mathcal{K}_H\subseteq\mathcal{F}_H$ and hence that $\hat{P}^\alpha=\bPhat^\CIa(\hat{x}_1,\hat{x}_2,\hat{P}_1,\hat{P}_2)=(\alpha H_1^\top\hat{P}_1^{-1}H_1 + (1-\alpha)H_2^\top\hat{P}_2^{-1}H_2)^{-1}$.
 It hence follows from $J$-optimality of $(\bK_1^\ast,\bK_2^\ast,\bPhat^\ast)$ that
 $J(\hat{P}^\ast)\leq J((\alpha H_1^\top\hat{P}_1^{-1}H_1 + (1-\alpha)H_2^\top\hat{P}_2^{-1}H_2)^{-1})$, a contradiction.

 In summary, it follows that $(\bK_1^\ast,\bK_2^\ast,\bPhat^\ast)\in\mathcal{K}_H$ and since it is $J$-optimal, it is in particular $J$-optimal over $\mathcal{K}_H$.
\qed

\textit{Proof of Theorem~\ref{th:CI_implies_Kahan}.}
 The proof proceeds case-by-case.
 Consider the case where $H_2^\top\hat{P}_2^{-1}H_2 \succ H_1^\top\hat{P}_1^{-1}H_1$.
 It then follows from Lemma~\ref{lem:corner_cases} \textit{(1)} and Lemma~\ref{lem:sandwich} that $H_2^\top\hat{P}_2^{-1}H_2 \succ 0$ and hence that
 \begin{align*}
  \alpha' H_1^\top\hat{P}_1^{-1}H_1 + (1-\alpha')H_2^\top\hat{P}_2^{-1}H_2 \succ 0.
 \end{align*}
 We can hence define
 $\hat{P}^{\alpha'}:=(\alpha' H_1^\top\hat{P}_1^{-1}H_1 + (1-\alpha')H_2^\top\hat{P}_2^{-1}H_2)^{-1}\in\PPn$.
 Assume, to arrive at a contradiction, that $\alpha'>0$.
 It follows that
 \begin{align*}
  (\hat{P}^{\alpha'})^{-1}
  &= \alpha' H_1^\top\hat{P}_1^{-1}H_1 + (1-\alpha')H_2^\top\hat{P}_2^{-1}H_2 \\
  &= H_2^\top\hat{P}_2^{-1}H_2 - \alpha'(H_2^\top\hat{P}_2^{-1}H_2-H_1^\top\hat{P}_1^{-1}H_1) \\
  &\prec H_2^\top\hat{P}_2^{-1}H_2.
 \end{align*}
 By \cite[Proposition~8.5.5]{bernstein} this implies that
 $\hat{P}^{\alpha'}\succ(H_2^\top\hat{P}_2^{-1}H_2)^{-1}$
 and hence that
 \begin{align*}
  \bar{J}(\hat{P}^{\alpha'})&=J(\hat{P}^{\alpha'}) \\
  &> J((H_2^\top\hat{P}_2^{-1}H_2)^{-1})=J(\hat{P}^\alpha)=\bar{J}(\hat{P}^\alpha)
 \end{align*}
 for $\alpha:=\balpha(\hat{x}_1,\hat{x}_2,\hat{P}_1,\hat{P}_2)$, a contradiction.
 It follows that $\alpha'=0$ in this case.

 Consider the case where $H_2^\top\hat{P}_2^{-1}H_2 = H_1^\top\hat{P}_1^{-1}H_1$ and note that $\alpha'\in[0,1]$ in this case.

 Consider the case where $H_2^\top\hat{P}_2^{-1}H_2 \prec H_1^\top\hat{P}_1^{-1}H_1$.
 It then follows from Lemma~\ref{lem:corner_cases} \textit{(2)} and Lemma~\ref{lem:sandwich} that $H_1^\top\hat{P}_1^{-1}H_1 \succ 0$ and hence that
 \begin{align*}
  \alpha' H_1^\top\hat{P}_1^{-1}H_1 + (1-\alpha')H_2^\top\hat{P}_2^{-1}H_2 \succ 0.
 \end{align*}
 We can hence define
 $\hat{P}^{\alpha'}:=(\alpha' H_1^\top\hat{P}_1^{-1}H_1 + (1-\alpha')H_2^\top\hat{P}_2^{-1}H_2)^{-1}\in\PPn$.
 Assume, to arrive at a contradiction, that $\alpha'<1$.
 It follows that
 \begin{align*}
  (\hat{P}^{\alpha'})^{-1}
  &= \alpha' H_1^\top\hat{P}_1^{-1}H_1 + (1-\alpha')H_2^\top\hat{P}_2^{-1}H_2 \\
  &= H_1^\top\hat{P}_1^{-1}H_1 \\
     &\qquad - (1-\alpha')(H_1^\top\hat{P}_1^{-1}H_1-H_2^\top\hat{P}_2^{-1}H_2) \\
  &\prec H_1^\top\hat{P}_1^{-1}H_1.
 \end{align*}
 By \cite[Proposition~8.5.5]{bernstein} this implies that
 $\hat{P}^{\alpha'}\succ(H_1^\top\hat{P}_1^{-1}H_1)^{-1}$
 and hence that
 \begin{align*}
  \bar{J}(\hat{P}^{\alpha'})&=J(\hat{P}^{\alpha'}) \\
  &> J((H_1^\top\hat{P}_1^{-1}H_1)^{-1})=J(\hat{P}^\alpha)=\bar{J}(\hat{P}^\alpha)
 \end{align*}
 for $\alpha:=\balpha(\hat{x}_1,\hat{x}_2,\hat{P}_1,\hat{P}_2)$, a contradiction.
 It follows that $\alpha'=1$ in this case.

 Consider the case where $H_2^\top\hat{P}_2^{-1}H_2 \not\succeq H_1^\top\hat{P}_1^{-1}H_1$ and
 $H_1^\top\hat{P}_1^{-1}H_1 \not\succeq H_2^\top\hat{P}_2^{-1}H_2$.
 Assume, to arrive at a contradiction, that $\alpha' H_1^\top\hat{P}_1^{-1}H_1 + (1-\alpha')H_2^\top\hat{P}_2^{-1}H_2$ is singular.
 It follows that
 \begin{align*}
  &\bar{J}((\alpha' H_1^\top\hat{P}_1^{-1}H_1 + (1-\alpha')H_2^\top\hat{P}_2^{-1}H_2)^{-1}) = +\infty \\
  &\quad > J(\hat{P}^\alpha)=\bar{J}(\hat{P}^\alpha)
 \end{align*}
 for $\alpha:=\balpha(\hat{x}_1,\hat{x}_2,\hat{P}_1,\hat{P}_2)$, a contradiction.
 It follows that $\alpha' H_1^\top\hat{P}_1^{-1}H_1 + (1-\alpha')H_2^\top\hat{P}_2^{-1}H_2$ is non-singular in this case.

 In summary, it follows that $\balpha'\in\mathcal{K}_H^{\balpha}$.

 Let $\balpha^\ast\in[0,1]$ be a solution of problem~\eqref{def:CI} then it follows that
 \begin{align*}
  &\bar{J}((\balpha^\ast H_1^\top\hat{P}_1^{-1}H_1 + (1-\balpha^\ast)H_2^\top\hat{P}_2^{-1}H_2)^{-1}) \leq \\
  &\quad \bar{J}(\bPhat^\CIa) = \bar{J}((\balpha H_1^\top\hat{P}_1^{-1}H_1 + (1-\balpha)H_2^\top\hat{P}_2^{-1}H_2)^{-1})
 \end{align*}
 for all $\balpha\in\mathcal{K}_H^{\balpha}\not=\emptyset$ and hence it follows from the preceding argument that $\balpha^\ast\in\mathcal{K}_H^{\balpha}$.
\qed

\textit{Proof of Corollary~\ref{cor:reduction_to_CI}.}
 Let $(\bK_1^\ast,\bK_2^\ast,\bPhat^\ast)\in\mathcal{F}_H$ be a solution of problem~\eqref{def:OF} then it follows from Theorem~\ref{th:opt_implies_Kahan} that $(\bK_1^\ast,\bK_2^\ast,\bPhat^\ast)\in\mathcal{K}_H$ and is $J$-optimal over $\mathcal{K}_H$.
 Hence there exists $\balpha^\ast\in\mathcal{K}_H^{\balpha}$ such that $(\bK_1^\ast,\bK_2^\ast,\bPhat^\ast)=(\bK_1^{\balpha^\ast},\bK_2^{\balpha^\ast},\bPhat^{\balpha^\ast})$.
 Assume, to arrive at a contradiction, that $\balpha^\ast$ is not a solution of problem~\eqref{def:CI} then there exists a function $\balpha'\in[0,1]$ such that
 $\bar{J}((\balpha' H_1^\top\hat{P}_1^{-1}H_1 + (1-\balpha')H_2^\top\hat{P}_2^{-1}H_2)^{-1})<\bar{J}(\bPhat^{\balpha^\ast})$.
 By Theorem~\ref{th:CI_implies_Kahan} this implies that $\balpha'\in\mathcal{K}_H^{\balpha}$.
 But then it follows from $J$-optimality of $(\bK_1^\ast,\bK_2^\ast,\bPhat^\ast)$ over $\mathcal{K}_H$ that
 \begin{align*}
  \bar{J}(\bPhat^{\balpha^\ast})&=J(\bPhat^{\balpha^\ast})=J(\bPhat^\ast)\\
  &\leq J(\bPhat^{\balpha'})=\bar{J}(\bPhat^{\balpha'}) \\
  &= \bar{J}((\balpha' H_1^\top\hat{P}_1^{-1}H_1 + (1-\balpha')H_2^\top\hat{P}_2^{-1}H_2)^{-1}),
 \end{align*}
 a contradiction.
 It follows that $\balpha^\ast$ is a solution of problem~\eqref{def:CI} and $(\bK_1^\ast,\bK_2^\ast,\bPhat^\ast)=(\bK_1^\CI,\bK_2^\CI,\bPhat^\CI)$.

 Conversely, let $\balpha^\ast\in[0,1]$ be a solution of problem~\eqref{def:CI} and let
 \begin{align*}
  (\bK_1^\ast,\bK_2^\ast,\bPhat^\ast):=(\bK_1^\CI,\bK_2^\CI,\bPhat^\CI).
 \end{align*}
 Since $\balpha^\ast$ is a solution of problem~\eqref{def:CI}, it follows in particular that
 $\bar{J}(\bPhat^\ast)\leq\bar{J}(\bPhat^\CIa)$ for all $(\bK_1^\CIa,\bK_2^\CIa,\bPhat^\CIa)\in\mathcal{K}_H$.
 It follows from Theorem~\ref{th:CI_implies_Kahan} that $\balpha^\ast\in\mathcal{K}_H^{\balpha}$ and hence that $(\bK_1^\ast,\bK_2^\ast,\bPhat^\ast)=(\bK_1^{\balpha^\ast},\bK_2^{\balpha^\ast},\bPhat^{\balpha^\ast})\in\mathcal{K}_H$.
 This implies that $J(\bPhat^\ast)=\bar{J}(\bPhat^\ast)\leq\bar{J}(\bPhat^\CIa)=J(\bPhat^\CIa)$ for all $(\bK_1^\CIa,\bK_2^\CIa,\bPhat^\CIa)\in\mathcal{K}_H$ and hence it follows that $(\bK_1^\ast,\bK_2^\ast,\bPhat^\ast)$ is $J$-optimal over $\mathcal{K}_H$.
 It follows from Theorem~\ref{th:Kahan_implies_opt} that $(\bK_1^\ast,\bK_2^\ast,\bPhat^\ast)$ is $J$-optimal and hence a solution of problem~\eqref{def:OF}.
\qed

\textit{Proof of Theorem~\ref{th:reduction_to_SDP}.}
 Let $\mathcal{S}_H$ denote the feasible set of problem~\eqref{def:SDP} and denote by
 \begin{align*}
  \mathcal{FS}_H := \{(\bK_1,\bK_2,\bPhat)\,|\,\exists_{\balpha\in\R}(\balpha,\bK_1,\bK_2,\bPhat)\in\mathcal{S}_H\}
 \end{align*}
 its projection onto the linear fusion rules.
 It follows from Definition~\ref{def:unbiased_fusion_rule}, Definition~\ref{def:optimal_fusion} and Theorem~\ref{th:conservative_3} that $\mathcal{FS}_H=\mathcal{F}_H$.

 Let the linear fusion rule $(\bK_1^\ast,\bK_2^\ast,\bPhat^\ast)\in\mathcal{F}_H$ be a solution of problem~\eqref{def:OF} then $(\bK_1^\ast,\bK_2^\ast,\bPhat^\ast)\in\mathcal{FS}_H$ and hence there exists $\balpha^\ast\in\R$ such that $(\balpha^\ast,\bK_1^\ast,\bK_2^\ast,\bPhat^\ast)\in\mathcal{S}_H$.
 Since $(\bK_1^\ast,\bK_2^\ast,\bPhat^\ast)$ is $J$-optimal, it follows that $(\balpha^\ast,\bK_1^\ast,\bK_2^\ast,\bPhat^\ast)$ is a solution of problem~\ref{def:SDP}.

 Conversely, let $(\balpha^\ast,\bK_1^\ast,\bK_2^\ast,\bPhat^\ast)\in\mathcal{S}_H$ be a solution of problem~\ref{def:SDP} then it follows that $(\bK_1^\ast,\bK_2^\ast,\bPhat^\ast)\in\mathcal{FS}_H=\mathcal{F}_H$ and that $(\bK_1^\ast,\bK_2^\ast,\bPhat^\ast)$ is $J$-optimal, the latter because $\mathcal{FS}_H=\mathcal{F}_H$.
 It hence follows that $(\bK_1^\ast,\bK_2^\ast,\bPhat^\ast)$ is a solution of problem~\ref{def:OF}.
\qed

\textit{Proof of Theorem~\ref{th:SDP_equiv_CI}.}
 Let $(\balpha^\ast,\bK_1^\ast,\bK_2^\ast,\bPhat^\ast)\in\DomSDP$ be a solution of problem \eqref{def:SDP}.
 It follows from Theorem~\ref{th:reduction_to_SDP} that $(\bK_1^\ast,\bK_2^\ast,\bPhat^\ast)$ is a solution of problem~\eqref{def:OF}.
 It follows from Corollary~\ref{cor:reduction_to_CI} and Theorem~\ref{th:CI_implies_Kahan} that there exists a solution $\balpha\in\mathcal{K}_H^{\balpha}$ of problem \eqref{def:CI} such that
 $(\bK_1^\ast,\bK_2^\ast,\bPhat^\ast)=(\bK_1^\CIa,\bK_2^\CIa,\bPhat^\CIa)$.
 By Theorem~\ref{th:alpha} this implies that $\balpha=\balpha^\ast$ and hence that $\balpha^\ast$ is a solution of problem~\eqref{def:CI} and $(\bK_1^\ast,\bK_2^\ast,\bPhat^\ast)=(\bK_1^\CI,\bK_2^\CI,\bPhat^\CI)$.

 Conversely, let $\balpha^\ast$ be a solution of problem~\eqref{def:CI}.
 By Theorem~\ref{th:CI_implies_Kahan} it follows that $\balpha^\ast\in\mathcal{K}_H^{\balpha}$.
 It follows from Corollary~\ref{cor:reduction_to_CI} that $(\bK_1^\CI,\bK_2^\CI,\bPhat^\CI)$ is a solution of problem~\eqref{def:OF} and hence it follows from Theorem~\ref{th:reduction_to_SDP} that there exists $\balpha\in\R$ such that $(\balpha,\bK_1^\CI,\bK_2^\CI,\bPhat^\CI)$ is a solution of problem~\eqref{def:SDP}.
 By Theorem~\ref{th:alpha} this implies that $\balpha^\ast=\balpha$ and hence that
 $(\balpha^\ast,\bK_1^\CI,\bK_2^\CI,\bPhat^\CI)=(\balpha,\bK_1^\CI,\bK_2^\CI,\bPhat^\CI)$ is a solution of problem~\eqref{def:SDP}.
\qed

\textit{Proof of Theorem~\ref{th:CI}.}
 For $S\in\mathcal{S}_n$, denote by $\lambda_{\min}(S)$ and $\lambda_{\max}(S)$ the smallest and largest eigenvalue, respectively.
 For $\alpha\in[0,1]$, define $\Sigma^\alpha:=\alpha H_1^\top\hat{P}_1^{-1}H_1+(1-\alpha)H_2^\top\hat{P}_2^{-1}H_2$.
 
 It follows from Lemma~\ref{lem:corner_cases} \textit{(3)} that $\Sigma^\alpha$ is non-singular for all $\alpha\in(0,1)$ and hence that $\bar{J}((\Sigma^{0.5})^{-1})=J((\Sigma^{0.5})^{-1})<+\infty$.
 
 By Lemma~\ref{lem:sandwich}, $H_1^\top\hat{P}_1^{-1}H_1\succeq 0$ and $H_2^\top\hat{P}_2^{-1}H_2\succeq 0$ and hence
 $\lambda_{\max}(H_1^\top\hat{P}_1^{-1}H_1)\geq 0$ and $\lambda_{\max}(H_2^\top\hat{P}_2^{-1}H_2)\geq 0$.
 
 By \cite[Proposition 8.5.13]{bernstein},
 \begin{align*}
  &\lambda_{\max}(\Sigma_\alpha)\\
  &\quad \leq \alpha\lambda_{\max}(H_1^\top\hat{P}_1^{-1}H_1) + (1-\alpha)\lambda_{\max}(H_2^\top\hat{P}_2^{-1}H_2)\\
  &\quad \leq \lambda_{\max}(H_1^\top\hat{P}_1^{-1}H_1) + \lambda_{\max}(H_2^\top\hat{P}_2^{-1}H_2) =: M
 \end{align*}
 for all $\alpha\in(0,1)$, and $\Sigma^{0.5}\succ 0$ implies that $M>0$.
 It follows that 
 \begin{align*}
  \lambda_{\min}((\Sigma^\alpha)^{-1})=\frac{1}{\lambda_{\max}(\Sigma^\alpha)} \geq \frac{1}{M}>0
 \end{align*}
 for all $\alpha\in(0,1)$.

 Case 1: $H_2^\top\hat{P}_2^{-1}H_2$ singular.
 In this case,
 \begin{align*}
  \lambda_{\min}((1-\alpha)H_2^\top\hat{P}_2^{-1}H_2)=0
 \end{align*}
 for all $\alpha\in[0,1]$.
 By \cite[Theorem~8.4.11]{bernstein} this implies that
 \begin{align*}
  &\lambda_{\min}(\Sigma^\alpha) \\
  &\quad \leq \lambda_{\max}(\alpha H_1^\top\hat{P}_1^{-1}H_1)+\lambda_{\min}((1-\alpha)H_2^\top\hat{P}_2^{-1}H_2)\\
  &\quad = \lambda_{\max}(\alpha H_1^\top\hat{P}_1^{-1}H_1)
 \end{align*}
 and hence that $\lim_{\alpha\to 0}\lambda_{\min}(\Sigma^\alpha)=0$ and
 \begin{align*}
  \lim_{\alpha\to 0}\lambda_{\max}((\Sigma^\alpha)^{-1})=\lim_{\alpha\to 0}\frac{1}{\lambda_{\min}(\Sigma^\alpha)}=+\infty.
 \end{align*}
 Since $J$ is spectrally unbounded on $\Pn$ it follows that $\lim_{\alpha\to 0}J((\Sigma^\alpha)^{-1})=+\infty$ and hence there exists $\alpha_l\in(0,0.5)$ such that
 \begin{equation}\label{eq:alpha_lower}
  \bar{J}((\Sigma^\alpha)^{-1})
  \begin{cases}
   > \bar{J}((\Sigma^{0.5})^{-1}), & \alpha\in[0,\alpha_l), \\
   < +\infty, & \alpha\in[\alpha_l,0.5],
  \end{cases}
 \end{equation}
 where the case $\alpha\in[\alpha_l,0.5]$ follows from Lemma~\ref{lem:corner_cases} \textit{(3)}.

 Case 2: $H_2^\top\hat{P}_2^{-1}H_2$ non-singular.
 In this case, define $\alpha_l:=0$ and observe that inequality~\eqref{eq:alpha_lower} still holds in this case, the case $\alpha\in[0,\alpha_l)$ trivially so because $[0,\alpha_l)=\emptyset$ and the case $\alpha\in[\alpha_l,0.5]$ by assumption (for $\alpha=0$) and Lemma~\ref{lem:corner_cases} \textit{(3)} (for $\alpha\in(0,0.5]$).

 In summary, there exists $\alpha_l\in[0,0.5)$ such that inequality~\eqref{eq:alpha_lower} holds.
 We proceed by distinguishing two cases for the matrix $H_1^\top\hat{P}_1^{-1}H_1\succeq 0$.

 Case 1: $H_1^\top\hat{P}_1^{-1}H_1$ singular.
 In this case,
 \begin{align*}
  \lambda_{\min}(\alpha H_1^\top\hat{P}_1^{-1}H_1)=0
 \end{align*}
 for all $\alpha\in[0,1]$.
 By \cite[Theorem~8.4.11]{bernstein} this implies that
 \begin{align*}
  &\lambda_{\min}(\Sigma^\alpha) \\
  &\quad \leq \lambda_{\min}(\alpha H_1^\top\hat{P}_1^{-1}H_1)+\lambda_{\max}((1-\alpha)H_2^\top\hat{P}_2^{-1}H_2)\\
  &\quad = \lambda_{\max}((1-\alpha)H_2^\top\hat{P}_2^{-1}H_2)
 \end{align*}
 and hence that $\lim_{\alpha\to 1}\lambda_{\min}(\Sigma^\alpha)=0$ and
 \begin{align*}
  \lim_{\alpha\to 1}\lambda_{\max}((\Sigma^\alpha)^{-1})=\lim_{\alpha\to 1}\frac{1}{\lambda_{\min}(\Sigma^\alpha)}=+\infty.
 \end{align*}
 Since $J$ is spectrally unbounded on $\Pn$ it follows that $\lim_{\alpha\to 1}J((\Sigma^\alpha)^{-1})=+\infty$ and hence there exists $\alpha_u\in(0.5,1)$ such that
 \begin{equation}\label{eq:alpha_upper}
  \bar{J}((\Sigma^\alpha)^{-1})
  \begin{cases}
   < +\infty, & \alpha\in[0.5,\alpha_u], \\
   > \bar{J}((\Sigma^{0.5})^{-1}), & \alpha\in(\alpha_u,1],
  \end{cases}
 \end{equation}
 where the case $\alpha\in[0.5,\alpha_u]$ follows from Lemma~\ref{lem:corner_cases} \textit{(3)}.

 Case 2: $H_1^\top\hat{P}_1^{-1}H_1$ non-singular.
 In this case, define $\alpha_u:=1$ and observe that inequality~\eqref{eq:alpha_upper} still holds in this case, the case $\alpha\in(\alpha_u,1]$ trivially so because $(\alpha_u,1]=\emptyset$ and the case $\alpha\in[0.5,\alpha_u]$ by assumption (for $\alpha=1$) and Lemma~\ref{lem:corner_cases} \textit{(3)} (for $\alpha\in[0.5,1)$).

 In summary, there exist $\alpha_l\in[0,0.5)$ and $\alpha_u\in(0.5,1]$ such that $\bar{J}((\Sigma^\alpha)^{-1})=J((\Sigma^\alpha)^{-1})<+\infty$ for all $\alpha\in[\alpha_l,\alpha_u]$ and $\bar{J}((\Sigma^\alpha)^{-1})>\bar{J}((\Sigma^{0.5})^{-1})$ for all $\alpha\in[0,1]\setminus[\alpha_l,\alpha_u]$.
 Since the function $\alpha\mapsto\bar{J}((\Sigma^\alpha)^{-1})=J((\Sigma^\alpha)^{-1})$ depends continuously on $\alpha$ on the compact interval $[\alpha_l,\alpha_u]$, it assumes a finite minimum for some $\alpha^\ast\in[a_l,a_u]\subseteq[0,1]$.
 Since $0.5\in[\alpha_l,\alpha_u]$, it follows that $\bar{J}((\Sigma^\alpha)^{-1})>\bar{J}((\Sigma^{0.5})^{-1})\geq\bar{J}((\Sigma^{\alpha^\ast})^{-1})$ also for all $\alpha\in[0,1]\setminus[\alpha_l,\alpha_u]$ and the result follows.
\qed

\textit{Proof of Theorem~\ref{th:CI_det}.}
 Let $(\hat{x}_1,\hat{x}_2,\hat{P}_1,\hat{P}_2)\in\Rpo\times\Rpt\times\PPpo\times\PPpt$, define
 \begin{align*}
  \Domalpha:=\{\alpha\in[0,1]\,|\,\Sigma_\alpha\text{ non-singular}\}
 \end{align*}
 and note that Lemma~\ref{lem:corner_cases} \textit{(3)} implies that $\Sigma_\alpha\succ 0$ for all $\alpha\in(0,1)$ and hence $(0,1)\subseteq\Domalpha$.
 Since $\Sigma_\alpha$ depends continuously on $\alpha$, it follows that $\Sigma_\alpha\succ 0$ for all $\alpha\in\Domalpha$.
 Define
 \begin{align*}
  g\colon\Domalpha\longrightarrow\R,\ \alpha\mapsto \widetilde{\det}((\Sigma_\alpha)^{-1})=\det((\Sigma_\alpha)^{-1})
 \end{align*}
 and compute $g'(\alpha)=-(\det(\Sigma_\alpha))^{-2}\Delta(\alpha)$ and
 $g''(\alpha)=(\det(\Sigma_\alpha))^{-3}((\Delta(\alpha))^2+\Trace((\adj(\Sigma_\alpha)(\Sigma_1-\Sigma_0))^2))$.
 Since $\Sigma_\alpha\succ 0$ implies that $\adj(\Sigma_\alpha)=\det(\Sigma_\alpha)(\Sigma_\alpha)^{-1}\succ 0$, it follows that
 \begin{align*}
  &\Trace((\adj(\Sigma_\alpha)(\Sigma_1-\Sigma_0))^2) = \\
  &\vectorize(\Sigma_1-\Sigma_0)^\top(\adj(\Sigma_\alpha)\otimes\adj(\Sigma_\alpha))\vectorize(\Sigma_1-\Sigma_0) > 0
 \end{align*}
 for all $\alpha\in\Domalpha$ and $\Sigma_0\not=\Sigma_1$, where the equality follows from \cite[Fact~7.4.8]{bernstein} and the inequality follows from \cite[Fact~7.4.15]{bernstein}.
 It follows that $g''(\alpha)>0$ for all $\alpha\in\Domalpha$ and $\Sigma_0\not=\Sigma_1$ and hence that $g$ is strictly convex on $\Domalpha$ as long as $\Sigma_0\not=\Sigma_1$.
 Using the extended calculus on $\mathcal{S}_n\cup\{+\infty\}$, it follows that the extended function
 \begin{align*}
  \tilde{g}\colon[0,1]\longrightarrow\R\cup\{+\infty\},\
  \alpha \mapsto \widetilde{\det}((\Sigma_\alpha)^{-1})
 \end{align*}
 is strictly convex as long as $\Sigma_0\not=\Sigma_1$ and hence problem~\eqref{def:CI} with $J=\det$ has at most one solution in the case $\Sigma_0\not=\Sigma_1$.

 By Theorem~\ref{th:CI}, problem~\eqref{def:CI} with $J=\det$ has a solution and as explained in the discussion preceding the statement of Theorem~\ref{th:CI_det}, any value for $\alpha^\ast\in[0,1]$ is a solution for the case $\Sigma_0=\Sigma_1$.

 If $\Delta(0)\leq 0$ and $\Sigma_0\not=\Sigma_1$, it follows that $g'(0)\geq 0$ and $g'(\alpha)>0$ for all $\alpha\in(0,1)$ and hence the unique minimum of $g$ must occur at $\alpha^\ast=0$.
 Similarly, if $\Delta(1)\geq 0$ and $\Sigma_0\not=\Sigma_1$, it follows that $g'(1)\leq 0$ and $g'(\alpha)<0$ for all $\alpha\in(0,1)$ and hence the unique minimum of $g$ must occur at $\alpha^\ast=1$.

 In the remaining case, either $g'(0)<0$ or $0\not\in\Domalpha$ and hence $\tilde{g}(0)=+\infty$, both of which imply that $\alpha^\ast>0$. Similarly, either $g'(1)>0$ or $1\not\in\Domalpha$ and hence $\tilde{g}(1)=+\infty$, both of which imply that $\alpha^\ast<1$.
 It follows that $\alpha^\ast$ is the unique solution of $g'(\alpha)=0\Leftrightarrow\Delta(\alpha)=0$ in $(0,1)$ in this case.
\qed

\textit{Proof of Theorem~\ref{th:CI_tr}.}
 Define the parameter set
 \begin{align*}
  \Domgamma := \{&\gamma=(A,B,c) \,|\, A\in\R^{n\times p_1},\ A\not=0, \\
                 &\quad B\in\R^{n\times p_2},\ B\not=0,\ c\in(0,1), \\
                 &\quad AH_1+BH_2=I
               \}
 \end{align*}
 and observe that $\Domgamma$ is equal to the intersection of an open set in a Euclidean space and an affine space defined by the constraint $AH_1+BH_2=I$.
 It follows that $\Domgamma$ is a valid domain for the Lagrange multiplier theorem.
 For every $\gamma\in\Domgamma$, define the linear fusion rule
 $(\bK_1^\gamma,\bK_2^\gamma,\bPhat^\gamma):=(A,B,\frac{1}{c}A\hat{P}_1A^\top+\frac{1}{1-c}B\hat{P}_2B^\top)$.
 It follows from Definition~\ref{def:unbiased_fusion_rule}, Theorem~\ref{th:conservative_3} \textit{(4)} applied with $\balpha:=c$ and Definition~\ref{def:optimal_fusion} that
 \begin{align*}
  \mathcal{F}_\gamma := \{(\bK_1^\gamma,\bK_2^\gamma,\bPhat^\gamma)\,|\,\gamma\in\Domgamma\}
  \subseteq \mathcal{F}_H,
 \end{align*}
 where we have used $\hat{P}_1\succ 0$ resp. $\hat{P}_2\succ 0$ to conclude that $\bQ_1=\bK_1^\gamma\hat{P}_1^{\frac{1}{2}}=A\hat{P}_1^{\frac{1}{2}}\not=0$ resp.
 $\bQ_2=\bK_2^\gamma\hat{P}_2^{\frac{1}{2}}=B\hat{P}_2^{\frac{1}{2}}\not=0$ for all $\gamma\in\Domgamma$.

 Let $\balpha^\ast\in[0,1]$ be a solution of problem~\eqref{def:CI} with $J=\Trace$.
 Fix $(\hat{x}_1,\hat{x}_2,\hat{P}_1,\hat{P}_2)\in\Rpo\times\Rpt\times\PPpo\times\PPpt$ and define
 \begin{align*}
  \ro &:= \sqrt{\Trace(K_1^{\alpha^\ast}\hat{P}_1(K_1^{\alpha^\ast})^\top)} \quad \text{and} \\
  \rt &:= \sqrt{\Trace(K_2^{\alpha^\ast}\hat{P}_2(K_2^{\alpha^\ast})^\top)}.
 \end{align*}
 We need to show that
 \begin{equation}\label{eq:CI_tr_proof}
  \alpha^\ast=\frac{\ro}{\ro+\rt}.
 \end{equation}
 The proof proceeds case-by-case.

 Case 1: $\alpha^\ast=0$.
 In this case it follows that $K_1^{\alpha^\ast}=\alpha^\ast\hat{P}^{\alpha^\ast}H_1^\top\hat{P}_1^{-1}=0$ and hence that $\ro=0$.
 It then follows from unbiasedness that $K_2^{\alpha^\ast}H_2=I$ and hence that $K_2^{\alpha^\ast}\not=0$ and thus $\rt\not=0$.
 But this implies that
 \begin{align*}
  \frac{\ro}{\ro+\rt} = \frac{0}{0+\rt} = 0 = \alpha^\ast
 \end{align*}
 in this case.

 Case 2: $\alpha^\ast=1$.
 In this case it follows that $K_2^{\alpha^\ast}=(1-\alpha^\ast)\hat{P}^{\alpha^\ast}H_2^\top\hat{P}_2^{-1}=0$ and hence that $\rt=0$.
 It then follows from unbiasedness that $K_1^{\alpha^\ast}H_1=I$ and hence that $K_1^{\alpha^\ast}\not=0$ and thus $\ro\not=0$.
 But this implies that
 \begin{align*}
  \frac{\ro}{\ro+\rt} = \frac{\ro}{\ro+0} = 1 = \alpha^\ast
 \end{align*}
 in this case.

 Case 3: $\alpha^\ast\in(0,1)$.
 In this case it follows that
 \begin{align*}
  \hat{P}^{\alpha^\ast}
  &= \hat{P}^{\alpha^\ast}(\hat{P}^{\alpha^\ast})^{-1}\hat{P}^{\alpha^\ast} \\
  &= \hat{P}^{\alpha^\ast}(\alpha^\ast H_1^\top\hat{P}_1^{-1}H_1+(1-\alpha^\ast)H_2^\top\hat{P}_2^{-1}H_2)\hat{P}^{\alpha^\ast} \\
  &= K_1^{\alpha^\ast}H_1\hat{P}^{\alpha^\ast} + K_2^{\alpha^\ast}H_2\hat{P}^{\alpha^\ast} \\
  &= \frac{1}{\alpha^\ast}K_1^{\alpha^\ast}\hat{P}_1(K_1^{\alpha^\ast})^\top
      + \frac{1}{1-\alpha^\ast}K_2^{\alpha^\ast}\hat{P}_2(K_2^{\alpha^\ast})^\top
 \end{align*}
 and hence that
 \begin{align*}
  (K_1^{\alpha^\ast},K_2^{\alpha^\ast},\hat{P}^{\alpha^\ast})
  \in (\mathcal{F}_\gamma)|_{\{(\hat{x}_1,\hat{x}_2,\hat{P}_1,\hat{P}_2)\}}.
 \end{align*}
 Since $\mathcal{F}_\gamma\subseteq\mathcal{F}_H$, it follows from Corollary~\ref{cor:reduction_to_CI} that $(K_1^{\alpha^\ast},K_2^{\alpha^\ast},\hat{P}^{\alpha^\ast})$ minimizes $\Trace(\hat{P}^\gamma)$ over this set.
 Defining the Lagrangian
 \begin{align*}
  L(\gamma,\lambda) := &\Trace(\frac{1}{c}A\hat{P}_1A^\top+\frac{1}{1-c}B\hat{P}_2B^\top) \\
                        &\qquad + \lambda^\top(AH_1+BH_2-I)
 \end{align*}
 it follows that there exists $\lambda^\ast\in\R^n$ such that
 \begin{align*}
  0 &= \frac{\partial L}{\partial c}((K_1^{\alpha^\ast},K_2^{\alpha^\ast},\alpha^\ast),\lambda^\ast) \\
    &= -\frac{1}{(\alpha^\ast)^2}(\ro)^2 + \frac{1}{(1-\alpha^\ast)^2}(\rt)^2
 \end{align*}
 or, equivalently,
 \begin{align*}
  ((\rt)^2-(\ro)^2)(\alpha^\ast)^2+2(\ro)^2\alpha^\ast-(\ro)^2 = 0.
 \end{align*}
 This quadratic equation has a unique solution $\alpha^\ast\in(0,1)$ given by equation~\eqref{eq:CI_tr_proof}.
\qed

\bibliographystyle{plain}
\bibliography{refs}
\end{document}